\documentclass[12pt]{amsart}

\usepackage{amssymb, amsmath, amsthm}
\usepackage{calligra, mathrsfs}
\usepackage{graphicx}
\usepackage{url}
\usepackage{mathtools}
\usepackage{enumerate}
\usepackage{verbatim}
\usepackage[retainorgcmds]{IEEEtrantools}
\usepackage{tikz-cd}
\usetikzlibrary{positioning}

\usepackage[margin=1in,marginparwidth=0.8in, marginparsep=0.1in]{geometry}

\newcommand{\Acal}{{\mathcal A}}
\newcommand{\Bcal}{{\mathcal B}}
\newcommand{\Ccal}{{\mathcal C}}
\newcommand{\Dcal}{{\mathcal D}}
\newcommand{\Ecal}{{\mathcal E}}
\newcommand{\Fcal}{{\mathcal F}}

\newcommand{\Ical}{{\mathcal I}}
\newcommand{\Jcal}{{\mathcal J}}
\newcommand{\Kcal}{{\mathcal K}}

\newcommand{\Mcal}{{\mathcal M}}
\newcommand{\Ncal}{{\mathcal N}}
\newcommand{\Ocal}{{\mathcal O}}
\newcommand{\Pcal}{{\mathcal P}}
\newcommand{\Qcal}{{\mathcal Q}}

\newcommand{\Scal}{{\mathcal S}}

\newcommand{\ZZ}{{\mathbb Z}}
\newcommand{\NN}{{\mathbb N}}

\newcommand{\Obb}{{\mathbb O}}

\newcommand{\Cfrak}{\mathfrak{C}}
\newcommand{\ffrak}{{\mathfrak f}}

\newcommand{\vfrak}{{\mathfrak v}}

\usepackage{scalerel,stackengine}
\stackMath
\newcommand\reallywidehat[1]{%
\savestack{\tmpbox}{\stretchto{%
  \scaleto{%
    \scalerel*[\widthof{\ensuremath{#1}}]{\kern-.6pt\bigwedge\kern-.6pt}%
    {\rule[-\textheight/2]{1ex}{\textheight}}
  }{\textheight}%
}{0.5ex}}%
\stackon[1pt]{#1}{\tmpbox}%
}

\makeatletter
\newcommand*\bigcdot{{\mathpalette\bigcdot@{.5}}}
\newcommand*\bigcdot@[2]{\mathbin{\vcenter{\hbox{\scalebox{#2}{$\m@th#1\bullet$}}}}}
\makeatother

\newcommand{\dsh}{{\text{\normalfont-}}}

\makeatletter
\g@addto@macro\bfseries{\boldmath}
\makeatother

\tikzset{shorten <>/.style={shorten >=#1,shorten <=#1}}

\newcommand{\kr}{\kern -2pt}

\DeclareMathOperator{\id}{id}
\DeclareMathOperator{\Hom}{Hom}

\DeclareMathOperator{\End}{End}
\DeclareMathOperator{\Funct}{Funct}
\DeclareMathOperator{\ev}{ev}

\newcommand{\cl}{{\text{\normalfont cl}}}

\newcommand{\enh}{{\text{\normalfont enh}}}

\newcommand{\op}{{\text{\normalfont op}}}
\newcommand{\rev}{{\text{\normalfont rev}}}

\DeclareMathOperator{\Spc}{Spc}

\DeclareMathOperator{\Set}{Set}

\DeclareMathOperator{\Cat}{Cat}
\newcommand{\twoCat}{{2\kr\Cat}}
\newcommand{\nCat}{{n\kr\Cat}}
\DeclareMathOperator{\omegaCat}{\omega\kr\Cat}

\DeclareMathOperator{\colim}{colim}

\newcommand{\rex}{{\text{\normalfont rex}}}

\newcommand{\cart}{{\text{\normalfont cart}}}
\newcommand{\cocart}{{\text{\normalfont cocart}}}

\DeclareMathOperator{\Adj}{Adj}

\newcommand{\radj}{{\text{\normalfont radj}}}
\newcommand{\ladj}{{\text{\normalfont ladj}}}

\newcommand{\cocompl}{{\text{\normalfont cocompl}}}
\DeclareMathOperator{\Catcocompl}{\widehat{\Cat}_\cocompl}
\DeclareMathOperator{\Ind}{Ind}
\let\Pr\relax
\DeclareMathOperator{\Pr}{Pr}
\newcommand{\pr}{{\text{\normalfont pr}}}

\DeclareMathOperator{\Sp}{Sp}
\newcommand{\St}{{\text{\normalfont St}}}

\DeclareMathOperator{\Op}{Op}
\DeclareMathOperator{\Fin}{Fin}
\DeclareMathOperator{\Assos}{Assos}
\DeclareMathOperator{\BM}{BM}
\DeclareMathOperator{\LM}{LM}
\DeclareMathOperator{\RM}{RM}

\DeclareMathOperator{\Env}{Env}

\DeclareMathOperator{\Alg}{{Alg}}
\DeclareMathOperator{\CAlg}{CAlg}

\DeclareMathOperator{\Algbrd}{Algbrd}

\DeclareMathOperator{\modd}{{-mod}}

\DeclareMathOperator{\LMod}{LMod}
\DeclareMathOperator{\BMod}{BMod}

\newcommand*{\shom}{{\mathscr{H}\kern -2pt om}}
\newcommand*{\Catscr}{{\mathscr{C}\kern -2pt at}}
\newcommand*{\Prscr}{{\mathscr{P}\kern -2pt r}}
\newcommand{\shext}{{\mathscr{E} \kern -2pt xt}}

\newtheorem{definition}[subsubsection]{Definition}
\newtheorem{proposition}[subsubsection]{Proposition}
\newtheorem{lemma}[subsubsection]{Lemma}
\newtheorem{theorem}[subsubsection]{Theorem}
\newtheorem{corollary}[subsubsection]{Corollary}
\newtheorem{conjecture}[subsubsection]{Conjecture}
\newtheorem*{theorem*}{Theorem}
\newtheorem*{conjecture*}{Conjecture}

\newtheoremstyle{claim}{\topsep}{\topsep}{}{}{\it}{.}{.5em}{}
\theoremstyle{claim}
\newtheorem*{claim*}{Claim}

\newtheoremstyle{note}
{8.0pt plus 2.0pt minus 4.0pt}{8.0pt plus 2.0pt minus 4.0pt}{}{}{\bf}{.}{.5em}{} 
\theoremstyle{note}
\newtheorem{example}[subsubsection]{Example}
\newtheorem{remark}[subsubsection]{Remark}
\newtheorem{notation}[subsubsection]{Notation}
\newtheorem{construction}[subsubsection]{Construction}
\newtheorem{warning}[subsubsection]{Warning}
\newtheorem*{example*}{Example}

\AtBeginDocument{\def\MR#1{}}

\usepackage{standalone} 

\title{Presentable $(\infty,n)$-categories}

\author{G. Stefanich}

\date{}

\begin{document}

\begin{abstract}
We define for each $n \geq 1$ a symmetric monoidal $(\infty,n+1)$-category $n\Prscr^L$ whose objects we call presentable $(\infty,n)$-categories, generalizing the usual theory of presentable $(\infty,1)$-categories. We show that each object $\Ccal$  in $n\Prscr^L$ has an underlying $(\infty,n)$-category $\psi_n(\Ccal)$ which admits all conical colimits, and that conical colimits of right adjointable diagrams in $\psi_n(\Ccal)$ can be computed in terms of conical limits after passage to right adjoints.
\end{abstract}

\maketitle

\tableofcontents


\section{Introduction}

The theory of presentable $(\infty,1)$-categories is one of the cornerstones of higher category theory. Many higher categories arising in nature are presentable, and many categorical constructions preserve presentability. The $(\infty,1)$-category $\Pr^L$ of presentable $(\infty,1)$-categories and colimit preserving functors can be given a symmetric monoidal structure, and is a convenient setting in which one can do algebra with $(\infty,1)$-categories.

This paper introduces an $(\infty, n)$-categorification of the theory of presentable $(\infty,1)$-categories. Our goal is to provide a context, analogous to $\Pr^L$, in which one can do algebra with $(\infty, n)$-categories. Our starting point is the observation that $\Pr^L$ itself is not a presentable $(\infty,1)$-category. We argue that $\Pr^L$ should instead be considered as a presentable $(\infty,2)$-category - we in fact define this notion so that $\Pr^L$ is the unit in the symmetric monoidal category of presentable $(\infty,2)$-categories. In general, we  find that for each $n \geq 1$ the collection of all presentable $(\infty,n)$-categories can be organized, not into a presentable $(\infty,n)$-category, but rather into a presentable $(\infty,n+1)$-category.

This article is a companion to an upcoming series of papers on the foundations of the theory of sheaves of $(\infty,n)$-categories in derived algebraic geometry. Non-categorical sheaf theories usually assign to each geometric object $X$ a stable presentable $(\infty,1)$-category of sheaves on $X$, subject to various functoriality properties. The $(\infty, n+1)$-category $n\Prscr^L$ (and in particular its full subcategory of stable objects) will be our target for various theories of sheaves of $(\infty,n)$-categories.

 Below we provide a more detailed description of the contents of this paper. We will use the convention where all objects are $\infty$-categorical by default, and suppress this from our notation from now on.

\subsection{Higher presentable categories} We begin by motivating the definition of the symmetric monoidal category $n\kr\Pr^L$ underlying $n\kr\Prscr^L$. Recall the following features of $\Pr^L$:
\begin{itemize}
\item The category $\Pr^L$ admits all small colimits.
\item For each pair of objects $\Ccal, \Dcal$ in $\Pr^L$, there is a presentable category $\Funct(\Ccal, \Dcal)$ parametrizing functors between $\Ccal$ and $\Dcal$. This endows $\Pr^L$ with the structure of a category enriched over itself.
\end{itemize}

Any presentable $2$-category will share the same features. If $\Pr^L$ were itself presentable, then any module over $\Pr^L$ inside $\Pr^L$ would have these properties. In that case, we could simply define $2\kr\Pr^L$ to be the category of modules for $\Pr^L$ inside itself.

Since $\Pr^L$ is not presentable, it does not make sense to consider its category of modules in $\Pr^L$. However, $\Pr^L$ is still a cocomplete category, and its symmetric monoidal structure is compatible with colimits. It therefore makes sense to consider the  category of modules for $\Pr^L$ inside the symmetric monoidal category $\widehat{\Cat}_{\cocompl}$ of large cocomplete categories admitting small colimits and colimit preserving functors.

The category $\Pr^L\modd(\widehat{\Cat}_{\cocompl})$ is however still not a good candidate for $2\kr\Pr^L$. A first issue is that its objects turn out to be enriched in $\widehat{\Cat}_{\cocompl}$, rather than $\Pr^L$. A related problem is that $\Pr^L\modd(\widehat{\Cat}_{\cocompl})$ is too big. Recall that presentable categories are controlled by a small amount of data, which in particular implies that $\Pr^L$ is a large category (as opposed to very large). However $\Pr^L\modd(\widehat{\Cat}_{\cocompl})$ is a very large category, and in particular it does not belong to $\widehat{\Cat}_{\cocompl}$, but to the category $\widehat{\operatorname{CAT}}_{\cocompl}$ of very large categories admitting large colimits. A theory of presentable $n$-categories which simply iteratively takes categories of modules would need a long sequence of nested universes, and would require one to carefully keep track of the relative sizes of different objects.

The main observation that leads to our definition of higher presentable categories is that, while $\Pr^L\modd(\widehat{\Cat}_\cocompl)$ is a very large category, it is controlled by a large subcategory of objects we call presentable. More precisely, $\Pr^L\modd(\widehat{\Cat}_\cocompl)$ is presentable as a very large category, and it is in fact $\kappa_0$-compactly generated, for $\kappa_0$ the smallest large cardinal. The category $2\kr\Pr^L$ is then defined to be the full subcategory of $\Pr^L\modd(\widehat{\Cat}_\cocompl)$ on the $\kappa_0$-compact objects.

Iterating the above discussion, one may obtain a symmetric monoidal category $n\kr\Pr^L$ for each $n \geq 1$:
\begin{definition}
We inductively define, for each $n \geq 2$, the symmetric monoidal category $n\kr\Pr^L$ to be the full subcategory of $(n-1)\kr\Pr^L\modd(\Cat_{\normalfont \cocompl})$ on the $\kappa_0$-compact objects.
\end{definition}

Our first main result is the construction of an $(n+1)$-categorical enhancement of $n\kr\Pr^L$, which we denote by $n\Prscr^L$. This is obtained as a consequence of the existence of a lax symmetric monoidal realization functor which maps objects in $(n-1)\kr\Pr^L\modd(\Cat_{\cocompl})$ to $n$-categories. In other words, although our approach to presentable $n$-categories starts out being $1$-categorical, we are in fact able to consider these as higher categories. To formulate the statement, it is convenient to set $0\kr\Pr^L$ to be the category of spaces.

\begin{theorem}\label{teo existence realization}
Let $n \geq 1$. There exists a lax symmetric monoidal functor 
\[
\psi_n: (n-1)\kr\Pr^L\modd(\widehat{\Cat}_{\normalfont \text{cocompl}}) \rightarrow \reallywidehat{\nCat}
\]
with the following properties:
\begin{enumerate}[\normalfont (i)]
\item For each $(n-1)\kr\Pr^L$-module $\Ccal$, the $1$-category underlying $\psi_n(\Ccal)$ is equivalent to the category underlying $\Ccal$.
\item Assume $n > 1$ and let $c, d$ be objects in $\Ccal$. Then
\[
\Hom_{\psi_n(\Ccal)}(c, d) = \psi_{n-1}(\shom_\Ccal(c, d) )
\]
where $\shom_\Ccal(c, d)$ denotes the Hom object between $c$ and $d$ obtained from the action of $(n-1)\kr\Pr^L$ on $\Ccal$.
\end{enumerate}
\end{theorem}

\begin{remark}
Although our main interest is in objects of $n\kr\Pr^L$, we need to have access to the realization functor on the bigger category $(n-1)\kr\Pr^L\modd(\widehat{\Cat}_{\normalfont \text{cocompl}})$. The main reason is that if $n > 2$ we do not know whether the Hom objects in part (ii) belong to $(n-1)\kr\Pr^L$, even if $\Ccal$ itself belongs to $n\kr\Pr^L$.
\end{remark}

\begin{corollary}
Let $n \geq 1$. There exists a symmetric monoidal $(n+1)$-category $n\Prscr^L  = \psi_{n+1}(n\kr\Pr^L)$ whose underlying category is $n\kr\Pr^L$, and such that for every pair of objects $\Ccal, \Dcal$, we have
\[
\Hom_{n\Prscr^L}(\Ccal, \Dcal) = \psi_n (\shom_{n\kr\Pr^L}(\Ccal, \Dcal)).
\]
In particular, there is an equivalence
\[
\End_{n\Prscr^L}((n-1)\kr\Pr^L) =  (n-1)\Prscr^L.
\]
\end{corollary}

\subsection{Conical colimits}
A direct consequence of item (i)  in theorem \ref{teo existence realization} is that for every object $\Ccal$ in $(n-1)\kr\Pr^L\modd(\Catcocompl)$, the category underlying $\psi_n(\Ccal)$ admits all small colimits. Our second main result states that these are compatible with the $n$-categorical structure.

\begin{theorem}\label{teo hay conical limits}
Let $n \geq 1$ and let $\Ccal$ be an object  in $(n-1)\kr\Pr^L\modd(\widehat{\Cat}_{\normalfont \text{cocompl}})$. Then for every small category $\Ical$ the diagonal functor
\[
\psi_n(\Ccal) \rightarrow \Funct(\Ical, \psi_n(\Ccal))
\]
admits a left adjoint.
\end{theorem}

Note that in the case $n = 1$, the existence of this adjoint is equivalent to the existence of colimits for diagrams in $\Ccal$. When $n > 1$ the condition in theorem \ref{teo hay conical limits} is stronger than the mere existence of colimits in the category underlying $\psi_n(\Ccal)$.  

We can reinterpret theorem \ref{teo hay conical limits} as saying that $\psi_n(\Ccal)$ admits all conical colimits. Here the word conical is used to distinguish these from the more broad class of weighted colimits. We expect in fact that $\psi_n(\Ccal)$ will admit all weighted colimits - although a general theory of weighted colimits does not seem to be present in the literature yet, one can expect that the existence of conical colimits and tensors in $\psi_n(\Ccal)$ should imply this.

The following result is a typical application of theorem \ref{teo hay conical limits}.
\begin{corollary}\label{coro colim adjointability}
Let $\Ccal$ be an object in $(n-1)\kr\Pr^L\modd(\Catcocompl)$, and let $\Ical$ be a small category. Let $F, G : \Ical \rightarrow \Ccal$ be functors, and $\eta: F \rightarrow G$ be a natural transformation. Assume that for every arrow $\alpha: i \rightarrow j$ in $\Ical$ the commutative square
\[
\begin{tikzcd}
F(i) \arrow{d}{F(\alpha)} \arrow{r}{\eta_i} & G(i) \arrow{d}{G(\alpha)} \\ F(j) \arrow{r}{\eta_j} & G(j)
\end{tikzcd}
\]
is horizontally right adjointable in $\psi_n(\Ccal)$ (in other words, the horizontal arrows admit right adjoints $\eta_i^R$ and $\eta_j^R$ in $\psi_n(\Ccal)$, and the induced $2$-cell $F(\alpha) \eta_i^R \rightarrow \eta_j^R G(\alpha)$ is an isomorphism). Then the morphism
\[
\colim_\Ical \eta : \colim_{\Ical} F \rightarrow \colim_{\Ical} G
\]
admits a right adjoint in $\psi_n(\Ccal)$.
\end{corollary}

In other words, a colimit of right adjointable arrows in $\psi_n(\Ccal)$ is right adjointable, as long as various base change properties hold. In future applications we will be interested in other consequences of theorem \ref{teo hay conical limits} which are similar to corollary \ref{coro colim adjointability}, in which we replace ``right adjointable arrow'' by more complicated diagrams satisfying various adjointability conditions.

\subsection{The passage to adjoints property}

A feature of the theory of presentable categories is that passage to adjoints interchanges colimits and limits. Concretely, if $F: \Ical \rightarrow \Pr^L$ is a diagram, then the colimit of $F$ is equivalent to the limit of the diagram $\Ical^\op \rightarrow \widehat{\Cat}$ obtained from $F$ by passing to right adjoints of every arrow. 

In the case when the right adjoint to $F(\alpha)$ belongs to $\Pr^L$ for every arrow $\alpha$ in $\Ical$, this is also the same as the limit of the resulting functor $F^R: \Ical^\op \rightarrow \Pr^L$.  The following definition abstracts the basic features of  the  $2$-categorical enhancement $2\Prscr^L$ which makes this fact  hold.

\begin{notation}
Let $\Dcal$ be an $n$-category. We denote by $\Dcal^{\leq 1}$ the $1$-category underlying $\Dcal$, and by $(\Dcal^{\leq 1})^{\text{radj}}$ (resp. $\Dcal^{\leq 1})^{\text{ladj}}$ the subcategory of $\Dcal^{\leq 1}$ containing all objects, and only those morphisms which are right (resp. left) adjointable in $\Dcal$.
\end{notation}

\begin{definition}\label{def passage right adj prop}
Let $\Dcal$ be an $n$-category. We say that $\Dcal$ satisfies the passage to adjoints property if the following conditions are satisfied:
\begin{itemize}
\item The category $(\Dcal^{\leq 1})^{\text{\normalfont radj}}$ has all colimits, and the inclusion $(\Dcal^{\leq 1})^{\text{\normalfont radj}} \rightarrow \Dcal$ preserves conical colimits.
\item The category $(\Dcal^{\leq 1})^{\text{\normalfont ladj}}$ has all limits, and the inclusion $(\Dcal^{\leq 1})^{\text{ \normalfont ladj}} \rightarrow \Dcal$ preserves  conical limits.
\end{itemize}
\end{definition}

The third main result of this paper is that the $n$-categories obtained from the realization functors of theorem \ref{teo existence realization} satisfy the passage to adjoints property.

\begin{theorem}\label{teo passage adj}
Let $n \geq 2$ and let $\Ccal$ be an object in $(n-1)\kr\Pr^L\modd(\widehat{\Cat}_{\normalfont \cocompl})$. Then $\psi_n(\Ccal)$ satisfies the passage to adjoints property.
\end{theorem}

In particular, although we do not know whether $n\kr\Pr^L$ admits all small limits when $n > 1$,  we are able to conclude that it has limits of left adjointable diagrams.

\subsection{The role of enriched category theory}
This paper and its planned applications make heavy use of the theory of enriched categories, as developed in \cite{GH} and \cite{Hinich}. Indeed, the realization functor from theorem \ref{teo existence realization} is ultimately obtained by iteratively applying the procedure of enrichment of a presentable module of a presentable category.  We have therefore chosen  to include a systematic treatment of the theory of enriched categories, adapted to our purposes.

We make along the way various contributions to the state of the art. First, we provide a new approach to the definition of the operad $\Assos_X$ that corepresents enrichment. We show that the assignment $X \mapsto \Assos_X$ is corepresented by an associative cooperad $\Cfrak$ internal to  $\Cat$, which is in turn determined by its underlying associative cooperad $\Cfrak^\cl$ internal to $0$-truncated categories. The cooperad $\Cfrak^\cl$ is classical and it  can therefore  be defined by specifying a finite amount of data. It is in fact uniquely characterized by a surprisingly small amount of data: namely the $0$-truncated categories of objects and operations, and the cosource and cotarget maps. This observation allows us to compare our definition to previous approaches in the literature. We also use similar arguments to provide concise definitions of the operads $\LM_X$ and $\BM_X$ that play a key role in \cite{Hinich}.

Our second contribution concerns the functoriality of the procedure of enrichment of presentable modules over presentable categories.\footnote{While finishing this article we learned about recent work of H. Heine \cite{Heinenr} which yields a similar functorial strengthening of this procedure.} As shown in \cite{GH} and \cite{Hinich}, to each presentable module $\Ccal$ over a presentable symmetric monoidal category $\Mcal$ there is associated an $\Mcal$ enriched category $\overline{\Ccal}$. In this text we upgrade the assignment $\Ccal \mapsto \overline{\Ccal}$ to a lax symmetric monoidal functor 
\[
\theta_\Mcal : \Mcal\modd(\Pr^L) \rightarrow \widehat{\Cat}^\Mcal
\]
 and show that the composition of $\theta_\Mcal$ with the forgetful functor $\widehat{\Cat}^\Mcal \rightarrow \widehat{\Cat}$ recovers the usual lax symmetric monoidal forgetful functor $\Mcal\modd(\Pr^L) \rightarrow \widehat{\Cat}$. This functoriality is key in carrying out the construction of theorem \ref{teo existence realization}.

Third, we introduce the notion of conical colimits in enriched categories, and prove general results regarding the existence of conical colimits in limits of enriched categories, and in functor enriched categories. We show that for every presentable module $\Ccal$ over a presentable symmetric monoidal category $\Mcal$, the enriched category $\theta_\Mcal(\Ccal)$ admits all conical limits and colimits - this is one of the inputs that goes into the proof of theorem \ref{teo hay conical limits}.  To prove these results on conical limits and colimits we develop the theory of adjunctions between enriched categories. We give local and global versions of the notion of adjunction between enriched categories, and show that a functor of enriched categories admits a globally defined adjoint if and only if it admits local adjoints at every object.

\subsection{Organization} We now describe the contents of the paper in more detail. We refer the reader also to the introduction of each section for an expanded outline of its contents.

Section 2 is a general introduction to enriched category theory. We begin by describing our approach to the definition of the operad $\Assos_X$, and the closely related operads $\LM_X$ and $\BM_{X, Y}$. We then review the notions of algebroids and enriched categories, and the general functoriality and multiplicativity properties of the theory. We finish this section with an introduction to the theory of $n$-categories and $\omega$-categories via iterated enrichment.

Section 3  deals with the theory of left modules and bimodules over algebroids. The bulk of this section is devoted to the construction of a functorial enhancement of the procedure of enrichment of modules over presentable categories.

In section 4 we study the theory of adjunctions between enriched categories, and conical limits and colimits in enriched categories. We introduce the notion of local right adjoint to a functor between enriched categories, and show that a right adjoint exists if and only if all local adjoints exist. We prove a general result regarding the existence of local adjoints to limits of functors, and study the case of local adjunctions between enriched categories of functors. As a special case, we obtain general results on the existence of conical limits and colimits in limits of enriched categories, and in enriched categories of functors. Finally, we study how adjunctions interact with the procedure of enrichment of presentable modules, and show that the enriched categories that arise this way admit all conical limits and colimits.

In section 5 we construct the symmetric monoidal categories $n\kr\Pr^L$ and the realization functor from theorem \ref{teo existence realization}. We then give a proof of theorem \ref{teo hay conical limits}. Finally, we study the procedure of passage to adjoints and prove theorem \ref{teo passage adj}.

\subsection{Conventions and notation}

We use the language of higher category theory and higher algebra as developed in \cite{HTT} and \cite{HA}. All of our notions will be assumed to be homotopical or $\infty$-categorical, and we suppress this from our notation - for instance, we use the word $n$-category to mean $(\infty,n)$-category. 

We work with a nested sequence of universes. Objects belonging to the first universe are called small, objects in the second universe are called large, and objects in the third universe are called very large. 

We denote by $\Spc$ and $\Cat$ the categories of (small) spaces and categories. For each $n \geq 2$ we denote by $\nCat$ the category of small $n$-categories.  We denote by $\Catscr$ the $2$-categorical enhancement of $\Cat$, and in general by $n\Catscr$ the $(n+1)$-categorical enhancement of $\nCat$. If $X$ is one of those objects (or related), we denote by $\widehat{X}$ its large variant. For instance, $\widehat{\Spc}$ denotes the category of large spaces.

We denote by $\Pr^L$ the category of presentable categories and colimit preserving functors. We usually consider this as a symmetric monoidal category, with the tensor product constructed in \cite{HA} section 4.8. By presentable (symmetric) monoidal category we mean a (commutative) algebra in $\Pr^L$. In other words, this is a presentable category equipped with a (symmetric) monoidal structure compatible with colimits.

 If $\Ccal$ is an $n$-category and $k \geq 0$, we denote by $\Ccal^{\leq k}$ the $k$-category obtained from $\Ccal$ by discarding all cells of dimension greater than $k$, and by $\prescript{\leq k}{}{\Ccal}$ the $k$-category obtained from $\Ccal$ by inverting all cells of dimension greater than $k$. In particular, for each category $\Ccal$ we denote by $\Ccal^{\leq 0}$ the space of objects of $\Ccal$.
 
For each category $\Ccal$ we denote by $\Hom_{\Ccal}(-, -)$ the Hom-functor of $\Ccal$. We denote by $\Pcal(\Ccal)$ its presheaf category. We usually identify $\Ccal$ with its image under the Yoneda embedding.

Given a pair of $n$-categories $\Ccal , \Dcal$ we denote by $\Funct(\Ccal, \Dcal)$ the $n$-category of functors between $\Ccal$ and $\Dcal$. When we wish to only consider the space of functors between them we will use the notation $\Hom_{\nCat}(\Ccal, \Dcal)$ instead.

 Given a right (resp. left) adjointable functor of categories $\beta: \Ccal \rightarrow \Dcal$, we will usually denote by $\beta^R$ (resp. $\beta^L$) its right (resp. left) adjoint. More generally, we use a similar notation for adjoint morphisms in a $2$-category. We say that a commutative square of categories
 \[
 \begin{tikzcd}
 \Ccal' \arrow{d}{\beta'} \arrow{r}{\alpha'} & \Ccal \arrow{d}{\beta} \\
 \Dcal' \arrow{r}{\alpha} & \Dcal
 \end{tikzcd}
 \]
 is vertically right adjointable if $\beta$ and $\beta'$ admit right adjoints and the induced natural transformation $\alpha' \beta'^R \rightarrow \beta^R \alpha$ is an isomorphism. Similarly, we can talk about horizontal left adjointability, or vertical left / right adjointability. We will at various points use the connections between adjointability of squares and the theory of two-sided fibrations, as explained in \cite{HSTI} section 2.

We make use at various points of the theory of operads. We use a language for speaking about these which is close in spirit to the classical language in terms of objects and operations which satisfy a composition rule. Namely, given an operad $\Ocal$ with associated category of operators $p: \Ocal^\otimes \rightarrow \Fin_*$, we call $p^{-1}(\langle 1 \rangle)$ the category of objects of $\Ocal$, and  arrows in $\Ocal^\otimes$ lying above an active arrow of the form $\langle n \rangle \rightarrow \langle 1 \rangle$ are called operations of $\Ocal$. We will for the most part work with $\Ocal$ without making explicit reference to the fibration $p$, and make it clear when we need to refer to the category of operators instead.

 We denote by $\Op$ the category of operads, and for each operad $\Ocal$ we denote by $\Op_\Ocal$ the category of $\Ocal$-operads. We denote by $\Assos, \LM, \BM$ the operads for associative algebras, left modules, and bimodules, respectively.

\subsection{Acknowledgments}

I would like to thank David Nadler and Nick Rozenblyum for useful conversations related to the subject of this paper.



\tableofcontents

\section{Enriched category theory}

Let $\Mcal$ be a monoidal category. An algebroid $\Acal$ in $\Mcal$ with space of objects $X$ consists of
\begin{itemize}
\item For every pair of objects $x, y$ an object $\Acal(y, x)$ in $\Mcal$.
\item For every object $x$ in $X$ a morphism $1_\Mcal \rightarrow \Acal(x, x)$.
\item For every triple of objects $x, y, z$ a morphism $\Acal(z, y) \otimes \Acal(y, x) \rightarrow \Acal(z, x)$.
\item Associativity and unit isomorphisms, and an infinite list of higher coherence data.
\end{itemize}

Given an algebroid $\Acal$ with space of objects $X$, there is a  Segal space underlying $\Acal$, with space of objects $X$ and for each pair of objects $y, x$ the space of morphism being given by $\Hom_\Mcal(1_\Mcal, \Acal(y, x))$. We say that $\Acal$ is an $\Mcal$-enriched category if its underlying Segal space is complete. Our goal in this section is to review the theory of algebroids and enriched categories, and the approach to $n$-category theory via iterated enrichment.

 For each space $X$ the assignment $\Mcal \mapsto \Algbrd_X(\Mcal)$ that sends each monoidal category $\Mcal$ to the category of algebroids with space of objects $X$ turns out to be corepresented by a nonsymmetric operad $\Assos_X$, to be thought of as a many object version of the associative operad. The assignment $X \mapsto \Assos_X$ determines a functor from spaces to the category of associative operads, which is in turn corepresented by an associative cooperad $\Cfrak$ internal to the category $\Cat$.
 
  We begin in \ref{subsection cfrak} by reviewing the notion of internal operads and cooperads, and presenting the definition of the cooperad $\Cfrak$. This is defined starting from a  cooperad $\Cfrak^\cl$ in the (classical) category of posets, which can be specified by a finite amount of data, namely the posets of objects and operations, with source, target, unit, and composition maps. We show that $\Cfrak$ is in fact uniquely determined from its categories of objects and operations, together with source and target maps - this uniqueness criterion allows us later on to compare our approach to enrichment with other approaches in the literature.
 
 In \ref{subsection assosx} we give the definition of the associative operad $\Assos_X$ for an arbitrary category $X$. Although for the purposes of enriched category theory the category $X$ will always be a space, we will use this extra generality later on to give a direct description of the equivalence between $\Cat$ and the category of categories enriched in $\Spc$. We also give definitions of the related operads $\LM_X$ and $\BM_{X,Y}$ which corepresent left modules and bimodules.

In \ref{subsection algebroids} we review the definition and functoriality of the category of algebroids $\Algbrd(\Mcal)$ in an associative operad $\Mcal$. We pay special attention to the case when $\Mcal$ is a presentable monoidal category - in this case we have that $\Algbrd(\Mcal)$ is also presentable. We introduce two basic examples of algebroids: the trivial algebroid, and the cells - together these generate $\Algbrd(\Mcal)$. 

In \ref{subsection enriched} we review the case of $\Spc$-algebroids with a space of objects, and its equivalence with the category of Segal spaces. We then define the category of $\Mcal$-enriched categories $\Cat^{\Mcal}$ as the full subcategory of $\Algbrd(\Mcal)$ on those algebroids with a space of objects and whose underlying Segal space is complete and reprove the basic fact that if $\Mcal$ is presentable monoidal then $\Cat^{\Mcal}$ is an accessible localization of $\Algbrd(\Mcal)$.

In \ref{subsection multiplicativity} we discuss the canonical symmetric monoidal structure in the category of algebroids over a symmetric monoidal category. In the presentable setting, this gives access in particular to a notion of functor algebroids and functor enriched categories. We prove here a basic result describing Hom objects in functor algebroids when the source algebroid is a cell, which will later on be used as a starting point for establishing various facts about general functor algebroids. We finish by studying the behavior of functor algebroids as we change the enriching category.

In \ref{subsection ncat} we review the approach to $n$-categories as categories enriched in $(n-1)$-categories. We discuss the various functors relating the categories $\nCat$ for different values of $n$. In the limit as $n$ tends to infinity we recover the category $\omega\kr\Cat$ of $\omega$-categories. Although for our purposes all of the $\omega$-categories we will encounter will be $n$-categories for some finite $n$, the theory of $\omega$-categories provides a convenient setting in which to work with $n$-categories in cases where the exact value of $n$ is irrelevant or may vary. We show that the theory of $\omega$-category is in fact a fixed point under enrichment: there is an equivalence between $\omega\kr\Cat$ and the category of categories enriched in $\omega\kr\Cat$.

\begin{remark}
The theory of algebroids and enriched categories was introduced in \cite{GH} and \cite{Hinich}. In this paper we introduce a new approach to the subject based on the internal cooperad $\Cfrak$, and show that this approach arrives at the same theory as that from \cite{Hinich}. 

Some of the basic facts about algebroids and enriched categories  that we discuss in \ref{subsection algebroids}-\ref{subsection multiplicativity} (for instance, claims about presentability, existence of symmetric monoidal structures, functoriality of the theory) appear already in some way in the references. We chose to include statements and proofs of most of those facts for completeness and ease of reference, as our notation and conventions differ from other sources. 

Another reason why we opted for a systematic treatment of the subject is that in many cases we in fact need tools that go beyond those which appear in the literature. For instance, in \ref{subsection multiplicativity} we show that the category of algebroids over a symmetric operad admits the structure of a symmetric operad, and that this structure is functorial under morphisms of operads - this functoriality will be necessary in section 5 to construct the realization functor $\psi_n$. We also pay special attention throughout the section to enriched cells. We are able to obtain a good understanding of products and functor categories for cells, which is a basic building block for proving results for arbitrary functor categories later on. In particular, this will be crucial in section 4 when we discuss adjunctions and conical limits in functor enriched categories.
\end{remark}

\subsection{The internal cooperad $\Cfrak$} \label{subsection cfrak}

We begin with a general discussion of the procedure of internalization of objects of arbitrary presentable categories.

\begin{definition}
Let $\Dcal$ be a presentable category and let $\Ccal$ be a complete category. A $\Dcal$-object internal to $\Ccal$ is a limit preserving functor $F: \Dcal^\op \rightarrow \Ccal$. We denote by $\Dcal(\Ccal)$ the full subcategory of $\Funct(\Dcal^\op, \Ccal)$ on the internal $\Dcal$-objects.
\end{definition}

\begin{example}
Let $\Dcal$ be a presentable category. Then it follows from \cite{HTT} proposition 5.5.2.2 that $\Dcal(\Spc)$ is equivalent to $\Dcal$.
\end{example}

\begin{example}\label{example presh}
Let $\Ccal$ be a complete category. Then $\Spc(\Ccal)$ is equivalent to $\Ccal$. More generally, if $\Dcal'$ is a small category then  $(\Pcal(\Dcal'))(\Ccal)$ is equivalent to $\Funct(\Dcal'^\op, \Ccal)$.
\end{example}

\begin{remark}\label{remark functor in the other direction}
Let $\Dcal$ be a presentable category and let $\Ccal$ be a locally small complete category. It follows from \cite{HTT} proposition 5.5.2.2 that a functor $F: \Dcal^\op \rightarrow \Ccal$ is limit preserving if and only if it has a left adjoint. In this context, the data of $F$ is equivalent to the data of a functor $G: \Ccal^\op \rightarrow \Dcal$ such that for every $d$ in $\Dcal$ the presheaf $\Hom_{\Dcal}(d, G(-))$ on $\Ccal$ is representable. If $\Ccal$ is presentable then this condition is equivalent to $G$ preserving limits, and so we conclude that $\Dcal(\Ccal) = \Ccal(\Dcal)$. In other words, if $\Ccal$ and $\Dcal$ are presentable then $\Dcal$-objects internal to $\Ccal$ are the same as $\Ccal$-objects internal to $\Dcal$. Indeed, in this case the category $\Dcal(\Ccal) = \Ccal(\Dcal)$ admits a symmetric presentation as $\Ccal \otimes \Dcal$ (see \cite{HA} proposition 4.8.1.17). 
\end{remark}

\begin{remark}\label{remark localization}
Let $L: \Dcal_1 \rightarrow \Dcal_2$ be a localization functor between presentable categories and let $\Ccal$ be a locally small complete category. Then the functor $\Dcal_2(\Ccal) \rightarrow \Dcal_1(\Ccal)$ given by precomposition with $L$ is a fully faithful embedding. A $\Dcal_1$-object $F:\Dcal_1^\op \rightarrow \Ccal$ belongs to $\Dcal_2(\Ccal)$ if and only if the associated functor $G: \Ccal^\op \rightarrow \Dcal_1$ factors through $\Dcal_2$.
\end{remark}

\begin{remark}\label{remark class}
Let $\Dcal$ be a presentable category and let $\Ccal$ be a classical locally small complete category. Let $\Dcal_{\leq 0}$ be the full subcategory of $\Dcal$ on the $0$-truncated objects and denote by $\tau_{\leq 0} : \Dcal \rightarrow \Dcal_{\leq 0}$ the truncation functor. Then it follows from remark \ref{remark localization} that precomposition with $\tau_{\leq 0}$ induces an equivalence $\Dcal_{\leq 0}(\Ccal) = \Dcal(\Ccal)$.
\end{remark}

\begin{example}
Let $\Ccal$ be a  classical locally small complete category. Then $\Set(\Ccal)$ is equivalent to $\Ccal$. If $\Ccal$ is presentable then $\Ccal(\Set)$ is also equivalent to $\Ccal$.
\end{example}

\begin{remark}\label{remark presentation}
Let $\Dcal$ be a presentable category. Let $\Dcal'$ be a small category equipped with a localization functor $L: \Pcal(\Dcal') \rightarrow \Dcal$, so that the right adjoint to $L$ embeds $\Dcal$ as a full subcategory of $\Pcal(\Dcal')$. Let $\Ccal$ be a locally small complete category. Then it follows from a combination of example \ref{example presh} and remark \ref{remark localization} that $\Dcal(\Ccal)$ is equivalent to the full subcategory of $\Funct(\Dcal'^\op, \Ccal)$ on those functors $F$ such that the presheaf $\Hom_{\Ccal}(c, F(-))$ belongs to $\Dcal$ for every $c$ in $\Ccal$.
\end{remark}

We now specialize the above discussion to the case of internal operads.

\begin{notation}
  Denote by $\Op$ the category of operads. For each operad  $\Ocal$ we  denote by $\Op_\Ocal$ the category of operads over $\Ocal$.
\end{notation}

\begin{definition}
Let $\Ocal$ be an operad, and $\Ccal$ be a complete category. An $\Ocal$-operad internal to $\Ccal$ is an $\Op_\Ocal$-object internal to $\Ccal$. If $\Ccal'$ is a cocomplete category  then an $\Ocal$-cooperad internal is an $\Ccal'$ is an $\Ocal$-operad internal to $(\Ccal')^\op$.
\end{definition}

\begin{remark}\label{remark access}
Let $\Ocal$ be an operad, and $\Ccal$ be a presentable category. Then following remark \ref{remark functor in the other direction}, we see that an $\Ocal$-cooperad internal to $\Ccal$ is the same data as an accessible, limit preserving functor $G: \Ccal \rightarrow \Op_\Ocal$.
\end{remark}

\begin{remark} \label{remark classical cooperad}
Let $\Ocal$ be an operad and $\Ccal$ be a classical locally small complete category. Then by virtue of remark \ref{remark class} we have an equivalence $\Op_\Ocal(\Ccal) = (\Op_{\Ocal})_{\leq 0}(\Ccal)$. In other words, $\Ocal$-operads internal to $\Ccal$ are the same as $0$-truncated $\Ocal$-operads internal to $\Ccal$.

Assume now that $\Ocal$ is a $0$-truncated object of $\Op$, so that it has a set $V$ of objects, and a set $M$ of operations. Consider the full subcategory $\Dcal'$ of $(\Op_\Ocal)_{\leq 0}$ on the following objects:
\begin{itemize}
\item The trivial $\Ocal$-operad $\vfrak_o$ for each $o$ in $V$.
\item For every operation $m$ in $M$, the free $\Ocal$-operad $\mathfrak{f}_m$ containing an $m$-operation. 
\item For every operation $m$ in $M$ with source objects $\lbrace o_i \rbrace_{i \in S}$, and every family of operations $\lbrace m_i \rbrace_{i \in S}$ where the target object of $m_i$ is $o_i$, the $\Ocal$-operad $\ffrak_{m_i, m}$ defined as the pushout
\[
\left(\bigsqcup_{i \in S}  \ffrak_{m_i} \right) \bigcup_{\bigsqcup_{i \in S} \vfrak_{o_i}} \ffrak_{m}.
\]
\end{itemize}

Then $(\Op_\Ocal)_{\leq 0}$ is a localization of the category of set valued presheaves on $\Dcal'$. Using remark \ref{remark presentation} are able to obtain a concrete description of the category $\Op_\Ocal(\Ccal)$. Namely, an $\Ocal$-operad $\Ocal'$ internal to  $\Ccal$ consist of the following data:
\begin{itemize}
\item For each $o$ in $V$ an object $\Ocal'(\vfrak_o)$ in $\Ccal$.
\item For each $m$ in $M$ with source objects $\lbrace o_i \rbrace_{i \in S}$ and target object $o$, an object $\Ocal'(\ffrak_m)$ in $\Ccal$, equipped with source and target maps $\prod_{i \in S} \Ocal'(\vfrak_{o_i}) \leftarrow \Ocal'(\ffrak_m) \rightarrow \Ocal'(\vfrak_{o})$.
\item For each $o$ in $V$ a unit map $\Ocal'(\vfrak_{o}) \rightarrow \Ocal'(\ffrak_{\id_o})$, where $\id_o$ denotes the identity $1$-ary operation of $o$.
\item For each $m$ in $M$ with source objects $\lbrace o_i \rbrace_{i \in S}$, and every family of operations $\lbrace m_i \rbrace_{i \in S}$ where the target object of $m_i$ is $o_i$, a composition map
\[
\prod_{i\in S} \Ocal'(\ffrak_{m_i}) \times_{\prod_{i \in S} \Ocal'(\vfrak_{o_i})} \Ocal'(\ffrak_{m}) \rightarrow \Ocal(\ffrak_l). 
\]
where $l$ denotes the composite in $\Ocal$ of family of operations $\lbrace m_i \rbrace_{i \in S}$ with $m$.
\end{itemize}
The above data is required to satisfy a finite list of standard compatibility conditions mimicking those of the category of $0$-truncated $\Ocal$-operads, built so that the data obtained from the above by applying a corepresentable functor $\Ccal \rightarrow \Set$ defines an $\Ocal$-operad in $\Set$. In other words:
\begin{itemize}
\item Composition and unit maps are required to be compatible with sources and targets.
\item Units are required to be compatible with compositions.
\item Composition is required to be associative.
\item For every $o$ in $V$ the unit map induces an isomorphism between $\Ocal'(\vfrak_o)$ and the subobject of isomorphisms inside $\Ocal'(\ffrak_{\id_o})$.
\end{itemize}
\end{remark}

\begin{example}\label{example nonsymm}
Let $\Ccal$ be a locally small complete category and let $\Assos$ be the operad governing associative algebras. Then $\Assos$-operads internal to $\Ccal$ will be called internal nonsymmetric operads. By virtue of remark \ref{remark classical cooperad}, in the case when $\Ccal$ is classical we can specify a nonsymmetric operad internal to $\Ccal$ by giving a finite amount of information. Namely, a nonsymmetric operad $\Ocal$ internal to $\Ccal$ consists of the following data:
\begin{itemize}
\item An object $V$ in $\Ccal$ parametrizing objects of $\Ocal$.
\item For each $n \geq 0$ an object  $M_n$ in $\Ccal$ parametrizing $n$-ary operations in $\Ocal$, equipped with source and target maps $V^n \leftarrow M_n \rightarrow V$.
\item A unit map $V \rightarrow M_1$.
\item For each $n \geq 0$ and each sequence $\lbrace n_i \rbrace_{1 \leq i \leq n}$ of nonnegative integers with sum $N$ a composition map \[
\prod_{1 \leq i \leq n} M_{n_i} \times_{V^n} M_n \rightarrow M_{N}
\]
subject to the conditions described in remark \ref{remark classical cooperad}.
\end{itemize}
\end{example}

We now present the construction of the internal nonsymmetric cooperad in strict categories $\mathfrak{C}^\cl$ which underlies  the assignment $X \mapsto \Assos_X$. As discussed in example \ref{example nonsymm}, we can do this by specifying a finite amount of information.  

\begin{construction}\label{construction cooperad}
Let $\Cat_{\leq 0}$ be the category of $0$-truncated categories. In other words, this is the category of strict categories with no nontrivial isomorphisms. We define a nonsymmetric cooperad $\mathfrak{C}^\cl$ internal to $\Cat_{\leq 0}$ as follows:
\begin{itemize}
\item The category $V$ of objects of $\mathfrak{C}^\cl$ is the set with two elements $\lbrace t, s \rbrace$.
\item For each $n \geq 0$ the category $M_n$ of $n$-ary operations of $\mathfrak{C}^\cl$ is in fact a poset, and has objects $t_i, s_i$ for $0 \leq i \leq n+1$, with $t_0 = s_0$ and $t_{n+1} = s_{n+1}$, and arrows $s_i \leftarrow t_{i+1}$ for $0 \leq i \leq n$. We depict this as follows:
\[
\begin{tikzcd}
t_0 = s_0 & \arrow{l}{} t_1 & s_1 & \arrow{l}{} t_2 & \ldots & s_{n} & t_{n+1} = s_{n+1} \arrow{l}{} 
\end{tikzcd}
\]
\item The cosource map $V^{\amalg n} \rightarrow M_n$ maps the $i$-th copy of $t$ and $s$ to $t_i$ and $s_i$ respectively. The cotarget map $M_n \leftarrow V$ maps $t, s$ to $t_0$ and $s_{n+1}$ respectively.
\item The counit map $M_1 \rightarrow V$ maps $t_0$ and $t_1$ to $t$ and $s_1$ and $s_2$ to $s$.
\item Let $n \geq 0$ and let $\lbrace n_j \rbrace_{1 \leq j \leq n}$ be a sequence of nonnegative integers with sum $N$. Denote the objects of $\bigsqcup_j M_{n_j}$ by $t_k^j, s_k^j$, where $1 \leq j \leq n$ and $0 \leq k\leq n_j+1$. The poset 
\[
\left( \bigsqcup_{1 \leq i \leq n} M_{n_i} \right) \bigcup_{V^{\amalg n}} M_n 
\] has two extra objects which are in the image of the cotarget map $M_n \leftarrow V$. We denote these by $t'_0$ and $s'_{n+1}$.  The cocomposition map 
\[
M_N \rightarrow \left( \bigsqcup_i M_{n_i} \right) \bigcup_{V^{\amalg n}} M_n
\]
 sends $t_0, s_{N+1}$ to $t'_0$ and $s'_{n+1}$, and for $i \neq 0, N$ sends $t_i$ and $s_i$ to $t_k^j$ and $s_k^j$ respectively, where $(j,k)$ is the unique pair with $1 \leq k \leq n_j$ such that $(\sum_{1 \leq l < j} n_l) + k = i$.
\end{itemize}
\end{construction}

\begin{remark}\label{remark uniqueness}
For each $n \geq 0$ the cosource and cotarget maps $V^n \rightarrow M_n \leftarrow V$ are jointly surjective. It follows from this, that $\mathfrak{C}^\cl$ is characterized uniquely by the first three items of construction \ref{construction cooperad}. In other words, there is a unique way in which we could have defined the counit and cocomposition maps for $\Cfrak^\cl$ once we are given the data of  $V, M_n$ and the cosource and cotarget maps.
\end{remark}

\begin{remark}\label{remark cosimplicial poset}
The internal nonsymmetric cooperad $\mathfrak{C}^\cl$ defines a colimit preserving functor $ (\Op_{\Assos})_{\leq 0} \rightarrow \Cat_{\leq 0}$ which we continue denoting by $\Cfrak^\cl$. 
 Composing this with the inclusion $\Delta \rightarrow (\Cat)_{\leq 0} \rightarrow (\Op_{\Assos})_{\leq 0}$ we obtain a cosimplicial $0$-truncated category $\Delta \rightarrow \Cat_{\leq 0}$. This satisfies the Segal conditions, and is in fact a cocategory object in $\Cat_{\leq 0}$. Inspecting construction \ref{construction cooperad} reveals that this is the functor that sends $[n]$ to the poset $[n] \bigsqcup [n]^\op$.
\end{remark}

Our next task is to extend $\Cfrak^\cl$ to a nonsymmetric cooperad internal to $\Cat$.

\begin{proposition}\label{lemma uniqueness}
There is a unique internal cooperad $\Cfrak: \Op_{\Assos} \rightarrow \Cat$ whose categories of objects and operations are $0$-truncated, and making the following triangle commute
\[
\begin{tikzcd}
\Op_{\Assos} \arrow[r, "\Cfrak"] \arrow[rd, swap, "\Cfrak^\cl"] & \Cat \arrow[d, "\tau_{\leq 0}"] \\
                                  & \Cat_{\leq 0} 
\end{tikzcd}
\]
where $\tau_{\leq 0}$ is left adjoint to the inclusion.
\end{proposition}
\begin{proof}
Recall the presentation of $\Op_{\Assos}$ in terms of complete Segal operads from \cite{Barwick}, as a localization of the presheaf category on the category $\Delta_{\Obb}$ of trees. This identifies $\Op_{\Assos}$ with the full subcategory of $\Pcal(\Delta_{\Obb})$ on those presheaves satisfying suitable Segal and completeness conditions. Let $L: \Pcal(\Delta_{\Obb}) \rightarrow \Op_{\Assos}$ be the localization functor, and $i : \Delta_\Obb \rightarrow \Pcal(\Delta_{\Obb})$ be the inclusion. Let $\vfrak$ be the terminal associative operad, thought of as an object of $\Delta_\Obb$, and for each $n \geq 0$ let $\ffrak_n$ be the free associative operad on an operation of arity $n$, again thought of as an object of $\Delta_\Obb$.

Since $L$ is a localization, and in particular an epimorphism, it suffices to show that there is a unique colimit preserving functor $\Cfrak' : \Pcal(\Delta_\Obb) \rightarrow \Cat$ which factors through $\Op_{\Assos}$, maps the objects $\vfrak$ and $\ffrak$ to $0$-truncated categories,  and makes the following triangle commute:
\[
\begin{tikzcd}
\Pcal(\Delta_\Obb) \arrow[r, "\Cfrak'"] \arrow[rd, swap, "\Cfrak^\cl L"] & \Cat \arrow[d, "\tau_{\leq 0}"] \\
                                  & \Cat_{\leq 0} 
\end{tikzcd}
\]
For this it suffices to show that there is a unique functor $G: \Delta_{\Obb} \rightarrow \Cat$ which satisfies the dual Segal and completeness conditions, maps $\vfrak$ and $\ffrak$ to $0$-truncated categories, and makes the following triangle commute:
\[
\begin{tikzcd}
\Delta_\Obb \arrow[r, "G"] \arrow[rd, swap, "\Cfrak^\cl Li"] & \Cat \arrow[d, "\tau_{\leq 0}"] \\
                                  & \Cat_{\leq 0} 
\end{tikzcd}
\] 

Let $j: \Cat_{\leq 0} \rightarrow \Cat$ be the inclusion. We claim that $G_0 = j \Cfrak^\cl L i$ satisfies the dual Segal and completeness conditions. The fact that $G_0$ satisfies the dual completeness conditions follows from the description of the simplicial category underlying $G_0$ from \ref{remark cosimplicial poset}. The fact that $G_0$ satisfies the dual Segal conditions follows from the fact that $j$ preserves the pushouts involved in them.

Note that $G_0$ comes equipped with an identification $\epsilon: \tau_{\leq 0} G_0 = \Cfrak^\cl L i$ given by the counit of the adjunction $\tau_{\leq 0} \dashv j$. It now suffices to show that pair $(G, \rho)$ of a functor $G$ and an identification $\rho: \tau_{\leq 0} G = \Cfrak^\cl L i$ as above is canonically equivalent to $(G_0, \epsilon)$. Let $\eta: G \rightarrow j\tau_{\leq 0}G$ be the unit map. The fact that $G$ maps $\vfrak$ and $\ffrak_n$ to $0$-truncated categories implies that $\eta$ is an isomorphism on the full subcategory of $\Delta_\Obb$ on the objects $\vfrak$ and $\ffrak_n$. Since both $G$ and $j \tau_{\leq 0} G = G_0$ satisfy the dual Segal conditions, we conclude that $(j\rho) \circ \eta$ gives us  an isomorphism $G = G_0$. Our claim now follows from the fact that the diagram of functors and natural isomorphisms
\[
\begin{tikzcd}[column sep = large]
\tau_{\leq 0} G \arrow{dr}{\rho}\arrow{r}{\tau_{\leq 0} \eta}  & \tau_{\leq 0}j\tau_{\leq 0} G \arrow{r}{\tau_{\leq 0}j\rho} & \tau_{\leq 0} G_0 \arrow{dl}[swap]{\epsilon} \\
 & \Cfrak^\cl L i &  
\end{tikzcd}
\]
commutes in a natural way.
\end{proof}

\begin{corollary}\label{coro uniqueness extended}
Let $C: \Op_{\Assos} \rightarrow \Cat$ be an internal associative cooperad. Assume that the category of objects and category of operations of $C$ are equivalent to those of $\Cfrak$, with an equivalence that commutes with the cosource and cotarget maps. The $C$ is equivalent to $\Cfrak$.
\end{corollary}
\begin{proof}
By remark \ref{remark uniqueness} we have that $\tau_{ \leq 0} C$ is equivalent to $\Cfrak^\cl$. The claim now follows from proposition \ref{lemma uniqueness}.
\end{proof}

\begin{remark}\label{remark passing op Cfrak}
Recall that $\Cat$ and $\Op_{\Assos}$ come equipped with involutions $(-)^\op$ and $(-)^\rev$, which correspond to actions of $\ZZ/ 2\ZZ$ on both categories. The internal cooperad $\Cfrak^\cl : \Op_{\Assos} \rightarrow \Cat_{\leq 0}$ can be given a $\ZZ/2\ZZ$- equivariant structure, by  switching  the role of $t$ and $s$ in the category of objects, and of $s_i, t_i$ with $t_{n+1 - i}$ and $s_{n+1-i}$ in the categories of operations. It follows from proposition \ref{lemma uniqueness} that the internal cooperad $\Cfrak$ inherits a $\ZZ/ 2\ZZ$ equivariant structure. In particular, we have a commutative square
\[
\begin{tikzcd}
\Op_{\Assos} \arrow{d}{\rev} \arrow{r}{\Cfrak} & \Cat \arrow{d}{\op} \\
\Op_{\Assos} \arrow{r}{\Cfrak} & \Cat.
\end{tikzcd}
\]
\end{remark}

\subsection{The operad $\Assos_X$}\label{subsection assosx}

We now introduce the operad $\Assos_X$ that corepresents the assignment $\Mcal \mapsto \Algbrd_X(\Mcal)$.

\begin{notation}
Let $\Assos_- : \Cat \rightarrow \Op_{\Assos}$ be the right adjoint to $\Cfrak$. This sends each category $X$ to an associative operad $\Assos_X$.
\end{notation}

\begin{remark}\label{remark structure assosx}
In the case when $X$ is the terminal category, the operad $\Assos_X$ coincides with the associative operad $\Assos$. In general, we think about $\Assos_X$ as a many object version of $\Assos$. The category of objects of $\Assos_X$ is  $X \times X^\op$. Given a nonempty sequence $\lbrace (y_i, x_i) \rbrace_{1 \leq i \leq n}$ of source objects, and a target object $(y, x)$, a multimorphism $\lbrace (y_i, x_i) \rbrace_{1 \leq i \leq n} \rightarrow (y, x)$ consists of a series of arrows $x_i \leftarrow y_{i+1}$  for $1 \leq i < n$, and arrows $x_n \leftarrow x$ and $y \leftarrow y_1$. A multimorphism from the empty sequence of objects to $(y, x)$ consists of an arrow $y \leftarrow x$.
\end{remark}

\begin{remark} \label{remark coincide hinich}
In \cite{Hinich}, Hinich works in the language of categories of operators, and defines an assignment $\Assos^H_- : \Cat \rightarrow \Op_{\Assos}$ to be corepresented by a certain functor $\Fcal: \Delta/\Delta^\op \rightarrow \Cat_{\leq 0}$. This functor is the categories of operators incarnation of the internal nonsymmetric cooperad $\mathfrak{C}$.

 Indeed, note that the functor $\Assos_-^H$ is accessible and preserves limits, so by virtue of remark \ref{remark access} it is corepresented by an internal nonsymmetric cooperad $\Cfrak^H: \Op_{\Assos} \rightarrow \Cat$.  Direct inspection of the definition of $\Fcal$ reveals that the category of objects and operations of $\Cfrak^H$ agree with those of  $\Cfrak^\cl$, in a way which is compatible with source and target maps. As observed in corollary \ref{coro uniqueness extended} this implies that $\Cfrak^H$ is equivalent to $\Cfrak$. It follows that the functor $\Assos^H$ defined in \cite{Hinich} is equivalent to our functor $\Assos$.
 \end{remark}
 
\begin{example}\label{example set w two elements}
Let $X = \lbrace a, b \rbrace$ be the set with two elements $a, b$. Then the associative operad $\Assos_X$ is classical, and can be computed explicitly from the definitions. We note that it has objects $(a, a), (b,b), (a, b), (b, a)$. The objects $(a, a)$ and $(b, b)$ are algebras in $\Assos_X$. The object $(a, b)$ is an $(a, a)-(b,b)$ bimodule and the object $(b, a)$ is a $(b,b)-(a,a)$ bimodule.
\end{example}

As we shall see below, for each pair of categories $X, Y$, the category $X \times Y^\op$ has compatible (weak) actions of the associative operads $\Assos_X$ and $\Assos_Y$ on the left and on the right. This is the basis for the theory of bimodules over algebroids.

\begin{notation}
Denote by $\BM, \LM, \RM$ be the associative operads governing bimodules left modules, and right modules, respectively. Recall that we have canonical inclusions $\LM \rightarrow \BM \leftarrow \RM$. We denote by $\Assos^+$ and $\Assos^-$ the copies of the associative operad in $\BM$  contained in $\LM$ and $\RM$, respectively.
\end{notation}

\begin{construction}\label{construction BMXY}
Let $X$ and $Y$ be categories. Consider the projection $X \sqcup Y \rightarrow \lbrace a, b \rbrace$ from the disjoint union of $X$ and $Y$ into the set with two elements $a, b$, that maps $X$ to $a$ and $Y$ to $b$. Applying the functor $\Assos_-$ we obtain a map of associative operads $\Assos_{X \sqcup Y} \rightarrow \Assos_{\lbrace a , b \rbrace}$. We let $\BM_{X,Y}$ be the $\BM$-operad obtained by pullback of $\Assos_{X \sqcup Y}$ along the map $\BM \rightarrow \Assos_{\lbrace a, b \rbrace}$ corresponding to the $(a, a)-(b,b)$ bimodule $(a, b)$ (see example \ref{example set w two elements}). The assignment $X \mapsto \BM_{X, Y}$ is functorial in $X$ and $Y$. We denote by $\BM_{-,-}: \Cat \times \Cat \rightarrow \Op_{\BM}$ the corresponding functor.
\end{construction}

\begin{remark}
Let $X, Y$ be categories. The functor $\Assos_-$ preserves limits since it is a right adjoint. Applying it to the cartesian square
\[
\begin{tikzcd}
X \arrow{d}{} \arrow{r}{} & X\sqcup Y \arrow{d}{} \\
\lbrace a \rbrace \arrow{r}{} & \lbrace a, b \rbrace
\end{tikzcd}
\]
we conclude that the $\Assos^-$-component of $\BM_{X,Y}$ coincides with $\Assos_X$. This equivalence is natural in $X$.

 Similarly, from the cartesian square
\[
\begin{tikzcd}
Y \arrow{d}{} \arrow{r}{} & X \sqcup Y \arrow{d}{} \\
\lbrace b \rbrace \arrow{r}{} & \lbrace a, b \rbrace
\end{tikzcd}
\]
we see that the $\Assos^+$-component of $\BM_{X,Y}$ coincides with $\Assos_Y$. 

Consider now the cartesian square
\[
\begin{tikzcd}
(\BM_{X, Y})_m \arrow{d}{} \arrow{r}{} & \Assos_{X \sqcup Y} \arrow{d}{} \\
\lbrace (a, b) \rbrace \arrow{r}{} & \Assos_{\lbrace a, b \rbrace}
\end{tikzcd}
\]
where the category $(\BM_{X, Y})_m$ is the fiber of $\BM_{X, Y}$ over the universal bimodule $m$ in $\BM$. Applying the (limit preserving) forgetful functor $\Op_{\Assos} \rightarrow \Cat$ we obtain a cartesian square
\[
\begin{tikzcd}
(\BM_{X, Y})_m \arrow{d}{} \arrow{r}{} & (X \times X^\op) \sqcup (Y \times Y^\op) \sqcup (X \times Y^\op) \sqcup (Y \times X^\op) \arrow{d}{} \\
\lbrace (a, b) \rbrace \arrow{r}{} & \lbrace (a, a), (b,b), (a,b), (b, a)\rbrace.
\end{tikzcd}
\]
We conclude that the category $(\BM_{X, Y})_m$ is equivalent to $X \times Y^\op$. This equivalence is also natural in $X, Y$.
\end{remark}

\begin{example}\label{remark diagonal BMX to assos}
Let $X$ be a category. Then we have an equivalence
\[
\Assos(X \sqcup X) = \Assos(X \times \lbrace a, b \rbrace ) = \Assos(X) \times \Assos(\lbrace a , b \rbrace)
\]
which is natural in $X$. It follows that we have an equivalence $\BM_{X, X} = \Assos_X \times \BM$, which is natural in $X$.
\end{example}

\begin{notation}\label{notation BMX}
Let $X$ be a category. Denote by $\BM_X$ the associative operad $\BM_{X, [0]}$. Let $\LM_X$  be the $\LM$-operad obtained by pullback of $\BM_X$ along the inclusion $\LM \rightarrow \BM$. We denote by $\BM_-: \Cat \rightarrow \Op_{\BM}$ the functor that assigns to each category $X$ the $\BM$-operad $\BM_X$, and by $\LM_-$  the composition of $\BM_-$ with the functor of base change to $\LM$.
\end{notation}

\begin{remark}
In \cite{Hinich}, Hinich defines an assignment $\BM^H_- : \Cat \rightarrow \Op_{\BM}$  in the language of categories of operators, by declaring it to be corepresented by a certain functor \[\Fcal_{\BM} : \Delta/{(\Delta_{/[1]})^\op} \rightarrow \Cat.
\] The functor $\BM^H_-$ is accessible and preserves limits, and therefore by remark \ref{remark access} it is corepresented by an internal $\BM$-cooperad $\Cfrak_{\BM}^H : \Op_{\BM} \rightarrow \Cat$.

Likewise, our functor $\BM_-$ can be obtained as the composite functor
\[
\Cat = \Cat_{/ \lbrace a \rbrace} \xrightarrow{- \sqcup \lbrace b \rbrace} \Cat_{/ \lbrace a, b \rbrace} \xrightarrow{\Assos_-} (\Op_{\Assos})_{/\Assos_{\lbrace a, b \rbrace} } \xrightarrow{ {}_{(a,a)}  (a, b)_{(b, b)} } (\Op_{\Assos})_{/\BM} = \Op_{\BM}
\]
and each of the functors in the composition preserves limits and is accessible, so we have that $\BM_-$ is also corepresented by an internal $\BM$-cooperad $\Cfrak_{\BM}: \Op_{\BM} \rightarrow \Cat$. Direct inspection of the functor $\Fcal_{\BM}$ reveals that $\Cfrak_{\BM}$ and $\Cfrak_{\BM}^H$ have equivalent categories of objects and operations, in a way compatible with sources and target. A variation of the arguments in proposition \ref{lemma uniqueness} and corollary \ref{coro uniqueness extended} (where we work with $(\Delta_\Obb)_{/\BM}$ ) shows that the cooperads $\Cfrak_{\BM}$ and $\Cfrak_{\BM}^H$ are equivalent, and thus our functor $\BM_-$ is equivalent to the functor $\BM_-^H$ from \cite{Hinich}.
\end{remark}

\subsection{Algebroids}\label{subsection algebroids}

We now present the definition of an algebroid in an associative operad (also known as categorical algebras in \cite{GH} and enriched precategories in \cite{Hinich}) and review the basic funtoriality properties of the theory.

\begin{definition}
Let $\Mcal$ be an associative operad and $X$ be a category. An algebroid on $\Mcal$ with category of objects $X$ is an $\Assos_X$-algebra in $\Mcal$.
\end{definition}

\begin{remark}
Let $\Mcal$ be an associative operad. An algebroid with category of objects $[0]$ is an algebra in $\Mcal$. In general, we think about an algebroid $\Acal$ in $\Mcal$ as a many-object associative algebra. Indeed, an algebroid with category of objects $X$ assigns to each pair of objects $y, x$ in $X$ an object $\Acal(y, x)$ in $\Mcal$ and to every $n \geq 0$ and every sequence of arrows $y_0 = x_0 \leftarrow y_1, x_1 \leftarrow y_2, \ldots, x_{n-1} \leftarrow y_n, x_n \leftarrow y_{n+1} = x_{n+1}$ it assigns a multimorphism 
\[
\lbrace \Acal(y_1, x_1) , \Acal(y_2, x_2) ,\ldots , \Acal(y_n, x_n) \rbrace \rightarrow \Acal(y_0, x_{n+1})
\]
in $\Mcal$. In the case when $\Mcal$ is a monoidal category, this is the same as a morphism
\[
\Acal(y_1, x_1) \otimes \Acal(y_2, x_2) \otimes \ldots, \otimes \Acal(y_n, x_n) \rightarrow \Acal(y_0, x_{n+1}).
\]

Specializing to the case $n = 0$ we obtain for every pair of objects $x, y$ a map 
\[
\Hom_X(y, x) \rightarrow \Hom_\Mcal(1_\Mcal, \Acal(y, x)).
\] In particular, starting from the identity in $\Hom_X(x , x)$ we obtain a map $1_\Mcal \rightarrow \Acal(x, x)$ (the unit at $x$). 

In the case when $n = 2$ and all the arrows are identities we obtain a map 
\[
\Acal(z, y) \otimes \Acal(y, x) \rightarrow \Acal(z, x)
\] (the composition map) for every triple of objects $x, y, z$ in $\Mcal$. It follows from the definition of the composition rule in the cooperad $\Cfrak$ that these maps satisfy the usual unitality and associativity rules up to homotopy.
\end{remark}

\begin{construction}\label{construction algbrd}
Let $X$ be a category and $\Mcal$ be an associative operad. We denote by $\Algbrd_X(\Mcal)$ the category of $\Assos_X$-algebras in $\Mcal$. Denote by $\Algbrd_-(-): \Cat^\op \times \Op_{\Assos} \rightarrow \Cat$ the composite functor
\[
\Cat^\op \times \Op_{\Assos} \xrightarrow{\Assos_-^\op \times \id} \Op_{\Assos}^\op \times \Op_{\Assos} \xrightarrow{\Alg_-(-)} \Cat.
\]
For each $\Mcal$ in $\Op_{\Assos}$ we denote by $\Algbrd(\Mcal)$ the total category of the cartesian fibration associated to the functor $\Algbrd_-(\Mcal): \Cat^\op \rightarrow \Cat$. We call $\Algbrd(\Mcal)$ the category of algebroids in $\Mcal$.

The assignment $\Mcal \mapsto (\Algbrd(\Mcal) \rightarrow \Cat)$ yields a functor 
\[
\Algbrd(-): \Op_{\Assos} \rightarrow \widehat{\Cat}
\] equipped with a natural transformation to the constant functor $\Cat$. We denote by $\Algbrd$ the total category of the cocartesian fibration associated to $\Algbrd(-)$. This comes equipped with a projection $\Algbrd \rightarrow \Cat \times \Op_{\Assos}$ whose fiber over a pair $(X, \Mcal)$ is the category $\Algbrd_X(\Mcal)$. This is the two-sided fibration associated to the functor $\Algbrd(-)_-$.\footnote{We refer to \cite{HSTI} for background on the theory of two-sided fibrations and bifibrations, and a general discussion of the Grothendieck construction which relates these to bifunctors into $\Cat$.}
\end{construction}

\begin{notation}
Let $i: X \rightarrow Y$ be a functor of categories, and let $j: \Mcal \rightarrow \Ncal$ be a map of associative operads. We denote by $i^!: \Algbrd_Y(\Mcal) \rightarrow \Algbrd_X(\Mcal)$ the functor induced by $i$, and by $j_!: \Algbrd(\Mcal) \rightarrow \Algbrd(\Ncal)$ the functor induced by $j$.
\end{notation}

\begin{example}
Let $\Mcal$ be a monoidal category. Then the unit $1_\Mcal$ has an algebra structure, and therefore defines an algebroid with category of objects $[0]$. Since $1_\Mcal$ is initial in $\Algbrd_{[0]}(\Mcal)$, the functor $\Algbrd(\Mcal) \rightarrow \Spc$ corepresented by $1_\Mcal$ is equivalent to the restriction along the projection $\Algbrd(\Mcal) \rightarrow \Cat$ of the functor $\Cat \rightarrow \Spc$ corepresented by $[0]$. In particular, for every $\Mcal$-algebroid $\Acal$ with category of objects $X$, the space $\Hom_{\Algbrd(\Mcal)}(1_\Mcal, \Acal)$ is equivalent to the space underlying $X$. 
\end{example}

\begin{example}\label{example cells}
Let $\Mcal$ be monoidal category, and let $m$ be an object in $\Mcal$. Assume that $\Mcal$ admits an initial object which is compatible with the monoidal structure. Let $X = \lbrace a, b \rbrace$ be the set with two elements. Then we may form the free $\Assos_X$-algebra $C_m$ in $\Mcal$ equipped with a map $m \rightarrow C_m(a, b)$. The description of free algebras from \cite{HA} definition 3.1.3.1 yields the following description of $C_m$:
\begin{itemize}
\item $C_m((a, a)) = C_m((b, b)) = 1_\Mcal$.
\item $C_m(a, b) = m$.
\item $C_m(b, a)$ is the initial object of $\Mcal$.
\end{itemize}

We call the algebroids of the form $C_m$ cells. These come equipped with two maps $1_\Mcal \rightarrow C_m$, which pick out the objects $(a, a)$ (the target) and $(b, b)$ (the source). These maps can be combined into a single map out of the coproduct $1_\Mcal \sqcup 1_\Mcal$ in $\Algbrd(\Mcal)$ (note that this coproduct indeed exists and is given by the initial object in $\Algbrd_{\lbrace a, b \rbrace}(\Mcal)$, which agrees with the cell associated to the initial object in $\Mcal$).

 The formation of cells is functorial in $m$: it underlies a colimit preserving functor 
 \[
 C_- : \Mcal \rightarrow \Algbrd_{\lbrace a, b \rbrace}(\Mcal)
 \] 
given by operadic left Kan extension along the inclusion $\lbrace (a, b) \rbrace \rightarrow \Assos_{\lbrace a, b \rbrace}$. This assignment is furthermore functorial in $\Mcal$. In other words, given another monoidal category $\Mcal'$ with compatible initial object and  a monoidal functor $F: \Mcal \rightarrow \Mcal'$ which preserves initial objects, then the commutative square of categories
\[
\begin{tikzcd}
\Algbrd_{\lbrace a, b \rbrace}(\Mcal) \arrow{d}{F_!} \arrow{r}{\ev_{(a, b)}} & \Mcal \arrow{d}{F} \\
\Algbrd_{\lbrace a, b \rbrace}(\Mcal') \arrow{r}{\ev_{(a, b)}} & \Mcal'
\end{tikzcd}
\]
is horizontally left adjointable.
\end{example}

\begin{remark}\label{remark hom functor}
Let $u: \Op_{\Assos} \rightarrow \Cat$ be the colocalization functor that attaches to each associative operad its category of objects. It follows from \ref{remark cosimplicial poset} that we have an equivalence $u\Assos_- = \id \times (\id)^\op$ as endofunctors of $\Cat$. 
In particular, for every category $X$ and associative operad $\Mcal$ we have a functor
\[
\Algbrd_X(\Mcal) \rightarrow \Funct(X \times X^\op, \Mcal)
\]
which is natural in $X$ and $\Mcal$. We think about this as the functor which attaches to each algebroid $\Acal$ with category of objects $X$, the hom-functor of $\Acal$ restricted to $X$.
\end{remark}

\begin{remark}\label{remark op algebroids}
Recall from remark \ref{remark passing op Cfrak} that the cooperad $\Cfrak$ intertwines the order reversing involution $(-)^\rev$ of $\Op_{\Assos}$ and the passing to opposites involution of $\Cat$. It follows that the same is true for the functor $\Assos_-: \Cat \rightarrow \Op_{\Assos}$.  We thus see that the functor  $\Algbrd_-(-)$ admits the structure of a fixed point for the involution $(-)^\op \times (-)^\rev$ on $\Cat^\op \times \Op_{\Assos}$, which implies that there is an involution $(-)^\op$ on $\Algbrd$ and an enhancement of the projection $\Algbrd \rightarrow \Cat \times \Op_{\Assos}$ to a $\ZZ / 2\ZZ$-equivariant map. 

In particular, for every category $X$ and associative operad $\Mcal$, we have an equivalence $\Algbrd_X(\Mcal) = \Algbrd_{X^\op}(\Mcal^\rev)$.  In the case when $X$ is a space and $\Mcal$ is the associative operad underlying a symmetric operad, then the pair $(X, \Mcal)$ is has the structure of fixed point for the involution $(-)^\op \times (-)^\rev$. It follows that the involution $(-)^\op : \Algbrd \rightarrow \Algbrd$ restricts to an involution on $\Algbrd_X(\Mcal)$. In other words, if $\Mcal$ underlies a symmetric operad, then any $\Mcal$-algebroid $\Acal$ with a space of objects $X$ has attached to it another $\Mcal$-algebroid $\Acal^\op$ with space of objects $X$. Examining the fixed point structure on $\Cfrak^\cl$ from remark \ref{remark passing op Cfrak} reveals that for every pair of objects $y, x$ in $X$ one has $\Acal^\op(y, x) = \Acal(x, y)$.
\end{remark}

\begin{example}
Let $\Mcal$ be a symmetric monoidal category compatible with initial object and let $m$ be an object of $\Mcal$. Then the cell $C_m$ is equivalent to its opposite. This equivalence interchanges sources and targets - namely, there is a commutative diagram of $\Mcal$-algebroids
\[
\begin{tikzcd}[column sep = huge]
1_{\Mcal} \bigsqcup 1_{\Mcal} \arrow{r}{((a,a), (b,b))} \arrow{d}{\id} & C_m \arrow{d}{=} \\
1_{\Mcal} \bigsqcup 1_{\Mcal} \arrow{r}{((b,b), (a,a))} \arrow{r}{} & (C_m)^\op.
\end{tikzcd}
\]
In general, for any $\Mcal$-algebroid $\Acal$ we can think about $\Acal^\op$ as being obtained from $\Acal$ by reversing the direction of the cells.
\end{example}

\begin{remark}\label{some infty two functoriality}
Let $\Mcal, \Mcal'$ be associative operads. Then for every category $X$ there is a functor
\[
\Alg_\Mcal(\Mcal') \times \Algbrd_X(\Mcal) \rightarrow \Algbrd_X(\Mcal').
\]
This is natural in $X$ and therefore defines a functor
\[
\Alg_\Mcal(\Mcal') \times \Algbrd(\Mcal) \rightarrow \Algbrd(\Mcal')
\]
which enhances the functoriality of construction \ref{construction algbrd} to take into account natural transformations between morphisms of associative operads. This is compatible with composition: namely, given a third associative operad $\Mcal''$, there is a commutative square
\[
\begin{tikzcd}
\Alg_{\Mcal'}(\Mcal'') \times \Alg_\Mcal(\Mcal') \times \Algbrd(\Mcal) \arrow{r}{} \arrow{d}{} & \Alg_{\Mcal'}(\Mcal'') \times \Algbrd(\Mcal') \arrow{d}{} \\
\Alg_{\Mcal}(\Mcal'') \times \Algbrd(\Mcal) \arrow{r}{} & \Algbrd(\Mcal'').
\end{tikzcd}
\]
This is part of the data that would arise from an enhancement of $\Algbrd(-)$ to a functor of $2$-categories. We do not construct this enhancement here; however note that the above property is already enough to conclude that if $\Mcal, \Mcal'$ are monoidal categories and $F: \Mcal \rightarrow \Mcal'$ is a monoidal functor admitting a (lax monoidal) right adjoint $F^R$, then we have an induced adjunction
\[
F_! : \Algbrd(\Mcal) \rightleftarrows \Algbrd(\Mcal') : (F^R)_!.
\]
\end{remark}

By working with monoidal envelopes and passing to presheaf categories, one can often reduce questions in enriched category theory to the case when the enriching category is a presentable monoidal category. We now study some of the features of this setting.

\begin{remark}
In construction \ref{construction algbrd} we implicitly assumed that all operads and categories were small. Passing to a larger universe, one can similarly discuss categories of algebroids in presentable monoidal categories. Given a presentable monoidal category $\Mcal$, we will denote by $\Algbrd(\Mcal)$ the category of algebroids in $\Mcal$ with a small category of objects. Its version where we allow large categories of objects will be denoted by $\widehat{\Algbrd}(\Mcal)$.
\end{remark}

\begin{proposition}\label{prop pres}
Let $\Mcal$ be a presentable category equipped with a monoidal structure which is compatible with colimits. Then
\begin{enumerate}[\normalfont (i)]
\item The category $\Algbrd(\Mcal)$ is presentable, and the projection $\Algbrd(\Mcal) \rightarrow \Cat$ is a limit and colimit preserving cartesian and cocartesian fibration.
\item For every colimit preserving monoidal functor $F: \Mcal \rightarrow \Mcal'$ into another presentable monoidal category, the induced functor $F_!: \Algbrd(\Mcal) \rightarrow \Algbrd(\Mcal')$ preserves colimits.
\end{enumerate}
\end{proposition}
\begin{proof}
  We note that item (ii) is a direct consequence of remark \ref{some infty two functoriality} together with the adjoint functor theorem. We now prove item (i). It follows from \cite{HA} corollary 3.2.3.5 that for every category $X$ the category $\Algbrd_X(\Mcal)$ is presentable. Moreover, using \cite{HA} corollary 3.1.3.5 we see that  for every functor $i: X \rightarrow Y$ the induced functor $i^!: \Algbrd_Y(\Mcal) \rightarrow \Algbrd_X(\Mcal)$ admits a left adjoint, so that the projection $\Algbrd(\Mcal) \rightarrow \Cat$ is both a cartesian and a cocartesian fibration.   Since the functors  $\Assos_-$ and $\Alg_-(\Mcal)$ are accessible, we conclude from \cite{GHN} theorem 10.3 that $\Algbrd(\Mcal)$ is a presentable category. The fact that the projection to $\Cat$ preserves limits and colimits is now a consequence of \cite{HTT} corollary 4.3.1.11.
  \end{proof}
  
\begin{notation}
For each associative operad $\Mcal$ denote by $\Algbrd(\Mcal)_{\Spc}$ the full subcategory of $\Algbrd(\Mcal)$ on those algebroids which have a space of objects.
\end{notation}

\begin{remark}\label{remark generators cells}
Let $\Mcal$ be a presentable monoidal category. Let $\kappa$ be a regular cardinal  and let   $\lbrace m_i \rbrace_{i \in \Ical}$ be a small family of $\kappa$-compact generators of  $\Mcal$. Then the cells $C_{m_i}$ together with the unit algebroid $1_\Mcal$ are a family of $\kappa$-compact generators of $\Algbrd(\Mcal)_{\Spc}$.
\end{remark}

\subsection{Enriched categories} \label{subsection enriched}

Our next goal is to review the notion of enriched category. In order to do this, we will first need to study the category of algebroids in the case $\Mcal = \Spc$ equipped with its cartesian monoidal structure.

\begin{construction}\label{construction inclusion cat}
Let $\Mcal = \Spc$ be the category of spaces, equipped with its cartesian monoidal structure. Then for every category $X$ the category $\Algbrd_X(\Spc)$ is presentable (\cite{HA} corollary 3.2.3.5), and in particular admits an initial object. Since the projection $\Algbrd(\Spc) \rightarrow \Cat$ is a cartesian fibration, there is a  unique section $s: \Cat \rightarrow \Algbrd(\Spc)$ such that for every category $X$ we have that $s(X)$ is initial in $\Algbrd_X(\Spc)$.

Consider the  cartesian square
\[
\begin{tikzcd}
 \Algbrd(\Spc)_{\Spc} \arrow{r}{i'} \arrow{d}{p'} & \Algbrd(\Spc) \arrow{d}{p} \\
 \Spc \arrow{r}{i} & \Cat.
\end{tikzcd}
\]
Since $p$ is a cartesian fibration and $i$ admits a right adjoint, we have that the above square is horizontally right adjointable. The right adjoint $i'^R: \Algbrd(\Spc) \rightarrow  \Algbrd(\Spc)_{\Spc}$ maps an algebroid $\Acal: \Assos_X \rightarrow \Spc$ to the algebroid defined by the composite map
\[
\Assos_{X^{\leq 0}} \rightarrow \Assos_X \xrightarrow{\Acal} \Spc.
\]
Denote by $\rho: \Cat \rightarrow \Algbrd(\Spc)_{\Spc}$ the composite map $(i')^R s$.
\end{construction}

\begin{example}
The algebroid $\rho([0])$ is the unit algebroid $1_{\Spc}$.
\end{example}

\begin{example}
Examining the description of free algebras from \cite{HA} definition 3.1.3.1 yields the following description of $\rho([1])$:
\begin{itemize}
\item $\rho([1])$ has a set of objects with two elements $0, 1$.
\item $\rho([1])(0, 0) = \rho([1])(1, 1) = \rho([1])(1, 0)$ are the singleton set.
\item $\rho([1])(0, 1)$ is empty.
\end{itemize}
In other words, we have that $\rho([1])$ is equivalent to the cell $C_{[0]}$.
\end{example}

\begin{lemma}\label{lemma left adjoint to s}
The section $s$ from  construction \ref{construction inclusion cat} admits a left adjoint.
\end{lemma}
\begin{proof}
We continue with the notation from construction \ref{construction inclusion cat}. It follows from remark \ref{remark generators cells} that $\Algbrd(\Spc)_{\Spc}$ is generated under colimits by the cell $C_{[0]}$ and the trivial algebroid $1_{\Spc}$. To obtain a set of generators for $\Algbrd(\Spc)$ it suffices to add the algebroid $s([1])$. Since $\Cat$ admits all colimits, in order to show that $s$ has a left adjoint, it suffices to show that for each generator $G$, there is a category $\Ccal$ and a morphism $\eta: G \rightarrow s(\Ccal)$ such that for every category $\Dcal$ the composite map
\[
\Hom_{\Cat}(\Ccal, \Dcal) \xrightarrow{s_*} \Hom_{\Algbrd(\Spc)}(s(\Ccal), s(\Dcal)) \xrightarrow{\eta^*} \Hom_{\Algbrd(\Spc)}(G, s(\Dcal)) 
\]
is an equivalence.  Note that the section $s$ is fully faithful, so the first map in the above composition is always an isomorphism. Since $1_{\Spc}$ and $s([1])$ belong to the image of $s$, the identity maps of $1_{\Spc}$ and $s([1])$ satisfy the desired condition. 

It remains to consider the case of the generator $C_{[0]}$. We take $\Ccal = [1]$, and the map $\eta: C_{[0]} \rightarrow s(\Ccal)$ to be the morphism of algebroids associated to the image of the map $1_{\Spc} \rightarrow s([1]) ( 1, 0)$ induced by the unique arrow $1 \leftarrow 0$ in $[1]$. Let $\Dcal$ be a category and consider the commutative triangle
\[
\begin{tikzcd}
\Hom_{\Cat}(\Ccal, \Dcal) \arrow{dr}{} \arrow{rr}{\eta^* s_*} & & \Hom_{\Algbrd(\Spc)}(C_{[0]}, s(\Dcal)) \arrow{dl}{} \\
& \Dcal^{\leq 0} \times \Dcal^{\leq 0} &
\end{tikzcd}
\]
where the diagonal maps are the source and target maps. In order to show that $\eta^*s_*$ is an equivalence, it suffices to show that it is an equivalence when restricted to the fiber over any point $(x, y)$ in $\Dcal^{\leq 0} \times \Dcal^{\leq 0}$. This restriction recovers the map $\Hom_{\Dcal}(x, y) \rightarrow s(\Dcal)(y, x)$ which assigns to each arrow $y \leftarrow x: \alpha$ in $\Dcal$, the image of the induced map $1_{\Spc} \rightarrow s(\Dcal)(y,x)$. Our claim now follows from the fact that $s(\Dcal)$ is the free $\Assos_{\Dcal}$-algebra in $\Spc$ on the unique algebra over the empty operad, together with the description of free algebras from \cite{HA} definition 3.1.3.1.
\end{proof}

The following proposition is a slight rephrasing of \cite{GH} theorem 4.4.7 and the discussion in \cite{Hinich} section 5.

\begin{proposition}\label{prop compare defs}
There is a unique equivalence between  $\Algbrd(\Spc)_{\Spc}$ and the category $\Pcal(\Delta)_{\text{\normalfont Seg}}$ of Segal spaces which intertwines the map $\rho$ and the canonical inclusion of $\Cat$ into $\Pcal(\Delta)_{\text{\normalfont Seg}}$ as the subcategory of complete Segal spaces. 
\end{proposition}
\begin{proof}
First we note that this equivalence is unique, if it exists. Indeed, the same method of proof of \cite{HTT} theorem 5.2.9.1 shows that the space of automorphisms of the category $\Pcal(\Delta)_{\text{Seg}}$ is a two element set, consisting of the identity and the orientation reversing automorphism. It follows that the space of automorphisms of $\Pcal(\Delta)_{\text{Seg}}$ that restrict to the identity on $\Cat$ is contractible.

The existence of an equivalence $F: \Algbrd(\Spc)_{\Spc} \rightarrow \Pcal(\Delta)_{\text{Seg}} $  is the subject of \cite{GH} theorem 4.4.7. Denote by $i: \Cat \rightarrow \Pcal(\Delta)_{\text{Seg}}$ the inclusion.  It remains to show that we have an equivalence $F\rho = i$. Note that by virtue of lemma \ref{lemma left adjoint to s}, the map $\rho$ admits a left adjoint. It therefore suffices to show that there is an equivalence  $i^L = \rho^L F^{-1}$. 

Both $i^L$ and $\rho^L F^{-1}$ are  colimit preserving functors $\Pcal(\Delta)_{\text{Seg}} \rightarrow \Cat$, and so they are determined by the data of a Segal cosimplicial category. In the case of $i^L$, this is the canonical inclusion $\Delta \rightarrow \Cat$. The proof of lemma \ref{lemma left adjoint to s} shows that $\rho^L$ maps $1_{\Spc}$ to $[0]$ and $C_{[0]}$ to $[1]$, in a way compatible with the source and target maps. Moreover, inspecting the construction of the equivalence $F$ from \cite{GH} reveals that $F^{-1}$ maps $[0]$ to $1_{\Spc}$ and $[1]$ to $C_{[0]}$, in a way compatible with source and target maps. It follows that the Segal cosimplicial category induced by $\rho^L F^{-1}$ is the identity on the full subcategory of $\Delta$ on the objects $[0]$ and $[1]$. The Segal conditions imply that $\rho^L F^{-1}([n])$ is equivalent to $i^L([n])$ for all $[n]$, and it is in particular a $0$-truncated category.  Our claim now follows from the fact that the source and target maps $[0] \rightarrow [1] \leftarrow [0]$ are jointly surjective, using  the same arguments as those that establish corollary \ref{coro uniqueness extended}
\end{proof}

\begin{remark}\label{remark op y compare defs}
It follows from remark \ref{remark op algebroids} that $\Algbrd(\Spc)_{\Spc}$ comes equipped with an involution $(-)^\op$. The map $\rho$ intertwines the involutions $(-)^\op$ on $\Cat$ and $\Algbrd(\Spc)_{\Spc}$. It follows from the uniqueness statement in proposition  \ref{prop compare defs} that the equivalence between $\Pcal(\Delta)_{\text{Seg}}$ and $\Algbrd(\Spc)_{\Spc}$ admits a $\ZZ/2\ZZ$-equivariant structure.
\end{remark}

We now review the definition of $\Mcal$-enriched categories. These are $\Mcal$-algebroids satisfying a suitable completeness condition.

\begin{notation}
Let $\Mcal$ be an associative operad and equip $\Spc$  with its cartesian monoidal structure. We denote by $\tau_\Mcal: \Mcal \rightarrow \Spc$ the morphism of associative operads which maps each object $m$ in $\Mcal$ to the space of operations from the empty list into $m$.
\end{notation}

\begin{remark}
Let $\Mcal$ be monoidal category. If $\Mcal$ is presentable monoidal then the lax symmetric monoidal functor $\tau_\Mcal: \Mcal \rightarrow \Spc$  is right adjoint to the unit map $\Spc \rightarrow \Mcal$. In general, $\tau_\Mcal$ can be obtained as the composition of the symmetric monoidal embedding $\Mcal \rightarrow \Pcal(\Mcal)$ together with the lax symmetric monoidal map $\tau_{\Pcal(\Mcal)} : \Pcal(\Mcal) \rightarrow \Spc$.
\end{remark}

\begin{definition}
Let $\Mcal$ be an associative operad. An object $\Acal$ in $\Algbrd(\Mcal)$ is said to be an enriched category if it has a space of objects, and the induced object $(\tau_\Mcal)_!\Acal$ in $\Algbrd(\Spc)_{\Spc}$ belongs to the image of $\rho$. We denote by $\Cat^\Mcal$ the full subcategory of $\Algbrd(\Mcal)$ on the enriched categories. Given an $\Mcal$-enriched category $\Acal$ and a pair of objects $x, y$ in $\Acal$, we will usually use the notation $\Hom_\Acal(x, y)$ instead of $\Acal(y, x)$.
\end{definition}

In other words, an $\Mcal$-enriched category is a $\Mcal$-algebroid whose underlying Segal space is a complete Segal space.

\begin{example}
Let $\Mcal$ be a monoidal category such that the  monoid $\End_\Mcal(1_\Mcal)$ does not  have nontrivial invertible elements (for instance, if $\Mcal$ is a cartesian closed presentable category).  Then the unit algebroid $1_\Mcal$ is an $\Mcal$-enriched category. If in addition $\Mcal$ has an initial object which is compatible with the monoidal structure, and the space of maps from the unit to the initial object is empty, then for every $m$ in $\Mcal$ the cell $C_m$ from example \ref{example cells} is an $\Mcal$-enriched category. 
\end{example}

\begin{remark}\label{remark op enriched}
It follows from remark \ref{remark op y compare defs} that an algebroid $\Acal$ is an enriched category if and only if $\Acal^\op$ is an enriched category. In other words, the involution $(-)^\op$ restricts to an involution on the full subcategory of $\Algbrd$ on the enriched categories.
\end{remark}

\begin{proposition}\label{prop properties enriched pres}
Let $\Mcal$ be a presentable monoidal category. Then
\begin{enumerate}[\normalfont (i)]
\item The inclusion $\Cat^\Mcal \rightarrow \Algbrd(\Mcal)_{\Spc}$ exhibits $\Cat^\Mcal$ as an accessible localization of $\Algbrd(\Mcal)_{\Spc}$. In particular, $\Cat^\Mcal$ is presentable.
\item Let $F: \Mcal \rightarrow \Mcal'$ be a colimit preserving monoidal functor into another presentable monoidal category. Then the functor $F_!: \Algbrd(\Mcal)_{\Spc} \rightarrow \Algbrd(\Mcal')_{\Spc}$ descends to a functor $\Cat^\Mcal \rightarrow \Cat^{\Mcal'}$.
\end{enumerate}
\begin{proof}
Recall that $\Cat$ embeds into the category of Segal spaces as the full subcategory of local objects for the projection $\alpha$ from the walking isomorphism to the terminal category.  Since the lax monoidal functor $\tau_\Mcal$ is right adjont to the unit map $1_\Mcal: \Spc \rightarrow \Mcal$, we obtain an adjunction
\[
\begin{tikzcd}
(1_\Mcal)_!: \Algbrd(\Spc)_{\Spc} \arrow[r, shift left=0.25em] &  \arrow[l, shift left = 0.25em] \Algbrd(\Mcal)_{\Spc} : (\tau_{\Mcal})_!
\end{tikzcd}.
\]
It follows that $\Cat^\Mcal$ is the full subcategory of $\Algbrd(\Mcal)_{\Spc}$ of $i_! \alpha$ local objects, which proves item (i). Item (ii) now follows from the fact that $F_! (1_\Mcal)_! \alpha$ is equivalent to $(1_{\Mcal'})_! \alpha$, which becomes an isomorphism upon projection to $\Cat^{\Mcal'}$. 
\end{proof}
\end{proposition}

\begin{remark}
Let $\Mcal$ be a monoidal category. Then $\Algbrd(\Mcal)_{\Spc}$ is a full subcategory of $\Algbrd(\Pcal(\Mcal))_{\Spc}$, and moreover $\Cat^{\Mcal}$ is the intersection of $\Algbrd(\Mcal)_{\Spc}$ with $\Cat^{\Pcal(\Mcal)}$. Let $\Ccal$ be an object in $\Algbrd(\Mcal)_{\Spc}$ and let $\Ccal'$ be its image in $\Cat^{\Pcal(\Mcal)}$. It follows from the description of local equivalences from \cite{GH} corollary 5.6.3 that $\Ccal'$ belongs to $\Cat^\Mcal$. It follows that $\Cat^\Mcal$ is a localization of $\Algbrd(\Mcal)_{\Spc}$, and moreover a map $\Acal \rightarrow \Acal'$ in $\Algbrd(\Mcal)_{\Spc}$ is local if and only if it is fully faithful (i.e., cartesian for the projection $\Algbrd(\Mcal)_{\Spc} \rightarrow \Cat$) and surjective on objects. 
\end{remark}

\begin{example}\label{remark unit inclusion}
Let $\Mcal$ be a presentable monoidal category. As a consequence of proposition \ref{prop properties enriched pres} the unit map $1_\Mcal: \Spc \rightarrow \Mcal$ induces a functor $(1_\Mcal)_!: \Cat = \Cat^{\Spc} \rightarrow \Cat^{\Mcal}$. In other words, any category defines an $\Mcal$-enriched category. 
\end{example}

\begin{remark}\label{remark shriek fully faithful}
Let $i : \Mcal \rightarrow \Mcal'$ be  a colimit preserving monoidal functor between presentable monoidal categories. Assume that $i$ is fully faithful, so that the functor $i_!: \Algbrd(\Mcal) \rightarrow \Algbrd(\Mcal')$ is fully faithful. Then for every $\Mcal$-algebroid $\Acal$ with a space of objects, the Segal space underlying $i_! \Acal$ is equivalent to the Segal space underlying $\Acal$. In particular, $\Acal$ is an $\Mcal$-enriched category if and only if $i_! \Acal$ is an $\Mcal'$-enriched category. This implies that the commutative square
\[
\begin{tikzcd}
\Algbrd(\Mcal)_{\Spc} \arrow{d}{i_!} \arrow{r}{} & \Cat^{\Mcal} \arrow{d}{i_!} \\
\Algbrd(\Mcal')_{\Spc} \arrow{r}{} & \Cat^{\Mcal'}
\end{tikzcd}
\]
arising from proposition \ref{prop properties enriched pres} item (ii), is horizontally right adjointable.

Assume now that $i$ admits a left adjoint, so that $\Mcal$ is a localization of $\Mcal'$. Then for every space $X$ the functor $i_!: \Algbrd_X(\Mcal) \rightarrow \Algbrd_X(\Mcal')$ preserves limits and is accessible. It follows from \cite{HTT} proposition 4.3.1.9 together with the fact that the projection $\Algbrd(\Mcal')_{\Spc} \rightarrow \Spc$ is both a cartesian and a cocartesian fibration, that the functor $i_!: \Algbrd(\Mcal) \rightarrow \Algbrd(\Mcal')$ preserves limits and is accessible. It now follows from the adjoint functor theorem that the above square is in fact also vertically left adjointable. In particular, we have that $\Cat^{\Mcal}$ is a localization of $\Cat^{\Mcal'}$.

We can describe this in more concrete terms. Let $\Acal$ be an $\Mcal'$-enriched category. Then $\Acal$ belongs to the image of $i_!$ if and only if for every pair of objects $x, y $ in $\Acal$ we have that $\Hom_\Acal(x, y)$ belongs to $\Mcal$. Equivalently, for every object $m'$ in $\Mcal'$, the map
\[
\Hom_{\Mcal'}(ii^L m' , \Hom_\Acal(x, y)) \rightarrow \Hom_{\Mcal'}(m', \Hom_{\Acal}(x, y))
\]
given by precomposition with the unit $m' \rightarrow i i^L m'$, is an equivalence. It follows that $\Acal$ belongs to the image of $i_!$ if and only if it is local for the class of maps $C_{m'} \rightarrow C_{i i^L m'}$. Note that we can simplify this further: we may take $m'$ to belong to a set of generators of $\Mcal'$.
\end{remark}

\subsection{Multiplicativity}\label{subsection multiplicativity}

We now discuss the notion of tensor product of algebroids and enriched categories.

\begin{proposition}\label{prop products of algebroids}
The category $\Algbrd$ has finite products, which are preserved by the projection to $\Cat \times \Op_{\Assos}$.
\end{proposition}
\begin{proof}
Let $X$ be a category. Then the functor $\Algbrd_X(-): \Op_{\Assos} \rightarrow \Cat$ is limit preserving. It follows from \cite{HTT} corollary 4.3.1.15, that the total category $\Algbrd_X$ of the associated cocartesian fibration has all finite products, which are preserved by the projection to $\Op_{\Assos}$. Furthermore, if  $\lbrace \Acal_i \rbrace_{i \in \Ical}$ is a finite family of objects of $\Algbrd_X$ lying above a finite family of associative operads $\lbrace \Mcal_i \rbrace_{i \in \Ical}$, then its product is the unique object $\Acal$ in $\Algbrd_X(\prod \Mcal_i)$ equipped with cocartesian  arrows to $\Acal_i$ lifting the projection $\prod \Mcal_i \rightarrow \Mcal_i$, for all $i$ in $\Ical$.

Assume now given a functor of categories $f: Y \rightarrow X$. Then $f^! \Acal$ is an object in $\Algbrd_Y(\prod \Mcal_i)$ which comes equipped with cocartesian arrows to $f^! \Acal_i$ lifting the projections $\prod\Mcal_i \rightarrow \Mcal_i$, for all $i$ in $\Ical$. It follows that $f^!$ preserves finite products. By a combination of  \cite{HTT} propositions 4.3.1.9 and 4.3.1.10 we conclude that the projection $\Algbrd \rightarrow \Cat$ has all relative finite products, which are preserved by the map $\Algbrd \rightarrow \Cat \times \Op_{\Assos}$. Our result now follows from the fact that $\Cat$ has all finite products.
\end{proof}

\begin{remark}
It follows from proposition \ref{prop products of algebroids} that the final object of $\Algbrd$ is the unique algebroid lying above the final object of $\Cat \times \Op_{\Assos}$. In other words, this is the unit algebroid of the final monoidal category.
\end{remark}

\begin{notation}
Let $\Acal, \Bcal$ be objects of $\Algbrd$. We denote by $\Acal \boxtimes \Bcal$ their product in $\Algbrd$.
\end{notation}

\begin{remark}\label{remark fact proj}
Let $X, Y$ be categories and $\Mcal, \Ncal$ be associative operads. Let $\Acal, \Bcal$ be objects in $\Algbrd_X(\Mcal)$ and $\Algbrd_Y(\Ncal)$, respectively. Denote by $p_1, p_2$ the projections from $\Algbrd$ to $\Cat$ and $\Op_{\Assos}$, respectively. It follows from the proof of proposition \ref{prop products of algebroids} that the span
\[
\Acal \leftarrow \Acal \boxtimes \Bcal \rightarrow \Bcal
\] is  the unique lift of the span
\[
(X, \Mcal) \leftarrow (X \times Y , \Mcal \times \Ncal) \rightarrow (Y, \Ncal)
\]
such that its left and right legs  can be written as the composition of a $p_2$-cocartesian followed by a $p_1$-cartesian morphism.  

It follows that $\Acal \boxtimes \Bcal$ is the algebroid defined by the map
\[
\Assos_{X \times Y} = \Assos_{X} \times \Assos_Y \xrightarrow{\Acal \times \Bcal} \Mcal \times \Ncal
\]
and the projections to $\Acal$ and $\Bcal$ are induced from the following commutative diagram:
\[
\begin{tikzcd}
\Assos_X \arrow{d}{\Acal} & \Assos_X \times \Assos_Y \arrow{d}{\Acal \times \Bcal} \arrow{l}{} \arrow{r}{} & \Assos_Y \arrow{d}{\Bcal} \\
\Mcal & \Mcal \times \Ncal \arrow{r}{} \arrow{l}{} & \Ncal
\end{tikzcd}
\]

In particular, $\Acal \boxtimes \Bcal$ is an $(\Mcal \times \Mcal')$-algebroid with category of objects $X \times Y$, and for every pair of objects $(x', y'), (x, y)$ we have an equivalence
\[
(\Acal \boxtimes \Bcal)((x', y'), (x, y)) = (\Acal(x' , x) , \Bcal(y', y)).
\]
The composition maps for $\Acal \boxtimes \Bcal$ are obtained by taking the product of the composition maps of $\Acal$ and $\Bcal$. 
\end{remark}

\begin{proposition}\label{prop products compatible with }
Let $f: \Acal \rightarrow \Acal'$ be a morphism in $\Algbrd$ and let $\Bcal$ be another object of $\Algbrd$.  Denote by $p = (p_1, p_2)$ the projection $\Algbrd \rightarrow \Cat \times \Op_{\Assos}$.
\begin{enumerate}[\normalfont (i)]
\item If $f$ is $p_1$-cartesian then $f \boxtimes \id_{\Bcal}$ is $p_1$-cartesian.
\item If $f$ is $p_2$-cocartesian then $f \boxtimes \id_{\Bcal}$ is $p_2$-cocartesian.
\end{enumerate}
\end{proposition}
\begin{proof}
Denote by $X, X', Y$ the categories of objects of $\Acal, \Acal'$ and $\Bcal$, respectively, and let $\Mcal, \Mcal', \Ncal$ be their underlying associative operads.  Consider the following commutative diagram in $\Cat \times \Op_{\Assos}$:
\[
\begin{tikzcd}[column sep = huge]
(X \times Y, \Mcal \times \Ncal) \arrow{d}{(p_1 f \times \id, p_2 f \times \id)} \arrow{r}{(\id, p_{\Mcal})} & (X \times Y, \Mcal) \arrow{d}{(p_1 f \times \id, p_2 f)} \arrow{r}{(p_X, \id_{\Mcal}) } & (X, \Mcal) \arrow{d}{(p_1 f, p_2 f)}\\
(X' \times Y, \Mcal' \times \Ncal) \arrow{r}{(\id, p_{\Mcal'})} & (X' \times Y, \Mcal') \arrow{r}{(p_{X'}, \id)} & (X', \Mcal') 
\end{tikzcd}
\]
This admits a lift to a commutative diagram
\[
\begin{tikzcd}
\Acal \boxtimes \Bcal \arrow{d}{f \boxtimes \id_{\Bcal}} \arrow{r}{} & \arrow{d}{\mu} \overline{\Acal}  \arrow{r}{} &  \Acal \arrow{d}{f} \\
\Acal' \boxtimes \Bcal \arrow{r}{} &\overline{\Acal'} \arrow{r}{} & \Acal'
\end{tikzcd}
\]
where the horizontal rows are the factorizations of the projections as $p_2$-cocartesian maps followed by $p_1$-cartesian maps.

Similarly, the commutative diagram
\[
\begin{tikzcd}[column sep = huge]
(X \times Y, \Mcal \times \Ncal) \arrow{d}{(p_1 f \times \id, p_2 f \times \id)} \arrow{r}{(\id, p_{\Ncal})} & (X \times Y, \Ncal) \arrow{d}{(p_1 f \times \id, \id)} \arrow{r}{(p_Y, \id) } & (Y, \Ncal) \arrow{d}{(\id, \id)}\\
(X' \times Y, \Mcal' \times \Ncal) \arrow{r}{(\id, p_{\Ncal})} & (X' \times Y, \Ncal) \arrow{r}{(p_{Y}, \id)} & (Y, \Ncal) 
\end{tikzcd}
\]
admits a lift to a commutative diagram
\[
\begin{tikzcd}
\Acal \boxtimes \Bcal \arrow{d}{f \boxtimes \id_{\Bcal}} \arrow{r}{} & \arrow{d}{\nu} \overline{\Bcal}  \arrow{r}{} &  \Bcal \arrow{d}{\id} \\
\Acal' \boxtimes \Bcal \arrow{r}{} & \widetilde{\Bcal} \arrow{r}{} & \Bcal
\end{tikzcd}
\]
where the horizontal rows consist of a  $p_2$-cocartesian  followed by a $p_1$-cartesian map.

Assume now that $f$ is $p_1$-cartesian, so that $\Mcal = \Mcal'$. Then $\mu$ and $\nu$ are $p_1$-cartesian. Write $f \boxtimes \id_{\Bcal} = \alpha \eta $ where $\alpha$ is $p_1$-cartesian and $\eta$ is such that $(p_1, p_2)\eta$ is invertible. We have 
\[
\mu = (p_{\Mcal})_! (f \boxtimes \id_{\Bcal}) = (p_{\Mcal})_!(\alpha ) (p_{\Mcal})_!(\eta).
\]
Since $(p_{\Mcal})_!$ is a morphism of cartesian fibrations, we have that $(p_{\Mcal})_!(\alpha)$ is cartesian and therefore $(p_{\Mcal})_!(\eta)$ is an isomorphism. Similarly, we can conclude that $(p_{\Ncal})_! (\eta)$ is an isomorphism. Item (i) now follows from the fact that the projections 
\[
\Algbrd_{X \times Y}( \Ncal) \xleftarrow{ (p_{\Ncal})_! } \Algbrd_{X \times Y}(\Mcal \times \Ncal) \xrightarrow{(p_{\Mcal})_!} \Algbrd_{X \times Y}(\Mcal)
\]
are jointly conservative.

We now prove item (ii). In this case, $f$ is $p_2$-cocartesian, so that $X = X'$. We therefore have that $\nu$ is invertible. Furthermore, we have that $\mu = p_X^! f$ is $p_2$-cocartesian. As before, write $f \boxtimes \id_{\Bcal} = \eta \alpha$ where $\alpha$ is $p_2$-cocartesian and $\eta$ is such that $(p_1, p_2)\eta$ is invertible. The composition of $\eta$ with the $p_2$-cocartesian map $p_{\overline{\Acal'}}: \Acal' \boxtimes \Bcal \rightarrow \overline{\Acal'}$ is a lift of the projection $(\id_{X \times Y}, p_{\Mcal'})$ whose composition with the $p_2$-cocartesian map $\alpha$ is $p_2$-cocartesian. It follows that $p_{\overline{\Acal'}} \eta$ is $p_2$-cocartesian, and therefore we have that  $(p_{\Mcal'})_! \eta$ is an isomorphism. A similar argument shows that $(p_{\Ncal})_! \eta$ is an isomorphism. Our result now follows from the fact that 
the projections 
\[
\Algbrd_{X \times Y}( \Ncal) \xleftarrow{ (p_{\Ncal})_! } \Algbrd_{X \times Y}(\Mcal' \times \Ncal) \xrightarrow{(p_{\Mcal'})_!} \Algbrd_{X \times Y}(\Mcal')
\]
are jointly conservative.
\end{proof}

\begin{construction} \label{construction sym mon structure}
We equip $\Algbrd$ and $\Cat \times \Op_{\Assos}$ with their cartesian symmetric monoidal structures, so that the projection $\Algbrd \rightarrow \Cat \times \Op_{\Assos}$ inherits a canonical symmetric monoidal structure by proposition \ref{prop products of algebroids}. It follows from proposition \ref{prop products compatible with } that the projection $\Algbrd \rightarrow \Op_{\Assos}$ is a cocartesian fibration of operads, which straightens to a lax symmetric monoidal structure on the functor $\Algbrd(-): \Op_{\Assos} \rightarrow \widehat{\Cat}$. Given $\Mcal$ and $\Ncal$ two associative operads, this produces a functor 
\[
\Algbrd(\Mcal) \times \Algbrd(\Ncal) \rightarrow \Algbrd(\Mcal \times \Ncal)
\]
 which sends a pair of algebroids $\Acal, \Bcal$ to $\Acal \boxtimes \Bcal$. 

Let $\Mcal$ be a symmetric monoidal category. We can think about $\Mcal$ as a commutative algebra object in $\Alg(\Cat)$, and hence as a commutative algebra object in $\Op_{\Assos}$. It follows that $\Algbrd(\Mcal)$ inherits a symmetric monoidal structure. We will usually denote by 
\[
\otimes : \Algbrd(\Mcal) \times \Algbrd(\Mcal) \rightarrow \Algbrd(\Mcal)
\]
 the resulting functor. We note that the assignment $\Mcal \mapsto (\Algbrd(\Mcal), \otimes)$ is part of a functor $\CAlg(\Cat) \rightarrow \CAlg(\widehat{\Cat})$.
\end{construction}

\begin{remark}\label{remark sym mon str}
Let $\Mcal$ be a symmetric monoidal category. The unit of the symmetric monoidal structure on $\Algbrd(\Mcal)$ is the unit algebra in $\Mcal$, thought of as an algebroid with category of objects $[0]$. The tensor product functor on $\Algbrd(\Mcal)$ can be computed as the composition
\[
\Algbrd(\Mcal) \times \Algbrd(\Mcal)\xrightarrow{\boxtimes} \Algbrd(\Mcal \times \Mcal) \xrightarrow{m_!} \Mcal
\]
where  $m: \Mcal \times \Mcal \rightarrow \Mcal$ is the tensoring map. In particular, if $\Acal$ and $\Bcal$ have category of objects $X$ and $Y$ respectively then $\Acal \otimes \Bcal$ has category of objects $X \times Y$. Moreover, if $x, x'$ are objects in $X$ and $y, y'$ are objects in $Y$, we have an equivalence 
\[
(\Acal \otimes \Bcal)( (x', y'), (x, y) ) = \Acal(x', x) \otimes \Acal(y', y).
\]
\end{remark}

\begin{proposition}\label{prop is cartesian}
Let $\Mcal$ be a category admitting finite products, equipped  with the cartesian symmetric monoidal structure. Then the symmetric monoidal structure on $\Algbrd(\Mcal)$ given by construction \ref{construction sym mon structure} is cartesian.
\end{proposition}
\begin{proof}
As observed in remark \ref{remark sym mon str}, the unit $1_{\Algbrd(\Mcal)}$ of $\Algbrd(\Mcal)$ is the unit algebra in $\Mcal$. To check that $1_{\Algbrd(\Mcal)}$ is final in $\Algbrd(\Mcal)$ we have to see that for every category $X$, the algebroid $\pi^!_X1_{\Algbrd(\Mcal)}$ is final in $\Algbrd_{X}(\Mcal)$, where $\pi_X: X \rightarrow [0]$ denotes the projection. Indeed, for every pair of objects $x, y$ in $X$ we have 
\[
\pi^!_X1_{\Algbrd(\Mcal)} ( y, x) = 1_{\Algbrd(\Mcal)}(\pi_X y, \pi_X x) = 1_\Mcal
\]
which is final in $\Mcal$. The fact that $\pi_X^!1_{\Algbrd(\Mcal)}$ is final then follows from \cite{HA} corollary 3.2.2.5.

It remains to check that for every pair of algebroids $\Acal, \Bcal$ in $\Mcal$, the projections
\[
\Acal = \Acal \otimes 1_{\Algbrd(\Mcal)} \leftarrow \Acal \otimes \Bcal \rightarrow 1_{\Algbrd(\Mcal)} \otimes \Bcal = \Bcal
\]
exhibit $\Acal \otimes \Bcal$ as the product of $\Acal$ and $\Bcal$ in $\Algbrd(\Mcal)$. Let $X, Y$ be the category of objects of $\Acal, \Bcal$ respectively. We have to show that for every category $Z$ equipped with functors $j: Z \rightarrow X$ and $j': Z \rightarrow Y$, the projections
\[
j^!\Acal \leftarrow (j \times j')^! (\Acal \otimes \Bcal) \rightarrow j'^! \Bcal
\]
exhibit $(j \times j')^! (\Acal \otimes \Bcal)$ as the product of $j^! \Acal$ and $j'^!\Bcal$ in the category $\Algbrd_Z(\Mcal)$. Let $z, w$ be objects in $Z$. The induced diagram
\[
j^!\Acal(z, w) \leftarrow (j \times j')^!(\Acal \otimes \Bcal)(z, w) \rightarrow j'^!\Bcal(z, w)
\]
is the equivalent to the diagram
\[
\Acal(jz, jw) = \Acal(jz, jw) \otimes 1_\Mcal \leftarrow \Acal(jz, jw) \otimes \Bcal(j'z, j'w) \rightarrow 1_\Mcal  \otimes \Bcal(j'z, j'w) = \Bcal(j'z, j'w)
\]
and therefore it exhibits $(j \times j')^!(\Acal \otimes \Bcal)(z, w)$ as the product of $j^!\Acal(z, w)$ and $j'^!\Bcal(z, w)$. Our result now follows from another application of \cite{HA} corollary 3.2.2.5.
\end{proof}

The symmetric monoidal structure on algebroids from construction \ref{construction sym mon structure} restricts to algebroids with a space of objects. The next proposition shows that this induces a symmetric monoidal structure on enriched categories.

\begin{proposition}\label{prop Cat mcal compat}
Let $\Mcal$ be a symmetric monoidal category. Then the localization functor $\Algbrd(\Mcal)_{\Spc} \rightarrow \Cat^\Mcal$ is compatible with the restriction of the symmetric monoidal structure of construction \ref{construction sym mon structure}.
\end{proposition}
\begin{proof}
Recall  that a morphism $F: \Acal \rightarrow \Bcal$ in $\Algbrd(\Mcal)_{\Spc}$ is local for the localization in  the statement if and only if it is fully faithful and  surjective on objects. Equivalently, this means that $F$ is $p_1$-cartesian and surjective on objects. 

Let $\Acal'$ be another object of $\Algbrd(\Mcal)_{\Spc}$. It follows from propositions \ref{prop products of algebroids} and \ref{prop products compatible with } that 
\[
F \boxtimes \id_{\Acal'} : \Acal \boxtimes \Acal' \rightarrow \Bcal \boxtimes \Acal'
\] is still fully faithful and surjective on objects. Therefore the map $F \otimes \id_{\Acal'} = m_! (F \boxtimes \id_{\Acal'})$ is also fully faithful and surjective on objects, so it is local for the localization in the statement, as desired.
\end{proof}

\begin{corollary}\label{coro enrich cat is sm}
Let $\Mcal$ be a  symmetric monoidal category. Then $\Cat^\Mcal$ inherits a symmetric monoidal structure from $\Algbrd(\Mcal)_{\Spc}$, and the localization $\Algbrd(\Mcal)_{\Spc} \rightarrow \Cat^\Mcal$ has a canonical symmetric monoidal structure.
\end{corollary}

\begin{example}\label{example cartesian closed}
Let $\Mcal$ be a category admitting finite products, equipped with its cartesian symmetric monoidal structure. Then it follows from proposition \ref{prop is cartesian} that the induced symmetric monoidal structure on $\Cat^{\Mcal}$ is cartesian.
\end{example}	

For later purposes, we will need a generalization of the functoriality of construction \ref{construction sym mon structure} which deals with lax symmetric monoidal functors between symmetric monoidal categories. In fact, it turns out that for any symmetric operad $\Mcal$ one can give $\Algbrd(\Mcal)$ and $\Cat^{\Mcal}$ the structure of a symmetric operad, in a way that depends functorially on $\Mcal$.

\begin{construction}\label{construction operadic structure}
Denote by $\Env: \Op \rightarrow \CAlg(\Cat)$ the functor that sends each symmetric operad to its enveloping symmetric monoidal category - in other words, this is left adjoint to the inclusion $\CAlg(\Cat) \rightarrow \Op$. Consider now the composite functor
\[
\xi: \Op \xrightarrow{\Env} \CAlg(\Cat) \xrightarrow{\Algbrd(-)} \CAlg(\widehat{\Cat}).
\]

Note that the composition of $\xi$ with the forgetful functor $\CAlg(\widehat{\Cat}) \rightarrow \widehat{\Cat}$ receives a natural transformation $\eta$ from the functor $\Algbrd(-)|_{\Op}: \Op \rightarrow \widehat{\Cat}$. For each symmetric operad $\Mcal$, this induces a functor 
\[
\eta(\Mcal) : \Algbrd(\Mcal) \rightarrow \xi(\Mcal) = \Algbrd(\Env(\Mcal)).
\]
Since the unit map $\Mcal \rightarrow \Env(\Mcal)$ is an inclusion of symmetric operads, we have that $\eta(\Mcal)$ is fully faithful. Its image consists of those $\Env(\Mcal)$-algebroids $\Acal$ such that $\Acal(y, x)$ belongs to $\Mcal$ for each pair of objects $y, x$ in $\Acal$. Since $\Algbrd(\Env(\Mcal))$ has a symmetric monoidal structure, the full subcategory $\Algbrd(\Mcal)$ inherits the structure of a symmetric operad. This is compatible with morphisms of symmetric operads, so we obtain a lift of $\Algbrd(-)|_{\Op}$ to a functor
\[
(\Algbrd(-)|_{\Op})^{\enh}: \Op \rightarrow \widehat{\Op}.
\]
\end{construction}

The following proposition shows that construction \ref{construction operadic structure} extends the functoriality of the theory of algebroids on symmetric monoidal categories from construction \ref{construction sym mon structure}.

\begin{proposition}
The restriction of the functor $(\Algbrd(-)|_{\Op})^{\enh}$  to $\CAlg(\Cat)$ factors through $\CAlg(\widehat{\Cat})$, and coincides with the functor arising from  construction \ref{construction sym mon structure}.
\end{proposition}
\begin{proof}
Let $\Mcal$ be a symmetric monoidal category. Then the inclusion $\Mcal \rightarrow \Env(\Mcal)$ exhibits $\Mcal$ as a symmetric monoidal localization of $\Env(\Mcal)$. It follows that the inclusion $\Algbrd(\Mcal) \rightarrow \Env(\Mcal)$ exhibits $\Algbrd(\Mcal)$ (with its symmetric monoidal structure from construction \ref{construction sym mon structure}) as a symmetric monoidal localization of $\Algbrd(\Env(\Mcal))$.  This shows that the operadic structure on $\Algbrd(\Mcal)$ from construction \ref{construction operadic structure} coincides with the operadic structure underlying the symmetric monoidal structure given to in construction \ref{construction sym mon structure}.

Assume now given a symmetric monoidal functor $F: \Mcal \rightarrow \Mcal'$ between symmetric monoidal categories. We have a commutative square of symmetric monoidal categories and symmetric monoidal functors
\[
\begin{tikzcd}
\Algbrd(\Env(\Mcal)) \arrow{r}{\Env(F)_!} \arrow{d}{} & \Algbrd(\Env(\Mcal')) \arrow{d}{} \\
\Algbrd(\Mcal) \arrow{r}{F_!} & \Algbrd(\Mcal').
\end{tikzcd}
\]
This is vertically right adjointable. Passing to right adjoints of the vertical arrows yields a commutative diagram of symmetric monoidal categories and lax symmetric monoidal functors
\[
\begin{tikzcd}
\Algbrd(\Env(\Mcal)) \arrow{r}{\Env(F)_!}  & \Algbrd(\Env(\Mcal'))  \\ \arrow{u}{}
\Algbrd(\Mcal) \arrow{r}{F_!} & \Algbrd(\Mcal') \arrow{u}{}.
\end{tikzcd}
\]
It follows that the structures of morphism of symmetric operads on $F_!$ arising from constructions \ref{construction sym mon structure} and \ref{construction operadic structure} agree. In particular, we conclude that the restriction of $(\Algbrd(-)|_{\Op})^{\enh}$  to $\CAlg(\Cat)$ factors through $\CAlg(\widehat{\Cat})$.

Consider now the lax commutative triangle
\[
\begin{tikzcd}[row sep = 3em]
\CAlg(\Cat) \arrow{r}{\Env} \arrow{dr}[swap, name = diag]{\Algbrd(-)} & \CAlg(\Cat) \arrow{d}{\Algbrd(-)} \arrow[Rightarrow, to=diag, shorten <>=8pt] \\
& \CAlg(\widehat{\Cat})
\end{tikzcd}
\]
obtained by applying $\Algbrd(-)$ to the counit of the defining adjunction for $\Env$. Passing to right adjoints yields a commutative triangle
\[
\begin{tikzcd}[row sep = 3em]
\CAlg(\Cat) \arrow{r}{\Env} \arrow{dr}[swap, name = diag]{\Algbrd(-)} & \CAlg(\Cat) \arrow{d}{\Algbrd(-)} \arrow[Rightarrow, from=diag, shorten <>=8pt] \\
& \widehat{\Op}.
\end{tikzcd}
\]
This identifies the diagonal arrow with the restriction of $(\Algbrd(-)|_{\Op})^{\enh}$  to $\CAlg(\Cat)$.
\end{proof}

\begin{remark}
Let $F: \Mcal \rightarrow \Mcal'$ be a symmetric monoidal functor between symmetric monoidal categories. Assume that $F$ admits a right adjoint $F^R$, and equip $F^R$ with its natural lax symmetric monoidal structure. Then $F$ and $F^R$ induce a symmetric monoidal adjunction
\[
\begin{tikzcd}
\Env(F) : \Env(\Mcal) \arrow[r, shift left = 0.25em] & \arrow[l, shift left = 0.25em] \Env(\Mcal') : \Env(F^R).
\end{tikzcd}
\]
This in turn induces a symmetric monoidal adjunction
\[
\begin{tikzcd}
\Env(F)_! : \Algbrd(\Env(\Mcal)) \arrow[r, shift left = 0.25em] & \arrow[l, shift left = 0.25em] \Algbrd(\Env(\Mcal')) :\Env(F^R)_!
\end{tikzcd}
\]
which restricts to an adjunction with symmetric monoidal left adjoint
\[
\begin{tikzcd}
F_! : \Algbrd(\Mcal) \arrow[r, shift left = 0.25em] & \arrow[l, shift left = 0.25em] \Algbrd(\Mcal') : F^R_!.
\end{tikzcd}
\]
It follows from this that the lax symmetric monoidal structure on $F^R_!$ arising from construction \ref{construction operadic structure} is equivalent to the one arising by passing to adjoints the  symmetric monoidal structure on $F_!$.
\end{remark}

\begin{remark}
Let $F: \Mcal \rightarrow \Mcal'$ be a lax symmetric monoidal functor between symmetric monoidal categories.  Then the lax symmetric monoidal functor
\[
F_! : \Algbrd(\Mcal) \rightarrow \Algbrd(\Mcal')
\]
restricts to a lax symmetric monoidal functor
\[
F_! : \Algbrd(\Mcal)_{\Spc} \rightarrow \Algbrd(\Mcal')_{\Spc}
\]
which in turn induces a lax symmetric monoidal functor $F_! : \Cat^\Mcal \rightarrow \Cat^{\Mcal'}$. This forms part of a functor
\[
(\Cat^{(-)}|_{\CAlg(\Cat)^{\text{lax}}})^{\enh}: \CAlg(\Cat)^{\text{lax}} \rightarrow \CAlg(\widehat{\Cat})^{\text{lax}}
\]
where $\CAlg(\Cat)^{\text{lax}}$ denotes the category of symmetric monoidal categories and lax symmetric monoidal functors.
\end{remark}

We now study the  presentable symmetric monoidal case. The following result follows from a version of \cite{GH} corollary 4.3.16 - we refer the reader there for a proof.
\begin{proposition} \label{prop algbrd is pres}
Let $\Mcal$ be a presentable symmetric monoidal category. Then the induced symmetric monoidal structures on $\Algbrd(\Mcal)$ and $\Cat^\Mcal$ are compatible with colimits.
\end{proposition}

\begin{corollary}\label{coro Algbrd y cat cart closed}
Let $\Mcal$ be a cartesian closed presentable category. Then $\Algbrd(\Mcal)$ and $\Cat^\Mcal$ are cartesian closed.
\end{corollary}
\begin{proof}
Combine proposition \ref{prop algbrd is pres} with example \ref{example cartesian closed}.
\end{proof}

\begin{notation}\label{notation funct}
Let $\Mcal$ be a presentable symmetric monoidal category. We denote by 
\[
\Funct(-, -): \Algbrd(\Mcal)^\op \times \Algbrd(\Mcal) \rightarrow \Algbrd(\Mcal)
\]
the internal Hom for the closed symmetric monoidal category $\Algbrd(\Mcal)$.
\end{notation}

\begin{proposition}
Let $\Mcal$ be a presentable symmetric monoidal category. Then $\Algbrd(\Mcal)_{\Spc}$ is both a symmetric monoidal localization and a symmetric monoidal colocalization of $\Algbrd(\Mcal)$.
\end{proposition}
\begin{proof}
 Consider the pullback square
\[
\begin{tikzcd}
 \Algbrd(\Mcal)_{\Spc} \arrow{r}{i'} \arrow{d}{p'} & \Algbrd(\Mcal) \arrow{d}{p} \\
 \Spc \arrow{r}{i} & \Cat.
\end{tikzcd}
\]
Recall from proposition \ref{prop pres} that $p$ is both a cartesian and cocartesian fibration. Since $i$ has both left and right adjoints, we conclude that $i'$ has both left and right adjoints as well. 

Concretely, given an $\Mcal$-algebroid $\Acal$ with category of objects $X$, the unit $\Acal \rightarrow i'i'^L \Acal$ is a $p$-cocartesian lift of the map $X \rightarrow \prescript{\leq 0}{}{X}$ from $X$ into its geometric realization, and the counit $i'i'^R \Acal \rightarrow \Acal$ is a $p$-cartesian lift of the map $X^{\leq 0} \rightarrow X$ which includes the space of objects of $X$ inside $X$.

It remains to show that the adjoints to $i'$ are compatible with the symmetric monoidal structure on $\Algbrd(\Mcal)$. Let $\Acal$ and $\Bcal$ be a pair of $\Mcal$-algebroids with category of objects $X$ and $Y$, respectively. Denote by $\eta_\Acal: \Acal \rightarrow i'i'^L \Acal$ and $\epsilon_\Acal: i'i'^R \Acal \rightarrow \Acal$ the localization and colocalization of $\Acal$.

 Applying propositions \ref{prop products of algebroids} and \ref{prop products compatible with } together with remark \ref{remark sym mon str} we see that  the map 
 \[
 \id_\Bcal \otimes \epsilon_\Acal : \Bcal \otimes i'i'^R \Acal \rightarrow \Bcal \otimes \Acal
 \] is $p$-cartesian and lies above an $i^R$-colocal map. This implies that it is $i'^R$-colocal, and therefore we have that $\Algbrd(\Mcal)_{\Spc}$ is a symmetric monoidal colocalization of $\Algbrd(\Mcal)$.

Consider now  the map 
\[
\id_\Bcal \otimes \eta_\Acal : \Bcal \otimes \Acal \rightarrow \Bcal \otimes i'i'^L \Acal.
\] 
Its image under $p$ is $i^L$-local by a combination of proposition \ref{prop products of algebroids} and remark \ref{remark sym mon str}, together with the fact that $i^L$ preserves products. To prove that it $\id_\Bcal \otimes \eta_\Acal$ is $i'$-local it now suffices to show that it is $p$-cocartesian. Using \cite{GH} lemma 3.6.15 we see that 
\[
\id_\Bcal \boxtimes \eta_\Acal : \Bcal \boxtimes \Acal \rightarrow \Bcal \boxtimes i'i'^L \Acal
\]
is cocartesian for the projection $\Algbrd(\Mcal \times \Mcal) \rightarrow \Cat$. In other words, $\id_\Bcal \boxtimes \eta_\Acal$ exhibits $\Bcal \boxtimes i'i'^L \Acal$ as the free $\Assos_{Y \times ({}^{\leq 0}{X})}$-algebra on the $\Assos_{Y \times X}$-algebra $\Bcal \boxtimes \Acal$. Our claim now follows from the fact that the multiplication map $m: \Mcal \times \Mcal \rightarrow \Mcal$ preserves the operadic colimits involved in the description of this free algebra.
\end{proof}

\begin{proposition}\label{prop funct cuando uno es especial}
Let $\Mcal$ be a presentable symmetric monoidal category and let $\Acal, \Bcal$ be two $\Mcal$-algebroids.
\begin{enumerate}[\normalfont (i)]
\item If $\Bcal$ has a space of objects then $\Funct(\Acal, \Bcal)$ has a space of objects.
\item If $\Bcal$ is an enriched category then $\Funct(\Acal, \Bcal)$ is an enriched category. 
\end{enumerate}
\end{proposition}
\begin{proof}
Item (i) follows directly from the fact that $\Algbrd(\Mcal)_{\Spc}$ is a symmetric monoidal localization of $\Algbrd(\Mcal)$. Similarly, item (ii) follows from proposition \ref{prop Cat mcal compat}.
\end{proof}

\begin{corollary}
Let $\Mcal$ be a presentable symmetric monoidal category. The functors $\Funct(-, -)|_{\Algbrd(\Mcal)_{\Spc}^\op \times \Algbrd(\Mcal)_{\Spc}}$ and
$\Funct(-, -)|_{(\Cat^{\Mcal})^\op \times (\Cat^{\Mcal})}$  are equivalent to the internal Homs of $\Algbrd(\Mcal)_{\Spc}$ and $\Cat^\Mcal$, respectively.
\end{corollary}
\begin{proof}
This is a direct consequence of proposition \ref{prop funct cuando uno es especial}.
\end{proof}

\begin{remark}
The involution $(-)^\op: \Algbrd \rightarrow \Algbrd$ is product preserving. It follows that if $\Mcal$ is a symmetric monoidal category, then the involution $(-)^\op$ on $\Algbrd(\Mcal)$ respects the symmetric monoidal structure. In particular, given $\Mcal$-algebroids $\Acal, \Bcal$, there is an equivalence
\[
\Funct(\Acal^\op, \Bcal^\op)^\op = \Funct(\Acal, \Bcal).
\]
\end{remark}

Our next goal is to provide a concrete description of the product of cells, and use it to study functor algebroids in the case when the source is a cell.

\begin{notation}\label{notation gluing}
Let $\Mcal$ be a monoidal category with an initial object compatible with the monoidal structure and let $m, m'$ be objects in $\Mcal$. Let $X = \lbrace a, b , c \rbrace$ be the set with three elements. Let $C_{m, m'}$ be the free $\Assos_X$-algebra equipped with maps $m \rightarrow C_{m, m'}(b, c)$ and $m' \rightarrow C_{m, m'}(a, b)$. This is characterized by the following properties:
\begin{itemize}
\item $C_{m, m'}(a, a) = C_{m, m'}(b, b) = C_{m, m'}(c, c) = 1_\Mcal$.
\item $C_{m, m'}(a, b) = m'$
\item $C_{m, m'}(b, c) = m$.
\item $C_{m, m'}(a, c) = m \otimes m'$
\item $C_{m, m'}(b, a) = C_{m, m'}(c, a) = C_{m, m'}(c, b)$ is the initial object of $\Mcal$.
\end{itemize}
We note that $C_{m,m'}$ fits into a pushout
\[
\begin{tikzcd}
1_{\Mcal} \arrow{r}{} \arrow{d}{} & C_m \arrow{d}{} \\ C_{m'} \arrow{r}{} & C_{m, m'}
\end{tikzcd}
\]
where the top horizontal arrow and left vertical arrows pick out the target and source objects, respectively.
\end{notation}

\begin{remark}\label{remark product cells}
Let $\Mcal$ be a  symmetric monoidal category with compatible initial object, and let $m, m'$ be objects in $\Mcal$. Then the algebroid $C_m \otimes C_{m'}$ has objects $(i, j)$ for $0 \leq i, j \leq 1$. Its morphisms can be depicted schematically as follows:
\[
\begin{tikzcd}[column sep = large, row sep = large]
(0, 0) \arrow{r}{m} \arrow{d}{m'} \arrow{dr}{m \otimes m'} & (1, 0) \arrow{d}{m'} \\
(0, 1) \arrow{r}{m} & (1, 1)
\end{tikzcd}
\]
Every Hom-object which is not associated to an arrow in the above diagram is the initial object in $\Mcal$. Note that $C_m \otimes C_{m'}$ fits into a commutative square 
\[
\begin{tikzcd}
C_{m \otimes m'} \arrow{r}{} \arrow{d}{} & C_{m, m'} \arrow{d}{} \\
C_{m', m} \arrow{r}{} & C_{m} \otimes C_{m'}
\end{tikzcd}
\]
where:
\begin{itemize}
\item The right vertical arrow picks out the $m$-cell between $(0,0)$ and $(1, 0)$ and the $m'$-cell between $(1,0)$ and $(1,1)$.
\item The bottom horizontal arrow picks out the $m'$-cell between $(0,0)$ and $(0, 1)$ and the $m$-cell between $(0,1)$ and $(1, 1)$. 
\item The cell $C_{m, \otimes m'} \rightarrow C_{m} \otimes C_{m'}$ the $m\otimes m'$-cell between $(0,0)$ and $(1,1)$.
\end{itemize}
\end{remark}

\begin{proposition}\label{prop product cells}
Let $\Mcal$ be a  symmetric monoidal category with compatible initial object, and let $m, m'$ be objects in $\Mcal$. Then the commutative square of $\Mcal$-algebroids 
\[
\begin{tikzcd}
C_{m \otimes m'} \arrow{r}{} \arrow{d}{} & C_{m, m'} \arrow{d}{} \\
C_{m', m} \arrow{r}{} & C_{m} \otimes C_{m'}
\end{tikzcd}
\]
from remark \ref{remark product cells}, is a pushout square.
\end{proposition}
\begin{proof}
Observe first that the induced square at the level of objects is given by
\[
\begin{tikzcd}
\lbrace (0, 0) , (1, 1) \rbrace  \arrow{r}{k} \arrow{d}{l} & \lbrace (0, 0), (1, 0) , (1, 1) \rbrace \arrow{d}{j} \\
\lbrace (0, 0), (0, 1) , (1, 1) \rbrace \arrow{r}{i} & \lbrace 0 , 1 \rbrace \times \lbrace 0, 1 \rbrace
\end{tikzcd}
\]
and is indeed a pushout square. Consider the following $\Mcal$-algebroids:
\begin{align*}
\overline{C_{m \otimes m'}} &= (jk)_!C_{m \otimes m'}\\
\overline{C_{m, m'}} &= j_! C_{m, m'} \\
\overline{C_{m', m}} &= i_! C_{m', m}
\end{align*}
Here we denote with $(-)_!$ the functors induced from the fact that the projection $\Algbrd(\Mcal) \rightarrow \Cat$ is a cocartesian fibration. In other words, the above algebroids are obtained from the previous ones by adding extra objects so that they all have set of objects $\lbrace 0, 1 \rbrace^2$, where the new Hom objects are declared to be the initial object of $\Mcal$.  Combining proposition \ref{prop pres} with \cite{HTT} proposition 4.3.1.9, we reduce to showing that the induced square
\[
\begin{tikzcd}
\overline{C_{m \otimes m'}} \arrow{r}{} \arrow{d}{} & \overline{C_{m, m'}} \arrow{d}{} \\
\overline{C_{m', m}} \arrow{r}{} & C_{m} \otimes C_{m'}
\end{tikzcd}
\]
is a pushout square in $\Algbrd_{\lbrace 0, 1 \rbrace ^2}(\Mcal)$. Our claim now follows from the fact that the images of the above square under the evaluation functors $\Algbrd_{\lbrace 0, 1 \rbrace ^2}(\Mcal) \rightarrow \Mcal$ are pushout squares.
\end{proof}

\begin{corollary}\label{coro homs en funct}
Let $\Mcal$ be a presentable symmetric monoidal category and let $m$ be an object in $\Mcal$. Let $\Acal$ be an $\Mcal$-algebroid and let $\mu, \nu: C_m \rightarrow \Acal$ be two $m$-cells in $\Acal$. Denote by $\shom_\Mcal$ the internal Hom functor for $\Mcal$. Then there is a pullback square
\[
\begin{tikzcd}
 \Funct(C_m, \Acal)(\nu, \mu) \arrow{r}{} \arrow{d}{} & \Acal(\nu(0), \mu(0))) \arrow{d}{} \\
\Acal(\nu(1), \mu(1)) \arrow{r}{} & \shom_\Mcal(m, \Acal(\nu(1), \mu(0)))
\end{tikzcd}
\]
where  the top horizontal and left vertical arrows are given by the source and target maps, and the right horizontal and bottom vertical arrows are induced by composition with the cells $\nu$ and $\mu$, respectively.
\end{corollary}
\begin{proof}
Let $m'$ be another object of $\Mcal$.  The pushout square of remark \ref{remark product cells} can be enhanced to a colimit diagram
\[
\begin{tikzcd}
1_{\Mcal} \sqcup 1_{\Mcal} \arrow[d] \arrow[r] \arrow[rd] & C_m \sqcup 1_{\Mcal} \arrow[rd]      &                       \\
1_{\Mcal} \sqcup C_{m}  \arrow[rd]                        & C_{m \otimes m'} \arrow[d] \arrow[r] & {C_{m, m'}} \arrow[d] \\
                                                          & {C_{m', m}} \arrow[r]                & C_{m} \otimes C_{m'} 
\end{tikzcd}
\]
where map $1_\Mcal \sqcup 1_\Mcal \rightarrow C_{m \otimes m'}$ picks out the source and target objects, and the maps $1_{\Mcal } \sqcup C_m \rightarrow C_{m', m}$ and $C_m \sqcup 1_\Mcal \rightarrow C_{m, m'}$ pick out the $m$-cell and the object which is not contained in it.

Let $X$ be the category of objects of $\Acal$. We have an induced limit diagram of spaces
\[
\begin{tikzcd}
X^{\leq 0} \times X^{\leq 0}                         & {\Hom(C_m, \Acal) \times X^{\leq 0}} \arrow[l]    &                                                           \\ { X^{\leq 0} \times \Hom(C_m, \Acal)} \arrow[u] & {\Hom(C_{m \otimes m'},\Acal)} \arrow[lu]    & {\Hom(C_{m, m'}, \Acal)} \arrow[l] \arrow[lu]           \\
                                            & {\Hom(C_{m', m},\Acal)} \arrow[lu] \arrow[u] & {\Hom(C_{m} \otimes C_{m'}, \Acal)} \arrow[l] \arrow[u]
\end{tikzcd}
\]
where all Homs are taken in $\Algbrd(\Mcal)$.

Consider now the commutative diagram of spaces
\[
\begin{tikzcd}
{[0] }\arrow[d, "{(\mu(0), \nu)}"] \arrow[r]               & {[0] }\arrow[d, "{(\mu(0), \nu(1))}"] & {[0]} \arrow[d , "{(\mu, \nu(1))}" ] \arrow[l]               \\
X^{\leq 0} \times \Hom_{\Algbrd(\Mcal)}(C_m, \Acal) \arrow[r] & X^{\leq 0} \times X^{\leq 0}  & {\Hom_{\Algbrd(\Mcal)}(C_m, \Acal) \times X^{\leq 0}} \arrow[l]
\end{tikzcd}
\]
where the bottom left horizontal arrow is the target map, and the bottom right horizontal arrow is the source map.
Pulling back our previous diagram along this yields a cartesian square of spaces
\[
\begin{tikzcd}
\Hom_\Mcal(m',\shom_\Mcal(m, \Acal(\nu(1), \mu(0))))  & \arrow{l}{} \Hom_{\Mcal}(m', \Acal(\nu(1), \mu(1))) \\
\Hom_{\Mcal}(m', \Acal(\nu(0), \mu(0))) \arrow{u}{} &  \arrow{l}{}  \arrow{u}{} \Hom_\Mcal(m', \Funct(C_m, \Acal)(\nu, \mu)).
\end{tikzcd}
\]
Our claim now follows from the fact that the above square is natural in $m'$.
\end{proof}

We finish by studying the behavior of functor algebroids under changes in the enriching categories.

\begin{proposition}\label{prop funct y changes}
Let $G: \Mcal \rightarrow \Mcal'$ be a colimit preserving symmetric monoidal functor between presentable symmetric monoidal categories. Let $\Ccal$ be an $\Mcal$-algebroid and let $\Dcal$ be an $\Mcal'$ algebroid. There is an equivalence of $\Mcal$-algebroids
\[
\Funct(\Ccal, (G^R)_!\Dcal) = (G^R)_!\Funct(G_! \Ccal , \Dcal)
\]
which is natural in $\Ccal$ and $\Dcal$.
\end{proposition}
\begin{proof}
Let $\Ecal$ be another $\Mcal$-algebroid. Then we have a chain of equivalences
\begin{align*}
\Hom_{\Algbrd(\Mcal)}(\Ecal , \Funct(\Ccal , (G^R)_! \Dcal)) &= \Hom_{\Algbrd(\Mcal)}(\Ecal \otimes \Ccal, (G^R)_!\Dcal) \\
 & = \Hom_{\Algbrd(\Mcal')}( G_!(\Ecal \otimes \Ccal), \Dcal) \\
 & = \Hom_{\Algbrd(\Mcal')}( G_!\Ecal \otimes G_!\Ccal, \Dcal) \\
 & = \Hom_{\Algbrd(\Mcal')}( G_! \Ecal , \Funct(G_!\Ccal, \Dcal))\\
 & = \Hom_{\Algbrd(\Mcal)}( \Ecal , (G^R)_! \Funct(G_! \Ccal, \Dcal) ).
\end{align*}
Our claim follows from the fact that the above equivalences are natural in $\Ecal, \Ccal$ and $\Dcal$.
\end{proof}

\begin{proposition}\label{coro funct y extensions}
Let $i: \Mcal \rightarrow \Mcal'$ be a lax symmetric monoidal functor between presentable symmetric monoidal categories. Assume that $i$  is fully faithful and admits a left adjoint which is strictly symmetric monoidal. Let $\Ccal$ be an $\Mcal'$-algebroid, and $\Dcal$ be an $\Mcal$ algebroid. Then $\Funct(\Ccal, i_!\Dcal)$ belongs to the image of $i_!: \Algbrd(\Mcal) \rightarrow \Algbrd(\Mcal')$.
\end{proposition}
\begin{proof}
Our conditions guarantee that $i_! : \Algbrd(\Mcal) \rightarrow \Algbrd(\Mcal')$ is fully faithful, and admits a symmetric monoidal left adjoint. To show that $\Funct(\Ccal, i_!\Dcal)$ belongs to the image of $i_!$, it suffices to show that it is local for the maps $q_\Ecal : \Ecal \rightarrow i_! i^L_! \Ecal$, for each $\Mcal'$-algebroid $\Ecal$. 
Indeed, the map
\[
\Hom_{\Algbrd(\Mcal')}(i_! i^L_! \Ecal, \Funct(\Ccal, i_! \Dcal)) \rightarrow \Hom_{\Algbrd(\Mcal')}(\Ecal, \Funct(\Ccal, i_! \Dcal))
\]
of precomposition with $q_{\Ecal}$ is equivalent to the map
\[
\Hom_{\Algbrd(\Mcal')}(i_!i^L_! \Ecal \otimes \Ccal, i_! \Dcal) \rightarrow \Hom_{\Algbrd(\Mcal')}(\Ecal \otimes \Ccal, i_! \Dcal)
\]
of precomposition with $q_{\Ecal} \otimes \id_{\Ccal}$. It therefore suffices to show that the induced map
\[
i^L_! (\Ecal \otimes \Ccal) \rightarrow i^L_!( i_! i_!^L \Ecal \otimes \Ccal) 
\] 
is an equivalence. This follows from the fact that $i^L_!$ is a symmetric monoidal localization.
\end{proof}

\begin{corollary}\label{coro compara internal functs}
Let $i: \Mcal \rightarrow \Mcal'$ be a  symmetric monoidal functor between presentable symmetric monoidal categories. Assume that $i$  is fully faithful and has a left adjoint which is strictly symmetric monoidal. Let $\Ccal, \Dcal$ be $\Mcal$-algebroids. Then there is an equivalence $i_!\Funct(\Ccal, \Dcal) = \Funct(i_! \Ccal, i_! \Dcal)$, which is natural in $\Ccal, \Dcal$.
\end{corollary}
\begin{proof}
This is a direct consequence of proposition \ref{coro funct y extensions}, as $i_!\Funct(\Ccal, \Dcal)$ and $\Funct(i_! \Ccal, i_! \Dcal)$ both corepresent the same functor on $\Algbrd(\Mcal)$.
\end{proof}

\subsection{$\omega$-categories}\label{subsection ncat}
The theory of $n$-categories can be obtained as a special case of the general notion of enriched categories. We refer the reader to \cite{GH} and \cite{Hinich} for proofs that the following definition agrees with other models for the theory of $n$-categories.

\begin{definition}
Let $1\kr\Cat$ be the category of categories. We inductively define for each $n \geq 2$ the cartesian closed presentable category $\nCat$ of  $n$-categories to be the category $\Cat^{(n-1)\kr\Cat}$ of categories enriched in the cartesian closed presentable category $(n-1)\kr \Cat$.
\end{definition}

\begin{construction}\label{construction inclusions ncat}
Let $i: \Spc \rightarrow \Cat$ be the inclusion. This admits a left adjoint, and in particular it  has a canonical symmetric monoidal structure, where we equip $\Spc$ and $\Cat$ with their cartesian symmetric monoidal structures. Specializing the discussion of remark \ref{remark shriek fully faithful} we obtain a commutative square of presentable categories and colimit preserving morphisms
\[
\begin{tikzcd}
\Algbrd(\Spc)_{\Spc} \arrow{d}{i_!} \arrow{r}{} & \Cat \arrow{d}{i_!} \\
\Algbrd(\Cat)_{\Spc} \arrow{r}{} & \twoCat
\end{tikzcd}
\]
which is both horizontally right adjointable and vertically left adjointable, and whose vertical arrows are fully faithful.

Denote by $i^{1, 2}$ the functor $i_!: \Cat \rightarrow \twoCat$. Arguing by induction, we obtain for all $n \geq 2$ a commutative square of presentable categories and colimit preserving functors 
\[
\begin{tikzcd}
\Algbrd((n-1)\kr\Cat)_{\Spc} \arrow{d}{(i^{n-1, n})_!} \arrow{r}{} & \nCat \arrow{d}{i^{n, n+1}} \\
\Algbrd(n\kr\Cat)_{\Spc} \arrow{r}{} & (n+1)\kr\Cat
\end{tikzcd}
\]
which is both horizontally right adjointable and vertically left adjointable, and whose vertical arrows are fully faithful. In particular, we have a sequence of presentable categories and left adjointable colimit preserving fully faithful functors
\[
0\kr\Cat \xrightarrow{i^{0, 1}} \Cat \xrightarrow{i^{1, 2}} \twoCat \xrightarrow{i^{2, 3}} 3\kr\Cat \ldots 
\]
where we set $0\kr\Cat = \Spc$ and $i^{0, 1} = i$. For each pair $m \geq n$ denote by $i^{n, m}: \nCat \rightarrow m\kr\Cat$ the corresponding inclusion.
\end{construction}

\begin{example}
Let $C_0$ be the terminal object in $\Spc$, and $C_1$ be the arrow category $[1]$. We inductively define for each $n \geq 2$ an $n$-category $C_n$ as the cell associated to the $(n-1)$-category $C_{n-1}$. We call $C_n$ the $n$-cell. It follows by induction that the category $\nCat$ is compactly generated by the object $C_n$. Note that for each $0 \leq k < n$ there are source and target maps $i^{k, n} C_k \rightarrow C_n$.
\end{example}

In practice, we usually identify $n$-categories with their images under the functors $i^{n, m}$. In other words, many times we implicitly work in the direct limit of the categories $\nCat$.

\begin{notation}
Let $\omega\kr\Cat$ be the direct limit in $\Pr^L$ of the sequence of construction \ref{construction inclusions ncat}. We call it the category of $\omega$-categories. For each $n \geq 0$ denote by $i^{n, \omega}: \nCat \rightarrow \omega\kr\Cat$ the induced map.
\end{notation}

\begin{remark} \label{remark descript omegacat}
We can alternatively think about $\omega\kr\Cat$ as the limit of the categories $\nCat$ under the functors $(i^{n, n+1})^R : (n+1)\kr\Cat \rightarrow \nCat$.  In other words, an $\omega$-category is a compatible family of $n$-categories for all $n \geq 0$. The resulting projections $\omega\kr\Cat \rightarrow \nCat$ are right adjoint to the maps $i^{n, \omega}$.

Note that $i^{n,n+1}$ preserves compact objects for every $n \geq 0$. Therefore $\omegaCat$ is compactly generated and the maps $i^{n, \omega}$ preserve compact objects. In particular, we have that the projections $(i^{n, \omega})^R$ preserve filtered colimits.
\end{remark}

\begin{remark}\label{remark omegacat as sequences}
The sequence of categories from construction \ref{construction inclusions ncat} yields a functor 
\[
(-)\kr\Cat : \NN \rightarrow \Cat
\]
 where $\NN$ is the poset of natural numbers. Let $p: \Ecal \rightarrow \NN$ be the associated cocartesian fibration. Since the functors $i^{n, n+1}$ admit right adjoints, this is also a cartesian fibration. By virtue of remark \ref{remark descript omegacat}, the category $\omega\kr\Cat$ is the category of cartesian sections of $p$. Note that  $\Ecal$ is a full subcategory of the product $\NN \times (\bigcup_{n \geq 0} \nCat)$. It follows that  $\omega\kr\Cat$ is the full subcategory of the functor category $\Funct(\NN, \bigcup_{n \geq 0} \nCat)$ on those sequences of objects
\[
\Ccal_0 \xrightarrow{i_0} \Ccal_1 \xrightarrow{i_1} \ldots
\]
such that $\Ccal_n$ is an $n$-category for each $n \geq 0$, and the map $(i^{n, n+1})^R i_n$ is an isomorphism for each $n \geq 0$.
\end{remark}

\begin{proposition}\label{proposition inclusions into omegacat}
For each $n \geq 0$ the inclusion $i^{n , \omega} : \nCat \rightarrow \omega\kr\Cat$ is fully faithful and admits both a left and a right adjoint.
\end{proposition}
\begin{proof}
The existence of a right adjoint was already observed in remark \ref{remark descript omegacat}.  By the same argument as in remark \ref{remark omegacat as sequences}, we have an equivalence between $\omega\kr\Cat$ and the full subcategory of the functor category $\Funct(\NN_{\geq n}, \bigcup_{m \geq n} m\kr\Cat)$ on those sequences
\[
\Ccal_n \xrightarrow{i_n} \Ccal_{n+1} \xrightarrow{i_{n+1}} \ldots
\]
such that $\Ccal_m$ is an $m$-category for each $m \geq n$, and the map $(i^{m, m+1})^R i_m$ is an isomorphism for each $m \geq n$.  In this language, the projection $(i^{n, \omega})^R$ is given by the corestriction to $\nCat$ of the composite map
\[
\omega\kr\Cat \hookrightarrow \Funct(\NN_{\geq n}, \bigcup_{m \geq n} m\kr\Cat) \xrightarrow{\ev_n} \bigcup_{m \geq n} m\kr\Cat.
\] 

The left adjoint to  $\ev_n$  is given by left Kan extension along the inclusion $\lbrace n \rbrace \rightarrow \NN_{\geq n}$, and is fully faithful. The fact that $i^{n, \omega}$ is fully faithful follows from the fact that the left adjoint to $\ev_n$ maps an $n$-category $\Ccal$ to the constant diagram 
\[
\Ccal \rightarrow \Ccal \rightarrow \ldots
\]
and this belongs to $\omega\kr\Cat$. 

It remains to show that $i^{n, \omega}$ admits a left adjoint. We have that $i^{n, \omega}$ is given by the corestriction to $\omega\kr\Cat$ of the composite map
\[
\nCat \hookrightarrow \bigcup_{m \geq n} m\kr\Cat \rightarrow \Funct(\NN_{\geq n}, \bigcup_{m \geq n} m\kr\Cat)
\]
where the second arrow is the functor of precomposition with the projection $\NN_{\geq n} \rightarrow \lbrace n \rbrace$. The composition of the two arrows above is limit and colimit preserving. Furthermore, the inclusion $\omega\kr\Cat \rightarrow \Funct(\NN_{\geq n}, \bigcup_{m \geq n} m\kr\Cat)$ is accessible and limit preserving. We conclude that $i^{n, \omega}$ is accessible and limit preserving, and the adjoint functor theorem guarantees that it admits a left adjoint, as desired.
\end{proof}

\begin{notation}
Let $n \geq 0$. We denote by $(-)^{\leq n}$ and $\prescript{{\leq n}}{}{(-)}$ the right and left adjoints to the inclusion $i^{n,\omega}$. We think about these as the functors that discard (resp. invert) cells of dimension greater than $n$.
\end{notation}

\begin{remark}\label{remark omega cats como union}
We have a diagram
\[
i^{0 , \omega} (-)^{\leq 0} \rightarrow i^{1, \omega} (-)^{\leq 1} \rightarrow i^{2 , \omega} (-)^{\leq 2} \rightarrow \ldots \rightarrow \id_{\omega\kr\Cat}
\]
of endofunctors of $\omega\kr\Cat$, where the transitions are induced by the counits of the adjunctions $i^{n, n+1} \dashv (i^{n, n+1})^R$. We note that for every $n \geq 0$ the composition of the above sequence with the functor $(-)^{\leq n}$ is eventually constant, and is therefore a filtered colimit diagram.

Recall from remark \ref{remark descript omegacat} that $\omega\kr\Cat$ is the limit of the sequence of categories
\[
0\kr\Cat \xleftarrow{(i^{0, 1})^R} 1\kr\Cat \xleftarrow{(i^{1, 2})^R} 2\kr\Cat \xleftarrow{(i^{2, 3})^R} \ldots \leftarrow \omega\kr\Cat.
\]
Each of the transition functors above preserve filtered colimits. We conclude that $\id_{\omega\kr\Cat}$ is the colimit of the endofunctors $i^{n, \omega} (-)^{\leq n}$. In other words, every $\omega$-category is the colimit of its truncations.
\end{remark}

\begin{proposition}
The category $\omega\kr\Cat$ is cartesian closed.
\end{proposition}
\begin{proof}
Let $\Dcal$ be an $\omega$-category and let $X^{\rhd}: \Ical^\rhd \rightarrow \omega\kr\Cat$ be a colimit diagram for $X = X^\rhd|_{\Ical}$. We have to show that the induced diagram $X^{\rhd} \times \Dcal$ is also a colimit diagram. For each $n \geq 0$ let $X_n^\rhd$ be an extension of $i^{n,\omega} (X^{\leq n})$ to a colimit diagram in $\omega\kr\Cat$. Thanks to remark \ref{remark omega cats como union}, we have $X^\rhd = \colim_n X^\rhd_n$. Since $\omegaCat$ is compactly generated we have that products distribute over filtered colimits, and therefore we have
\[
X^\rhd \times \Dcal = \colim_n( X_n^\rhd \times \Dcal).
\]
It suffices to show that $X^\rhd_n \times \Dcal$ is a colimit diagram for each $n \geq 0$. Using remark \ref{remark omega cats como union} again, we have
\[
X^\rhd_n \times \Dcal = \colim_m (X^\rhd_n \times i^{m,\omega}\Dcal^{\leq m}).
\]
It therefore suffices to show that for each $n, m \geq 0$, the diagram $X^\rhd_n \times i^{m, \omega} \Dcal^{\leq m}$ is a colimit diagram in $\omega\kr\Cat$. This follows from the fact that $(\max(n, m))\kr\Cat$ is cartesian closed and closed under products and colimits inside $\omegaCat$.
\end{proof}

\begin{notation}
We denote by 
\[
\Funct(-,-): \omegaCat^\op \times \omegaCat \rightarrow \omegaCat
\]
the internal Hom of $\omegaCat$.
\end{notation}

\begin{remark}\label{remark enrichement}
Equip the category $\omega\kr\Cat$ with its cartesian symmetric monoidal structure. For each category $X$ we have that $\Algbrd_X(\omega\kr\Cat)$ is equivalent to the limit of the categories $\Algbrd_X(\nCat)$ under the transition functors $(i^{n, n+1})^R_!$. Integrating over all such $X$ we see that $\Algbrd(\omega\kr\Cat)$ is the limit of the categories $\Algbrd(\nCat)$ under the functors $(i^{n, n+1})^R_!$. 

Using proposition  \ref{proposition inclusions into omegacat} we see that the inclusion $i^{0,\omega}$ is the unit map for the presentable monoidal category $\omegaCat$. It follows that that an object in $\Algbrd(\omega\kr\Cat)$ is an enriched category if and only if its image in $\Algbrd(n\kr\Cat)$ is an enriched category for all $n \geq 0$. Therefore the category $\Cat^{\omega\kr\Cat}$ is the limit of the categories $\Cat^{n\kr\Cat}$ under the transition functors $(i^{n, n+1})^R_!$. Passing to left adjoints we conclude that the colimit in $\Pr^L$ of the diagram
\[
\Cat^{\Spc} \xrightarrow{i^{0, 1}_!} \Cat^{\Cat} \xrightarrow{i^{1, 2}_!} \Cat^{\twoCat} \xrightarrow{i^{2, 3}_!} \ldots
\]
is $\Cat^{\omega\kr\Cat}$. The above is however equivalent to the diagram in construction \ref{construction inclusions ncat}. It follows that there is an equivalence
\[
\Cat^{\omega\kr\Cat} = \omega\kr\Cat
\]
which makes the following diagram commute for all $n \geq 0$:
\[
\begin{tikzcd}
\Cat^{\nCat} \arrow{d}{ = } \arrow{r}{i^{n , \omega}_!} & \Cat^{\omega\kr\Cat} \arrow{d}{} \\
(n+1)\kr\Cat \arrow{r}{ i^{n+1, \omega} } & \omega\kr\Cat
\end{tikzcd}
\]
In other words, an $\omega$-category can be thought of as a category enriched in $\omega$-categories, in a way which is compatible with the definition of $(n+1)$-categories as categories enriched in $n$-categories. In particular, the internal Hom functor for $\omegaCat$ fits into the framework of functor enriched categories from \ref{subsection multiplicativity}.
\end{remark}

\begin{proposition}\label{prop left adj preserva products}
The functor $\prescript{\leq n}{}{(-)}: \omegaCat \rightarrow \nCat$ preserves finite products for all $n \geq 0$.
\end{proposition}
\begin{proof}
We first observe that the final object of $\omegaCat$ is the image of the final object of $\nCat$ under the (limit preserving) inclusion $i^{n, \omega}$. Hence $\prescript{\leq n}{}{(-)}$ preserves final objects. 

We now show that $\prescript{\leq n}{}{(-)}$ preserves binary products. Let $\Ccal, \Dcal$ be two $\omega$-categories. Thanks to remark \ref{remark omega cats como union}, we can write $\Ccal = \colim_{m\geq 0} \Ccal_m$ and $\Dcal = \colim_{m \geq 0} \Dcal_m$ where $\Ccal_m$ and $\Dcal_m$ are $m$-categories for each $m \geq 0$. Since products commute with filtered colimits in $\omegaCat$, we have 
\[
\Ccal \times \Dcal = \colim_{m \geq 0} \Ccal_m \times \Dcal_m.
\]
Since $\prescript{\leq n}{}{(-)}$ preserves colimits, we reduce to showing that for each $m \geq n$ the left adjoint to the inclusion $i^{n, m}$ preserves products. Arguing inductively, we may furthermore reduce to the case $n = 0$, $m = 1$, which follows from the fact that the geometric realization functor $\Cat \rightarrow \Spc$ preserves finite products.
\end{proof}

\begin{corollary}
Let $\Ccal$ and $\Dcal$ be $\omega$-categories, and assume that $\Dcal$ is an $n$-category for some $n \geq 0$. Then $\Funct(\Ccal, \Dcal)$ is an $n$-category. Moreover, if $\Ccal$ is also an $n$-category then $\Funct(\Ccal, \Dcal)$ can be identified with the internal Hom between $\Ccal$ and $\Dcal$ in $\nCat$.
\end{corollary}
\begin{proof}
Combine proposition \ref{prop left adj preserva products}, proposition \ref{coro funct y extensions} and corollary \ref{coro compara internal functs}.
\end{proof}

\begin{remark}\label{remark compatibility enrichement y localization}
Let $n \geq 0$. Passing to right adjoints in the commutative diagram of remark \ref{remark enrichement} yields a commutative square
\[
\begin{tikzcd}
\Cat^{\nCat} \arrow{d}{ = } &  \arrow{l}[swap]{(-)^{\leq n}_!} \Cat^{\omega\kr\Cat} \arrow{d}{} \\
(n+1)\kr\Cat & \arrow{l}[swap]{(-)^{\leq n+1}} \omega\kr\Cat.
\end{tikzcd}
\]

It follows that for every $\omega$-category $\Ccal$, the category underlying $\Ccal$ (thought of as an object of $\Cat^{\omega\kr\Cat}$) is $\Ccal^{\leq 1}$. In particular, its space of objects is $\Ccal^{\leq 0}$. Furthermore, for each par of objects $x, y$ in $\Ccal$, we have an equivalence
\[
\Hom_{\Ccal^{\leq n+1}}(x, y) = \Hom_{\Ccal}(x, y)^{\leq n}.
\]

Similarly, passing  to left adjoints in the commutative diagram of remark \ref{remark enrichement} yields a commutative square
\[
\begin{tikzcd}
\Cat^{\nCat} \arrow{d}{ = } &  \arrow{l}[swap]{\prescript{\leq n}{}{(-)}_!} \Cat^{\omega\kr\Cat} \arrow{d}{} \\
(n+1)\kr\Cat & \arrow{l}[swap]{\prescript{\leq n+1}{}{(-)}} \omega\kr\Cat.
\end{tikzcd}
\]

It follows that for every $\omega$-category $\Ccal$, the $(n+1)$-category $\prescript{\leq n+1}{}{\Ccal}$  is the image under the localization functor $\Algbrd(\nCat)_{\Spc} \rightarrow (n+1)\kr\Cat$ of an algebroid $\prescript{\leq n+1}{}{\Ccal}^{\text{pre}}$ with space of objects $\Ccal^{\leq 0}$, and such that for every pair of objects $x, y$ in $\Ccal^{\leq 0}$ we have an equivalence
\[
\Hom_{\prescript{\leq n+1}{}{\Ccal}^{\text{pre}}}(x, y) = \prescript{\leq n}{}{\Hom_{\Ccal}(x, y)}.
\]
\end{remark}

\begin{example}
Let $m >  n \geq 1$. It follows by induction that $\prescript{\leq n}{}{C_m}$ is the singleton set. On the other hand, $C_m^{\leq n}$ is the boundary of the $(n+1)$-cell $\partial C_{n+1}$, defined inductively by the fact that $\partial C_0$ is empty, and for $n \geq 1$ we have $\partial C_n = C_{\partial C_{n-1}}$.
\end{example}

\begin{remark}\label{remark detectar n categories}
Looking at the unit of the adjunction $i^{(n+1), \omega} \dashv (-)^{\leq n+1}$ through the equivalence given by the first commutative square in remark \ref{remark compatibility enrichement y localization} shows that an $\omega$-category $\Ccal$ is an $(n+1)$-category for some $n \geq 0$ if and only for every pair of objects $x, y$ in $\Ccal$, the $\omega$-category $\Hom_{\Ccal}(x, y)$ is an $n$-category. 
\end{remark}

\begin{remark}
Recall   from remarks \ref{remark op algebroids} and \ref{remark op enriched} that we have an involution $(-)^\op$ on the full subcategory of $\Algbrd$ on the enriched categories, which restricts to an involution on the category of algebroids over any symmetric operad. In particular, for each $n \geq 1$ we have an induced involution on $\nCat$ by virtue of its description  as $\Cat^{(n-1)\kr\Cat}$.

 It follows by induction that $\nCat$ comes equipped with $n$ commuting involutions $(-)^{k\text{-op}}$ for $1 \leq k \leq n$. We think about $(-)^{k\text{-op}}$ as the involution that inverts the direction of all $k$-cells. These involutions are compatible with the inclusions $i^{n, n+1}$, and they therefore induce an infinite family of commuting involutions on the category $\omega\kr\Cat$.
\end{remark}



\tableofcontents

\section{Modules over algebroids}

Let $\Mcal$ be a monoidal category and let $\Ccal$ be a category left tensored over $\Mcal$. A left module in $\Ccal$ for an $\Mcal$-algebroid $\Acal$ consists of:
\begin{itemize}
\item For each object $x$ in $\Acal$ an object $\Pcal(x)$ in $\Ccal$.
\item For every pair of objects $x, y$ in $\Acal$ a morphism $\Acal(y, x) \otimes \Pcal(x) \rightarrow \Pcal(y)$.
\item An infinite list of compatibility data between the above morphisms and the structure maps for $\Acal$.
\end{itemize}

More generally, assume given another monoidal category $\Mcal'$, and an $\Mcal-\Mcal'$-bimodule category $\Ccal$. If $\Acal$ and $\Bcal$ are algebroids in $\Mcal$ and $\Mcal'$ respectively, an $\Acal-\Bcal$-bimodule in $\Ccal$ consists of:
\begin{itemize}
\item For each pair of objects $x$ in $\Acal$ and $y$ in $\Bcal$, an object $\Pcal(x, y)$ in $\Ccal$.
\item For every pair of objects $x, x'$ in $\Acal$ and  object $y$ in $\Bcal$, a morphism
\[
\Acal(x', x) \otimes \Pcal(x, y) \rightarrow \Pcal(x', y).
\]
\item For every pair of objects $y, y'$ in $\Bcal$ and object $x$ in $\Acal$, a morphism
\[
\Pcal(x, y) \otimes \Bcal(y, y') \rightarrow \Pcal(x, y')
\]
\item An infinite list of compatibility data between the above morphisms and the structure maps for $\Acal$ and $\Bcal$.
\end{itemize}

Our goal in this section is to review the theory of left modules and bimodules, and provide a functorial enhancement of the  procedure of enrichment of presentable modules over presentable monoidal categories.

We begin in \ref{subsection leftmod} by using the operads $\LM_X$ from section \ref{subsection assosx} to define the category of left modules over an algebroid. We show that there is a well behaved procedure of restriction of scalars along morphisms of algebroids. We record here two basic results regarding the multiplicativity properties of the theory of left modules, analogous the ones obtained in section \ref{subsection multiplicativity} for the theory of algebroids.

In \ref{subsection enrichement} we construct, for each presentable symmetric monoidal category $\Mcal$, a lax symmetric monoidal functor 
\[
\theta_\Mcal : \Mcal\modd(\Pr^L) \rightarrow \widehat{\Cat}^\Mcal.
\]
For each presentable $\Mcal$-module $\Ccal$, the enriched category $\theta_\Mcal(\Ccal)$ has $\Ccal$ as its underlying category, and for each pair of objects $x, y$ in $\Ccal$ one has an isomorphism between $\Hom_{\theta_\Mcal(\Ccal)}(x, y)$ and the Hom object  $\shom_\Ccal(x, y)$ obtained from the action of $\Mcal$ on $\Ccal$. We show that the functor $\theta_\Mcal$ is compatible with changes in the enriching category. In the particular case when $\Mcal$ is the category of spaces with its cartesian symmetric monoidal structure, we prove that the functor $\theta_\Mcal$ is equivalent, as a lax symmetric monoidal functor, to  the usual forgetful functor from $\Mcal\modd(\Pr^L)$ to $\widehat{\Cat}$. As a first consequence of the existence and properties of $\theta_\Mcal$ we show that $\Mcal$ admits a canonical enrichment over itself. This allows us in particular to construct an $(n+1)$-category of $n$-categories for each $n \geq 0$, and in the limit it provides a definition of the $\omega$-category of $\omega$-categories.

In \ref{subsection bimod} we review the notion of bimodule over an algebroid. We recall here the approach to the construction of the Yoneda embedding via the diagonal bimodule and the folding construction from \cite{Hinich}, and record a basic result regarding the procedure of restriction of scalars in the context of bimodules over algebroids.

\subsection{Left modules}\label{subsection leftmod}

We begin by reviewing the concept of left module over an algebroid.

\begin{notation}
For each $\LM$-operad $\Mcal$ we denote by $\Mcal_l$ the its $\Assos$-component, and by $\Mcal_m$ the fiber of $\Mcal$ over the module object in $\LM$.
\end{notation}

\begin{definition}
Let $\Mcal$ be a $\LM$-operad. Let $\Acal$ be an algebroid in $\Mcal_l$ with category of objects $X$. A left $\Acal$-module is a $\LM_X$-algebra in $\Mcal$, whose $\Assos_X$-component is identified with $\Acal$.
\end{definition}

\begin{remark}
Let $\Mcal$ be a $\LM$-operad. Let $\Acal$ be an algebroid in $\Mcal_l$ with category of objects $X$. A left $\Acal$-module $\Pcal$ assigns to each object $x$ in $X$ an object $\Pcal(x)$  in $\Mcal_m$. For every $n \geq 0$ and every sequence of objects and  arrows 
\[
y_0 = x_0 \leftarrow y_1,  x_1 \leftarrow y_2, \ldots, x_{n-1} \leftarrow y_{n} , x_{n} \leftarrow y_{n+1}
\]
in $X$, the left $\Acal$- module $\Pcal$ induces a multimorphism
\[
\lbrace \Acal( y_1, x_1 ),\ldots,  \Acal(y_n, x_n), \Pcal(y_{n+1}) \rbrace \rightarrow \Pcal(y_0) 
\]
in $\Mcal$. In the case when $\Mcal$ is a $\LM$-monoidal category (in other words, $\Mcal_l$ is a monoidal category and $\Mcal_m$ is a left module for it), this induces  a morphism
\[
 \Acal( y_1, x_1 ) \otimes \ldots \otimes \Acal(y_n, x_n) \otimes \Pcal(y_{n+1}) \rightarrow \Pcal(y_0) .
\]
In particular, in the case when $n = 1$ and the arrows are identities we obtain, for every pair of objects $x_0, x_1$ in $X $ a map
\[
\Acal(x_0, x_1) \otimes \Pcal(x_1) \rightarrow \Pcal(x_0).
\]
This is compatible with the units and composition of $\Acal$, up to homotopy.
\end{remark}

\begin{example}
Let $X$ be a category and  let $x$ be an object of $X$. Then the functor $(\id, x): X \sqcup [0] \rightarrow X$ induces a morphism of associative operads $\LM_X \rightarrow \Assos_X$. It follows that for every associative operad $\Mcal$ and every $\Mcal$-algebroid $\Acal$ with category of objects $X$, we have an induced left module $\Pcal$ in $\Mcal$. This has the following properties:
\begin{itemize}
\item For every object $x'$ in $X$ we have $\Pcal(x') = \Acal(x', x)$.
\item For every pair of objects $x', x''$ in $X$ the action map
\[
\lbrace \Acal(x'', x') , \Pcal(x') \rbrace \rightarrow \Pcal(x'')
\]
is equivalent, under the identifications of the previous item, to the composition map
\[
\lbrace \Acal(x'', x') , \Acal(x', x) \rbrace \rightarrow \Acal(x'', x).
\]
\end{itemize}
We call $\Pcal$ the left module corepresented by $x$. We will usually use the notation $\Acal(-, x)$ for $\Pcal$, and in the case when $\Acal$ is an $\Mcal$-enriched category, we instead write $\Hom_{\Acal}(x, -)$.
\end{example}

\begin{construction}\label{construction functoriality lmod}
Consider the functor $\Alg_{\LM_{-}}(-)$ defined by the composition
\[
\Cat^\op \times \Op_{\LM} \xrightarrow{\LM_{-} \times \id_{\Op_{\LM}}} \Op_{\LM}^\op \times \Op_{\LM} \xrightarrow{\Alg_-(-)} \Cat.
\]
For each object $\Mcal$ in $\Op_{\LM}$ we denote by $\LMod(\Mcal)$ the total category of the cartesian fibration associated to the functor $\Alg_{\LM_{-}}(\Mcal):  \Cat^\op \rightarrow \Cat$. This comes equipped with a forgetful functor $\LMod(\Mcal) \rightarrow \Algbrd(\Mcal_l) $. For each algebroid $\Acal$ in $\Mcal_l$ we denote by $\LMod_\Acal(\Mcal)$  the fiber  of $\LMod(\Mcal)$ over $\Acal$, and call it the category of left $\Acal$-modules.

The assignment $\Mcal \mapsto \LMod(\Mcal)$ defines a functor $\LMod(-): \Op_{\BM} \rightarrow \widehat{\Cat}$. Let $\LMod$ be the total category of the associated cocartesian fibration. In other words, $\LMod$ is the total category of the two-sided fibration associated to $\Alg_{\LM_-}(-)$.
\end{construction}

\begin{warning}
Our usage of the terminology $\LMod(\Mcal)$ conflicts with that of \cite{HA}. There only left modules over associative algebras are considered - this corresponds to the fiber of the projection $\LMod(\Mcal) \rightarrow \Cat$ over $[0]$.
\end{warning}

\begin{remark}\label{remark lmod fits into}
The category $\LMod$  fits into a commutative square
\[
\begin{tikzcd}
\operatorname{Arr}_{\text{oplax}}(\Cat) \arrow{d}{} & \arrow{l}{} \LMod \arrow{r}{} \arrow{d}{} & \Algbrd \arrow{d}{} \\
 \Cat \times \Cat & \arrow{l}{} \Cat \times \Op_{\LM} \arrow{r}{} & \Cat \times \Op_{\Assos} .
\end{tikzcd}
\]
Here the vertical arrows are the two-sided fibrations classified by the functors $\Funct(-, -)$,  $\Alg_{\LM_{-}}(-)$ and $\Algbrd_-(-)$, and the horizontal arrows are the functors of ``forgetting the algebra'' and  ``forgetting the left module".
\end{remark}

\begin{proposition}\label{prop restrict scalars lmod}
Let $\Mcal$ be a $\LM$-operad. Then the projection $p: \LMod(\Mcal) \rightarrow \Algbrd(\Mcal_l)$ is a cartesian fibration. Moreover, a morphism $F: (\Acal, \Pcal) \rightarrow (\Bcal, \Qcal)$ in $\LMod(\Mcal)$ is $p$-cartesian if and only if for every object $x$ in $\Acal$ the induced map $\Pcal(x) \rightarrow \Qcal(F(x))$ is an isomorphism.
\end{proposition}
\begin{proof}
Let $\Env(\Mcal)$ be the $\LM$-monoidal envelope of $\Mcal$, and let $\Pcal(\Env(\Mcal))$ be the image of $\Env(\Mcal)$ under the symmetric monoidal functor $\Pcal : \Cat \rightarrow \Pr^L$. We have a commutative square of categories
\[
\begin{tikzcd}
\LMod(\Mcal) \arrow{d}{} \arrow{r}{} & \LMod(\Pcal(\Env(\Mcal))) \arrow{d}{} \\
\Algbrd(\Mcal_l) \arrow{r}{} & \Algbrd(\Pcal(\Env(\Mcal))_l).
\end{tikzcd}
\]
Note that the horizontal arrows are fully faithful. Our result would follow if we are able to show that the right vertical arrow is a cartesian fibration, and that cartesian morphisms are given by the condition in the statement. In other words, it suffices to prove the  result in the case when $\Mcal$ is a presentable $\LM$-monoidal category. We assume that this is the case from now on.

Let $X$ be a category. Recall from \cite{Hinich} that $\LM_X$ is a flat $\LM$-operad. It follows that there is a universal $\LM$-operad $\Mcal_X$ equipped with a morphism of $\LM$-operads  
\[
\Mcal_X \times_{\LM} \LM_X \rightarrow \Mcal.
\] 
In particular, we have equivalences
\[
\Alg_{\LM_X}(\Mcal) = \Alg_{\LM}(\Mcal_X)
\]
and
\[
\Algbrd_X(\Mcal) = \Alg_{\Assos}((\Mcal_X)_l).
\]
The projection $p_X : \Alg_{\LM_X}(\Mcal) \rightarrow \Algbrd_X(\Mcal_l)$ becomes identified, under this dictionary, with the canonical projection
\[
p'_X: \Alg_{\LM}(\Mcal_X) \rightarrow \Alg_{\Assos}((\Mcal_X)_l).
\]

As discussed in \cite{Hinich} corollary 4.4.9, $\Mcal_X$ is a presentable $\LM$-monoidal category. It now follows from \cite{HA} corollary 4.2.3.2 that $p_X$ is a cartesian fibration, and moreover a morphism $F: (\Acal, \Pcal) \rightarrow (\Bcal, \Qcal)$ in $\Alg_{\LM_X}$ is $p_X$-cartesian if and only if for every object $x$ in $X$ the induced map $\Pcal(x) \rightarrow \Qcal(x)$ is an isomorphism.

Assume now that $F$ is $p_X$-cartesian and let $g: Y \rightarrow X$ be a functor of categories. Consider the induced morphism
\[
g^!F: (g^!\Acal, g^!\Pcal) \rightarrow (g^!\Bcal, g^!\Qcal)
\]
in $\Alg_{\LM_Y}(\Mcal)$. Let $y$ be an object in $Y$. Then the induced map $g^!\Pcal(y) \rightarrow g^!\Qcal(y)$ is equivalent to the map $\Pcal(g(y)) \rightarrow \Qcal(g(y))$, and is therefore an isomorphism. It follows that $g^!F$ is also $p_Y$-cartesian, and hence $g^!$ is a morphism of cartesian fibrations. Combining \cite{HTT} propositions 2.4.2.8 and 2.4.2.11 we conclude that $p$ is a cartesian fibration.

Our characterization of $p$-cartesian morphisms follows from the above characterization of $p_X$-cartesian morphisms together with item (iii) in \cite{HTT} proposition 2.4.2.11.
\end{proof}

\begin{notation}
Let $\Mcal$ be an associative operad. We denote by $(\Op_{\LM})|_{\Mcal}$ the fiber of the projection $\Op_{\LM} \rightarrow \Op_{\Assos}$ over $\Mcal$, and by $\LMod|_{\Mcal}$ the fiber over $\Mcal$ of the projection $\LMod \rightarrow \Op_{\Assos}$. 
\end{notation}
\begin{corollary}\label{coro lmod is two sided}
Let $\Mcal$ be an associative operad. Then the projection 
\[
\LMod|_{\Mcal} \rightarrow (\Op_{\LM})|_{\Mcal} \times \Algbrd(\Mcal)
\]
 is a two-sided fibration from $(\Op_{\LM})|_{\Mcal}$ to $\Algbrd(\Mcal)$.
\end{corollary}
\begin{proof}
By construction, the projection $\LMod \rightarrow \Algbrd$ is a morphism of cocartesian fibrations over the functor $\Op_{\LM} \rightarrow \Op_{\Assos}$. It follows that the projection in the statement is a morphism of cocartesian fibrations over $(\Op_{\LM})|_{\Mcal}$. Its fiber over a given $\Mcal$-module is a cartesian fibration, thanks to proposition \ref{prop restrict scalars lmod}. Our claim now follows from \cite{HSTI} proposition 1.1.9.
\end{proof}

\begin{proposition}\label{prop lmod tiene prods}
The categories  $\LMod$ and $\operatorname{Arr}_{\text{\normalfont oplax}}(\Cat)$ admit finite products. Moreover, all the maps in the diagram of remark \ref{remark lmod fits into} preserve finite products.
\end{proposition}
\begin{proof}
The fact that $\LMod$ admits finite products which are preserved by the projection to $\Cat \times \Op_{\LM}$ follows by the same arguments as those from proposition \ref{prop products of algebroids}. One similarly shows that $\operatorname{Arr}_{\text{oplax}}(\Cat)$ admits finite products which are preserved by the projection to $\Cat \times \Cat$. It remains to show that the projections from $\LMod$ to $\operatorname{Arr}_{\text{oplax}}(\Cat)$ and $\Algbrd$ preserve finite products. Both claims can be proven using similar arguments - below we present the case of  $\Algbrd$.

Observe first that the final object for $\LMod$ is the unique object lying above the final object in $\Cat \times \Op_{\LM}$. Its image in $\Algbrd$ is the unique algebroid lying above the final object in $\Cat \times \Op_{\Assos}$, which is indeed the final object of $\Algbrd$.

 It remains to show that the projection $\LMod \rightarrow \Algbrd$ preserves binary products. Let $\Mcal, \Ncal$ be two $\LM$-operads, and let $X, Y$ be two categories. Let $(\Acal, \Pcal)$ and $(\Bcal, \Qcal)$  be objects of $\LMod$ lying above $(X, \Mcal)$  and $(Y, \Ncal)$, respectively. A variant of the discussion from remark \ref{remark fact proj} shows that their product $(\Acal, \Pcal) \boxtimes (\Bcal, \Qcal)$ fits into a diagram
\[
(\Acal, \Pcal) \xleftarrow{ \beta_{(\Acal,\Pcal)} } \overline{(\Acal, \Pcal)}   \xleftarrow{\alpha_{(\Acal,\Pcal)}} (\Acal, \Pcal) \boxtimes (\Bcal, \Qcal) \xrightarrow{\alpha_{(\Bcal,\Qcal)}} \overline{(\Bcal, \Qcal)} \xrightarrow{\beta_{(\Bcal,\Qcal)}} (\Bcal, \Qcal)
\]
where $\alpha_{{(\Acal,\Pcal)}}$ and $\alpha_{{(\Bcal,\Qcal)}}$ are cocartesian for the projection $\LMod \rightarrow \Op_{\LM}$ and $\beta_{{(\Acal,\Pcal)}}$ and $\beta_{{(\Bcal,\Qcal)}}$ are cartesian for the projection $\LMod \rightarrow \Cat$. The image of the above diagram under the projection to $\Algbrd$ recovers a diagram
\[
\Acal \xleftarrow{\beta_\Acal} \overline{\Acal} \xleftarrow{\alpha_\Acal} W \xrightarrow{\alpha_{\Bcal}} \overline{\Bcal} \xrightarrow{\beta_{\Bcal}} \Bcal
\]
where $\alpha_\Acal$ and $\alpha_{\Bcal}$ are cocartesian for the projection $\Algbrd \rightarrow \Op_{\Assos}$ and $\beta_{\Acal}$ and $\beta_{\Bcal}$ are cartesian for the projection $\Algbrd \rightarrow \Cat$. Using remark \ref{remark fact proj} we conclude that the above diagram exhibits $W$ as the product of $\Acal$ and $\Bcal$ in $\Algbrd$, as desired.
\end{proof}

\begin{proposition}\label{prop Lmod products comp with}
Let $f: (\Acal, \Pcal) \rightarrow (\Acal', \Pcal')$ be a morphism in $\LMod$ and let $(\Bcal, \Qcal)$ be another object of $\LMod$.  Denote by $p = (p_1, p_2)$ the projection $\LMod \rightarrow \Cat \times \Op_{\LM}$.
\begin{enumerate}[\normalfont (i)]
\item If $f$ is $p_1$-cartesian then $f \boxtimes \id_{\Bcal}$ is $p_1$-cartesian.
\item If $f$ is $p_2$-cocartesian then $f \boxtimes \id_{\Bcal}$ is $p_2$-cocartesian.
\end{enumerate}
\end{proposition}
\begin{proof}
Follows from the same arguments as those of proposition \ref{prop products compatible with }.
\end{proof}

\subsection{Enrichment of presentable modules} \label{subsection enrichement} Our next goal is to discuss the procedure of enrichment of presentable modules over presentable monoidal categories.

\begin{notation}
Recall the projection $\LMod \rightarrow \operatorname{Arr}_{\text{oplax}}(\Cat)$ from remark \ref{remark lmod fits into}. Note that we have an inclusion $\Funct([1], \Cat) \rightarrow \operatorname{Arr}_{\text{oplax}}(\Cat)$ which is surjective on objects, which arises from straightening the natural transformation $\Hom_{\Cat}(-,-) \rightarrow \Funct(-,-)$ (see \cite{HSTI} proposition 2.2.4). In other words, $\Funct([1], \Cat)$ is the total category of the maximal bifibration contained inside the two-sided fibration $\operatorname{Arr}_{\text{oplax}}(\Cat) \rightarrow \Cat \times \Cat$. 

We denote by $\LMod'$ the fiber product $\LMod \times_{\operatorname{Arr}_{\text{oplax}}} \Funct([1],\Cat)$. For each $\Mcal$ in $\Op_{\LM}$ we denote by $\LMod'(\Mcal)$ the fiber over $\Mcal$ of the projection $\LMod' \rightarrow \Op_{\LM}$.
\end{notation}

\begin{remark}\label{remark prods in LModprime}
A variation of the argument in \ref{prop lmod tiene prods} shows that the inclusion of $\Funct([1],\Cat)$ inside $\operatorname{Arr}_{\text{oplax}}(\Cat)$ preserves finite products. It follows that $\LMod'$ admits finite products, which are preserved by its inclusion inside $\LMod$.
\end{remark}

\begin{proposition}\label{prop existence final}
Let $\Mcal$ be a presentable $\LM$-monoidal category (in other words, a pair of a monoidal category $\Mcal_l$ and a presentable module $\Mcal_m$). Then $\widehat{\LMod'}(\Mcal)$ has a final object. Moreover, a pair $(\Acal, \Pcal)$ of an $\Mcal_l$-algebroid $\Acal$ with category of objects $X$ and a left $\Acal$-module $\Pcal$ in $\Mcal_m$ is final if and only if the functor $X \rightarrow \Mcal_m$ underlying $\Pcal$ is an equivalence, and for every pair of objects $x, y$ in $X$, the action map
\[
\Acal(y, x) \otimes \Pcal(x) \rightarrow \Pcal(y)
\]
exhibits $\Acal(y, x)$ as the Hom object between $x$ and $y$.
\end{proposition}
\begin{proof}
Note that the composition of the projection $\widehat{\LMod'}(\Mcal) \rightarrow \Funct([1],\widehat{\Cat})$ with the target map $\Funct([1],\widehat{\Cat}) \rightarrow \widehat{\Cat}$ is canonically equivalent to the constant functor $\Mcal_m$. It follows that we have a commutative diagram
\[
\begin{tikzcd}
\widehat{\Cat}_{/\Mcal_m} \arrow{dr}{} & & \arrow{ll}{} \widehat{\LMod'}(\Mcal) \arrow{dl}{} \\
& \widehat{\Cat} & 
\end{tikzcd}
\]
where the left vertical map is the forgetful functor, and the right vertical arrow picks out the category of objects of the underlying algebroid. Since the right vertical arrow is a cartesian fibration and the left vertical arrow is a right fibration, we have that the horizontal arrow is a cartesian fibration.

Let $X$ be a category equipped with a map $f: X \rightarrow \Mcal_m$ and recall the presentable $\LM$-monoidal category $\Mcal_X$ from the proof of proposition \ref{prop restrict scalars lmod}. The fiber of $\widehat{\LMod'}(\Mcal)$ over $X$ is the category of $\LM_X$-algebras in $\Mcal_m$ whose underlying functor $X \rightarrow \Mcal_m$ is $f$. This can equivalently be described as the category of associative algebras in $(\Mcal_X)_l$ equipped with an action on $f$ (thought of as an object of $(\Mcal_X)_m$). By \cite{Hinich} corollary 6.3.4, we conclude that $(\widehat{\LMod'}(\Mcal))_X$ admits a final object, and moreover a pair $(\Acal, \Pcal)$ of an $\Mcal_l$-algebroid with category of objects $X$ and a left $\Acal$-module $\Pcal$ whose underlying functor $X \rightarrow \Mcal_m$ is $f$  is final if and only if for every pair of objects $x, y$ in $X$ the action map 
\[
\Acal(y, x) \otimes f(x) \rightarrow f(y)
\]
exhibits $\Acal(y, x)$ as the Hom object between $f(x)$ and $f(y)$.

This description implies that if $(\Acal, \Pcal)$ is final in $(\widehat{\LMod'}(\Mcal))_X$ and $g: Y \rightarrow X$ is a functor of categories, then $(g^!\Acal, g^!\Pcal)$ is final in $(\widehat{\LMod'}(\Mcal))_Y$. The result now follows from an application of \cite{HTT} proposition 4.3.1.10.
\end{proof}

\begin{corollary}\label{coro repr right fibr}
Let $\Mcal$ be a presentable $\LM$-monoidal category. Then the projection $\widehat{\LMod'}(\Mcal) \rightarrow \widehat{\Algbrd}(\Mcal_l)$ is a representable right fibration.
\end{corollary}
\begin{proof}
Combine propositions \ref{prop restrict scalars lmod} and \ref{prop existence final}.
\end{proof}

\begin{notation}
Let $\Mcal$ be a presentable monoidal category. We denote by $\widehat{\LMod'}|_{\Mcal\modd(\Pr^L)}$ the base change of the projection $\widehat{\LMod'} \rightarrow \widehat{\Op}_{\LM}$ along the inclusion $\Mcal\modd(\Pr^L) \rightarrow \widehat{\Op}_{\LM}$.
\end{notation}

For our next result, we need the notion of representable bifibration (see for instance \cite{HSTI} definition 2.2.7).
\begin{corollary}\label{coro Lmodprime repr bifibr}
Let $\Mcal$ be a presentable monoidal category. Then the projection 
\[
\widehat{\LMod'}|_{\Mcal\modd(\Pr^L)} \rightarrow \Mcal \modd(\Pr^L) \times \widehat{\Algbrd}(\Mcal)
\]
is a representable bifibration from $\Mcal\modd(\Pr^L)$ to $\widehat{\Algbrd}(\Mcal)$.
\end{corollary}
\begin{proof}
It follows from corollary \ref{coro lmod is two sided} together with the description of cartesian arrows from proposition \ref{prop restrict scalars lmod} that the projection in the statement is the maximal bifibration contained inside the base change of the projection from corollary \ref{coro lmod is two sided} along the inclusion
\[
\Mcal \modd(\Pr^L) \times \widehat{\Algbrd}(\Mcal) \rightarrow (\widehat{\Op}_{\LM})|_{\Mcal} \times \widehat{\Algbrd}(\Mcal).
\]
The fact that it is representable is the content of corollary \ref{coro repr right fibr}.
\end{proof}

Our next goal is to study the dependence in $\Mcal_l$ of the algebroid from proposition \ref{prop existence final}.

\begin{construction}\label{construction thetaprime}
Consider the commutative diagram of categories
\[
\begin{tikzcd}
\widehat{\LMod'}|_{\Alg_{\LM}(\Pr^L)} \arrow{r}{i'} \arrow{d}{p'} & \widehat{\LMod'} \arrow{d}{p} \arrow{r}{r} & \widehat{\Algbrd}  \arrow{d}{q}
\\
\Alg_{\LM}(\Pr^L) \arrow{r}{i} & \widehat{\Op_{\LM}} \arrow{r}{j} & \widehat{\Op_{\Assos}}
\end{tikzcd}
\]
where the left square is cartesian.  We equip all four categories in the right square with their cartesian symmetric monoidal structure. By propositions \ref{prop products of algebroids} and \ref{prop lmod tiene prods} together with remark \ref{remark prods in LModprime} we see that the right square  has a   canonical lift to a commutative square  of cartesian symmetric monoidal categories. It follows from propositions \ref{prop products compatible with }  and \ref{prop Lmod products comp with}  that $p$ and $q$ are in fact cocartesian fibrations of operads.

Equip $\Alg_{\LM}(\Pr^L)$ with its canonical symmetric monoidal structure, so that $i$  inherits a lax symmetric monoidal structure. It follows from the above that $p'$ has a canonical structure of cocartesian fibration of operads, and $i'$ of lax symmetric monoidal morphism.
 
 Using  proposition \ref{prop existence final} we see that that $p'$ admits a fully faithful right adjoint $p'^R$, which comes equipped with a canonical lax symmetric monoidal structure. We therefore have a commutative diagram of symmetric monoidal categories and lax symmetric monoidal functors as follows:
 \[
 \begin{tikzcd}
 & \widehat{\Algbrd} \arrow{d}{q}
 \\
 \Alg_{\LM}(\Pr^L) \arrow{r}{ji} \arrow{ur}{r i' p'^R} & \widehat{\Op_{\Assos}}
 \end{tikzcd}
 \]
 Note that $ji$ factors through the lax symmetric monoidal inclusion $\Alg_{\Assos}(\Pr^L) \rightarrow \widehat{\Op_{\Assos}}$, so we have an induced diagram of symmetric monoidal categories and lax symmetric monoidal functors
\[
\begin{tikzcd}
 & \widehat{\Algbrd}|_{\Alg_{\Assos}(\Pr^L)} \arrow{d}{q'}
 \\
 \Alg_{\LM}(\Pr^L) \arrow{r}{u} \arrow{ur}{\theta'} & \Alg_{\Assos}(\Pr^L).
 \end{tikzcd}
\]

 Observe that the maps $u$ and $q'$ are cocartesian fibrations of operads. If $\Mcal_l$ is a presentable symmetric monoidal category, thought of as a commutative algebra object in $\Alg_{\Assos}(\Pr^L)$, we obtain in particular an induced lax symmetric monoidal functor
\[
\theta'_\Mcal: \Mcal \modd(\Pr^L) \rightarrow \widehat{\Algbrd}(\Mcal).
\]
\end{construction}

\begin{proposition}\label{prop thetaprime pres cartesian}
The map $\theta'$ from construction \ref{construction thetaprime} is a morphism of cartesian fibrations over $\Alg_{\Assos}(\Pr^L)$.
\end{proposition}
\begin{proof}
We continue with the notation from construction \ref{construction thetaprime}. Observe that $u$ and $q'$ are indeed cartesian fibrations. For $q'$ this follows from the adjoint functor theorem combined with remark \ref{some infty two functoriality}, and for $u$ this follows from \cite{HA} corollary 4.2.3.2. The rest of the proof is devoted to showing that $\theta'$ maps $u$-cartesian arrows to $q'$-cartesian arrows.

Let $F: \Mcal \rightarrow \Mcal'$ be a  $u$-cartesian arrow in $\Alg_{\LM}(\Pr^L)$, whose components consist of a morphism of presentable monoidal categories $F_l : \Mcal_l \rightarrow \Mcal'_l$ and an isomorphism of modules $F_m : \Mcal_m \rightarrow \Mcal_m'$. Let $(\Acal, \Pcal) = p'^R(\Mcal)$ and $(\Acal', \Pcal') = p'^R(\Mcal')$. The morphism 
\[
p'^R F : (\Acal, \Pcal) \rightarrow (\Acal', \Pcal')
\]
 can be factored as $\eta \alpha$ where $\alpha$ is a $p'$-cocartesian lift of $F$, and $\eta: F_! (\Acal, \Pcal) \rightarrow (\Acal', \Pcal')$ is the unique map in $\widehat{\LMod'}(\Mcal')$ from $F_!(\Acal, \Pcal)$ to the final object. Consider the morphism of algebroids $r i' \eta : (F_l)_!\Acal \rightarrow \Acal'$. We have to show that the induced map $\mu:\Acal \rightarrow  (F^R_l)_!\Acal'$ is an isomorphism. 

Let $X$ and $X'$ be the categories of objects of $\Acal$ and $\Acal'$, respectively. We have a commutative square
\[
\begin{tikzcd}
\Mcal_m \arrow{r}{F_m} & \Mcal'_m \\
X \arrow{u}{\Pcal} \arrow{r}{g} & X' \arrow{u}{\Pcal'}
\end{tikzcd}
\]
where the bottom horizontal arrow is the projection of $p'^RF$ to $\Cat$. It follows from proposition \ref{prop existence final} that the vertical arrows are equivalences. Since $F_m$ is also an isomorphism, we conclude that $g$ is an equivalence. It follows that $\eta$ and $\mu$ also induce equivalences at the level of categories of objects.
To simplify notation, in the rest of the proof we identify $X$ with $\Mcal_m$ and $X'$ with $\Mcal'_m$ via the isomorphisms $\Pcal$ and $\Pcal'$.

Let $x, y$ be two objects in $\Mcal_m$. We have to show that the induced map 
\[
\mu_*: \Acal(y, x) \rightarrow F_l^R \Acal'(F_m y, F_m x)
\] 
is an isomorphism. The morphism $\eta$ induces a commutative square in $\Mcal_m'$ as follows:
\[
\begin{tikzcd}
F_l\Acal (y, x) \otimes F_m x \arrow{r}{}  \arrow{d}{ri' \eta_* \otimes \id} & F_m y \arrow{d}{\id}\\
\Acal'(F_m y, F_m x) \otimes F_m x\arrow{r}{} & F_m y
\end{tikzcd}
\]
Applying the (lax symmetric monoidal) right adjoint to $F$ yields a commutative square in $\Mcal_m$ as follows:
\[
\begin{tikzcd}
F_l^R F_l\Acal (y, x) \otimes x \arrow{r}{}  \arrow{d}{F_l^R ri' \eta_* \otimes \id} & y \arrow{d}{\id }\\
F_l^R\Acal'(F_m y, F_m x) \otimes x \arrow{r}{} & y
\end{tikzcd}
\]

Composing with the unit map $\Acal(y, x) \rightarrow F_l^R F_l \Acal(y, x)$ yields a commutative square
\[
\begin{tikzcd}
\Acal (y, x) \otimes x \arrow{r}{}  \arrow{d}{\mu_* \otimes \id} & y \arrow{d}{\id }\\
F_l^R\Acal'(F_m y, F_m x) \otimes x \arrow{r}{} & y
\end{tikzcd}
\]
where  the top horizontal arrow exhibits $\Acal(y, x)$ as the Hom object between $x$ and $y$.

It now suffices to show that the bottom horizontal arrow exhibits $F^R_l\Acal'(F_m y, F_m x)$ as the Hom object between $x$ and $y$. This follows  from the right adjointability of the following commutative square of categories:
\[
\begin{tikzcd}[column sep = large]
\Mcal_l \arrow{d}{F_l} \arrow{r}{-\otimes x} & \Mcal_m \arrow{d}{F_m} \\
\Mcal_l' \arrow{r}{- \otimes F_m x} & \Mcal_m'
\end{tikzcd}
\]
\end{proof}

\begin{corollary}\label{coro restriction}
Let $F: \Mcal \rightarrow \Mcal'$ be a colimit preserving symmetric monoidal functor between presentable symmetric monoidal categories. Then there is a commutative square of symmetric monoidal categories and lax symmetric monoidal functors
\[
\begin{tikzcd}
\Mcal'\modd(\Pr^L) \arrow{r}{\theta'_{\Mcal'}} \arrow{d}{F^*} & \widehat{\Algbrd}(\Mcal') \arrow{d}{F^R_!} \\
\Mcal\modd(\Pr^L) \arrow{r}{\theta'_{\Mcal}} & \widehat{\Algbrd}(\Mcal)
\end{tikzcd}
\]
where $F^*$ denotes the functor of restriction of scalars along $F$.
\end{corollary}
\begin{proof}
For each symmetric monoidal category $X$ denote by $X^\otimes$ its category of operators. The map $F$ corresponds to a morphism $F^\otimes: [1] \times \Fin_* \rightarrow \Alg_{\Assos}(\Pr^L)^\otimes$ over $\Fin_*$. Base change of $\theta'$ along $F^\otimes$ yields a commutative diagram
\[
\begin{tikzcd}
& (\widehat{\Algbrd}|_{\Alg_{\Assos}(\Pr^L)})^\otimes_{F^\otimes} \arrow{d}{v=(v_1, v_2)} \\
\Alg_{\LM}(\Pr^L)^\otimes_{F^\otimes}\arrow{r}{h = (h_1, h_2)} \arrow{ur}{\theta'_{F^\otimes}} & {[1]} \times \Fin_*
\end{tikzcd}
\]
where $h$ and $v$ are cocartesian fibrations. Observe that the maps $h$ and $v$ are also two-sided fibrations - in other words, the associated functors $[1] \times \Fin_* \rightarrow \Cat$ are right adjointable in the $[1]$ coordinate. In particular, we see that $h_1$ and $v_1$ are cartesian fibrations, and $h$ and $v$ are morphisms of cartesian fibrations over $[1]$. It follows from proposition \ref{prop thetaprime pres cartesian} that $\theta'_{F^\otimes}$ is a morphism of cartesian fibrations over $[1]$. Straightening it yields a commutative square in $\widehat{\Cat}_{/\Fin_*}$.  Tracing the construction of this square reveals that it is actually a commutative square of commutative operads, and satisfies the desired conditions.
\end{proof}

\begin{proposition}\label{prop el caso de spc}
The lax symmetric monoidal functor
\[
\theta'_{\Spc} : \Pr^L = \Spc\modd(\Pr^L) \rightarrow \widehat{\Algbrd}(\Spc)
\]
factors through the image of the section $s$ from construction \ref{construction inclusion cat}. Furthermore, the composition of $\theta'_{\Spc}$ with the symmetric monoidal projection  $\widehat{\Algbrd}(\Spc) \rightarrow \widehat{\Cat}$ is equivalent to the forgetful functor $\Pr^L \rightarrow \widehat{\Cat}$ with its canonical lax symmetric monoidal structure.
\end{proposition}
\begin{proof}
Let $\Ccal$ be a presentable category. Then $\theta'_{\Spc}(\Ccal)$ is a $\Spc$-algebroid  equipped with a left module in $\Ccal$ which identifies its category of objects with $\Ccal$, and such that for every pair of objects $x, y$ the map
\[
\theta'_{\Spc}(\Ccal)(y, x) \otimes x \rightarrow y
\]
exhibits $\theta'_{\Spc}(\Ccal)(y, x)$ as the Hom object between $x$ and $y$. Inspecting the unit map $\Hom_\Ccal(x, y) \rightarrow \theta'_{\Spc}(\Ccal)(y, x)$ one obtains a commutative square
\[
\begin{tikzcd}
\Hom_\Ccal(x, y) \otimes x \arrow{r}{\ev} \arrow{d}{} & y \arrow{d}{\id} \\ 
\theta'_{\Spc}(\Ccal)(y, x) \otimes x \arrow{r}{} & y.
\end{tikzcd}
\]
The top horizontal arrow exhibits $\Hom_{\Ccal}(x, y)$ as the Hom object between $x$ and $y$, so we see that the unit map $\Hom_{\Ccal}(x, y) \rightarrow \theta'_{\Spc}(\Ccal)(y , x)$ is an isomorphism. As discussed in the proof of lemma \ref{lemma left adjoint to s}, this is also the case for the algebroid $s(\Ccal)$. It follows that the canonical map $s(\Ccal) \rightarrow \theta'_{\Spc}(\Ccal)$ is an equivalence, and therefore $\theta'_{\Spc}$ factors through the image of $s$.

Consider now the following diagram:
\[
\begin{tikzcd}
\Funct([1],\Cat) \arrow{d}{(t_0, t_1)} & \arrow{l}[swap]{v} \widehat{\LMod'} \arrow{r}{r} \arrow{d}{(\alpha, p)} & \widehat{\Algbrd} \arrow{d}{(\beta, q)} \\
\widehat{ \Cat} \times \widehat{\Cat} & \arrow{l}[swap]{(\id, w)} \widehat{\Cat} \times \widehat{\Op_{\LM}} \arrow{r}{(\id, j)} & \widehat{\Cat} \times \widehat{\Op_{\Assos}} 
\end{tikzcd}
\]
We equip all categories above with their cartesian symmetric monoidal structure, so that all functors inherit a canonical symmetric monoidal structure.

The composition of $\theta'_{\Spc}$ with the symmetric monoidal projection $\widehat{\Algbrd}(\Spc) \rightarrow \widehat{\Cat}$ is equivalent to the lax symmetric monoidal functor obtained by taking the fiber over $\Spc$ of the lax symmetric monoidal functor $ (\beta,q) r i' p'^R$. This is equivalent to the composite lax symmetric monoidal functor
\[
\Pr^L = \Spc\modd(\Pr^L) \hookrightarrow \Alg_{\LM}(\Pr^L) \xrightarrow{p'^R} \widehat{\LMod'}_{\Alg(\Pr^L)} \xrightarrow{i'} \widehat{\LMod'} \xrightarrow{r} \widehat{\Algbrd} \xrightarrow{\beta} \widehat{\Cat}
\]
which is in turn equivalent to the following composition:
\[
\Pr^L = \Spc\modd(\Pr^L) \hookrightarrow \Alg_{\LM}(\Pr^L) \xrightarrow{p'^R} \widehat{\LMod'}_{\Alg(\Pr^L)} \xrightarrow{i'} \widehat{\LMod'} \xrightarrow{v} \Funct([1],\Cat) \xrightarrow{t_0} \widehat{\Cat}
\]
Meanwhile, the  lax symmetric monoidal forgetful functor $\Pr^L \rightarrow \widehat{\Cat}$ can be obtained as the following composition:
\[
\Pr^L = \Spc\modd(\Pr^L) \hookrightarrow \Alg_{\LM}(\Pr^L) \xrightarrow{p'^R} \widehat{\LMod'}_{\Alg(\Pr^L)} \xrightarrow{i'} \widehat{\LMod'} \xrightarrow{v} \Funct([1],\Cat) \xrightarrow{t_1} \widehat{\Cat}.
\]
We have to show that these agree. Note that they are both obtained by composing the lax symmetric monoidal functor  
\[
F: \Pr^L = \Spc\modd(\Pr^L) \hookrightarrow \Alg_{\LM}(\Pr^L) \xrightarrow{p'^R} \widehat{\LMod'}_{\Alg(\Pr^L)} \xrightarrow{i'} \widehat{\LMod'} \xrightarrow{v}  \Funct([1],\Cat)
\]
with either $t_0$ or $t_1$.  However, thanks to the characterization of the image of $p'^R$ from proposition \ref{prop existence final}, we have that the image of $F$ belongs to the full subcategory $\Funct([1], \Cat)_{\text{iso}}$ of $\Funct([1],\Cat)$ on the isomorphisms. Hence we can factor the lax symmetric monoidal functor $F$ as follows:
\[
\Pr^L \xrightarrow{F'} \Funct([1], \Cat)_{\text{iso}} \hookrightarrow \Funct([1],\Cat).
\] 
Our claim now follows from the fact that the restrictions of $s_0$ and $s_1$ to $\Funct([1],\Cat)_{\text{iso}}$ are equivalent.
\end{proof}

\begin{corollary}\label{prop is enriched}
Let $\Mcal$ be a presentable $\LM$-monoidal category, and let $\gamma: \Mcal_m^{\leq 0} \rightarrow \Mcal_m$ be the inclusion. Then $\gamma^! \theta'(\Mcal)$ is an $\Mcal_l$-enriched category.
\end{corollary}
\begin{proof}
We have an equivalence $(\tau_{\Mcal})_! \gamma^! \theta'(\Mcal) = \gamma^! (\tau_{\Mcal})_! \theta'(\Mcal)$. Thanks to corollary \ref{coro restriction} we have that $(\tau_\Mcal)_!(\theta'(\Mcal))$ is equivalent to $\theta'(\Ncal)$, where $\Ncal$ is the presentable $\LM$-monoidal category obtained by restricting the action of $\Mcal_l$ on $\Mcal_m$ along the unit map $\Spc \rightarrow \Mcal_l$. Our claim now follows directly from proposition \ref{prop el caso de spc}.
\end{proof}

We now construct a variant of the functor $\theta'$ which takes values in enriched categories.
\begin{construction}\label{construction theta}
We continue with the notation of construction \ref{construction thetaprime}. Consider the following commutative diagram:
\[
 \begin{tikzcd}
 & \widehat{\Algbrd} \arrow{d}{q}  & \widehat{\Algbrd}_{\Spc} \arrow{l}[swap]{h} \arrow{d}{q_{\Spc}} 
 \\
 \Alg_{\LM}(\Pr^L) \arrow{r}{ji} \arrow{ur}{r i' p'^R} & \widehat{\Op_{\Assos}} & \arrow{l}[swap]{\id} \widehat{\Op_{\Assos}}
 \end{tikzcd}
 \]
  Here $\widehat{\Algbrd}_{\Spc}$ is the full subcategory of $\widehat{\Algbrd}$ on those algebroids with a space of objects. We equip all categories on the right square with their cartesian symmetric monoidal structure. Note that $q$ and $q_{\Spc}$ have canonical structures of cocartesian fibrations of operads, and $h$ is a morphism of cocartesian fibrations of operads.
  
 We observe that the category  $\widehat{\Algbrd}_{\Spc}$ is obtained by base change of the cartesian fibration $\widehat{\Algbrd} \rightarrow \widehat{\Cat}$ along the inclusion $\widehat{\Spc} \rightarrow \widehat{\Cat}$. The latter admits a right adjoint, which implies that $h$ admits a right adjoint $h^R$ such that for every object $\Acal$ in $\widehat{\Algbrd}$ the canonical map $hh^R \Acal \rightarrow \Acal$ is cartesian. It follows that the right square in the above diagram is horizontally right adjointable, so we have a commutative diagram of symmetric monoidal categories and lax symmetric monoidal functors as follows:
 \[
 \begin{tikzcd}
 & \widehat{\Algbrd} \arrow{d}{q}  \arrow{r}{h^R} & \widehat{\Algbrd}_{\Spc}  \arrow{d}{q_{\Spc}} \\
 \Alg_{\LM}(\Pr^L) \arrow{r}{ji} \arrow{ur}{r i' p'^R} & \widehat{\Op_{\Assos}} \arrow{r}{\id} &  \widehat{\Op_{\Assos}}
 \end{tikzcd}
 \]
 As before, we observe that $ji$ factors through the lax symmetric monoidal inclusion of $\Alg_{\Assos}(\Pr^L)$ inside $\widehat{\Op_{\Assos}}$, so we have an induced diagram of symmetric monoidal categories and lax symmetric monoidal functors
\[
\begin{tikzcd}
 & \widehat{\Algbrd}|_{\Alg_{\Assos}(\Pr^L)} \arrow{d}{q'}  \arrow{r}{h'^R} & (\widehat{\Algbrd}_{\Spc})|_{\Alg_{\Assos}(\Pr^L)} \arrow{d}{q'_{\Spc}} 
 \\
 \Alg_{\LM}(\Pr^L) \arrow{r}{u} \arrow{ur}{\theta'} & \Alg_{\Assos}(\Pr^L) \arrow{r}{\id} & \Alg_{\Assos}(\Pr^L).
 \end{tikzcd}
\]

We denote by
\[
\theta: \Alg_{\LM}(\Pr^L) \rightarrow \widehat{\Algbrd}_{\Spc}
\]
the lax symmetric monoidal functor obtained by composing $h'^R$ and $\theta'$. It follows from corollary \ref{prop is enriched} that $\theta'$ factors through the full subcategory of $\widehat{\Algbrd}_{\Spc}$ on the enriched categories.

 If $\Mcal_l$ is a presentable symmetric monoidal category, thought of as a commutative algebra object in $\Alg_{\Assos}(\Pr^L)$, we obtain in particular a lax symmetric monoidal functor
\[
\theta_\Mcal: \Mcal \modd(\Pr^L) \rightarrow \widehat{\Cat}^{\Mcal}.
\]
\end{construction}

\begin{proposition}\label{prop theta pres cartesian}
The map $\theta$ from construction \ref{construction theta} is a morphism of cartesian fibrations over $\Alg_{\Assos}(\Pr^L)$.
\end{proposition}
\begin{proof}
Using proposition \ref{prop thetaprime pres cartesian}, we reduce to showing that the morphism
\[
h'^R: \widehat{\Algbrd}_{\Alg_{\Assos}(\Pr^L)} \rightarrow (\widehat{\Algbrd}_{\Spc})|_{\Alg_{\Assos}(\Pr^L)}
\]
is a morphism of cartesian fibrations over $\Alg_{\Assos}(\Pr^L)$. This is a direct consequence of the fact that it is right adjoint to a morphism of cocartesian fibrations.
\end{proof}

\begin{corollary}\label{coro restriction 2}
Let $F: \Mcal \rightarrow \Mcal'$ be a colimit preserving symmetric monoidal functor between presentable symmetric monoidal categories. Then there is a commutative square of symmetric monoidal categories and lax symmetric monoidal functors
\[
\begin{tikzcd}
\Mcal'\modd(\Pr^L) \arrow{r}{\theta_{\Mcal'}} \arrow{d}{F^*} & \widehat{\Cat}^{\Mcal'} \arrow{d}{F^R_!} \\
\Mcal\modd(\Pr^L) \arrow{r}{\theta_{\Mcal}} & \widehat{\Cat}^{\Mcal}
\end{tikzcd}
\]
where $F^*$ denotes the functor of restriction of scalars along $F$.
\end{corollary}
\begin{proof}
This is deduced from proposition \ref{prop theta pres cartesian} using similar arguments as those from the proof of corollary \ref{coro restriction}.
\end{proof}

\begin{proposition}\label{prop theta el caso de spc}
The lax symmetric monoidal functor
\[
\theta_{\Spc} : \Pr^L = \Spc\modd(\Pr^L) \rightarrow \widehat{\Cat}
\]
is equivalent to the  (lax symmetric monoidal) forgetful functor $\Pr^L \rightarrow \widehat{\Cat}$.
\end{proposition}
\begin{proof}
This is a direct consequence of proposition \ref{prop el caso de spc}.
\end{proof}	
	
\begin{notation}
Let $\Mcal$ be a presentable symmetric monoidal category. We denote by  $\overline{\Mcal}$ the image of the unit of $\Mcal\modd(\Pr^L)$ under $\theta_{\Mcal}$. This is a commutative algebra in $\widehat{\Cat}^\Mcal$ whose underlying symmetric monoidal category is equivalent to $\Mcal$.

In the special case $\Mcal = \Cat$ we use the notation $\Catscr = \overline{\Cat}$. This is the symmetric monoidal $2$-category of categories. More generally, for each $n \geq 1$ we set $n\Catscr = \overline{\nCat}$. This is the symmetric monoidal $(n+1)$-category of $n$-categories. We also set $\omega\Catscr = \overline{\omegaCat}$. This is the symmetric monoidal $\omega$-category of $\omega$-categories.
\end{notation}

We finish by considering a variant of the functor $\theta_\Mcal$ from construction \ref{construction theta} which admits a left adjoint.

\begin{notation}
Let $\kappa$ be an uncountable regular cardinal. We denote by $\Pr^L_\kappa$ the subcategory of $\Pr^L$ on the $\kappa$-compactly generated categories and functors which preserve $\kappa$-compact objects.  We equip $\Pr^L_\kappa$ with the restriction of the symmetric monoidal structure from $\Pr^L$.

Let $\Mcal$ be a commutative algebra in $\Pr^L_\kappa$. We denote by $\widehat{\LMod'}|_{\Mcal\modd(\Pr^L_\kappa)}$ the base change of the projection $\widehat{\LMod'} \rightarrow \widehat{\Op}_{\LM}$ along the inclusion $\Mcal\modd(\Pr^L_\kappa) \rightarrow \widehat{\Op}_{\LM}$. We denote by $\LMod'^\kappa|_\Mcal$ the full subcategory of $\widehat{\LMod'}|_{\Mcal\modd(\Pr^L_\kappa)}$ consisting of those triples $(\Acal, \Pcal, \Ccal)$ of an object $\Ccal$ in $\Mcal\modd(\Pr^L_\kappa)$, a small algebroid $\Acal$ in $\Mcal$ with category of objects $X$, and a left module $\Pcal$ such that for every $x$ in $X$ the object $\Pcal(x)$ in $\Ccal$ is $\kappa$-compact.
\end{notation}

\begin{proposition}\label{prop repr y corepr bifbr}
Let $\kappa$ be an uncountable regular cardinal and let $\Mcal$ be a commutative algebra in $\Pr^L_\kappa$. Then the projection 
\[
p^\kappa = (p^\kappa_1, p^\kappa_2): \LMod'^\kappa|_\Mcal \rightarrow \Mcal\modd(\Pr^L_\kappa) \times \Algbrd(\Mcal)
\]
 is a representable bifibration.
\end{proposition}
\begin{proof}
Consider first the projection 
\[
\widehat{\LMod'}|_{\Mcal\modd(\Pr^L_\kappa)} \rightarrow \Mcal \modd(\Pr^L_\kappa) \times \widehat{\Algbrd}(\Mcal).
\]
This arises by base change from the projection of corollary \ref{coro Lmodprime repr bifibr} so we conclude that it is a representable bifibration. It follows directly from the definition that $\LMod'^\kappa|_\Mcal$ is still a cocartesian fibration over $ \Mcal \modd(\Pr^L_\kappa)$, and for every $\Ccal$ in $\Mcal\modd(\Pr^L_\kappa)$ the projection $\LMod'^\kappa|_\Mcal(\Ccal) \rightarrow \Algbrd(\Mcal)$ is a right fibration. Note that this right fibration is represented by $j^! \theta'_\Mcal(\Ccal)$, where $j$ is the inclusion of the full subcategory of $\kappa$-compact objects inside $\Ccal$. We conclude that $p^\kappa$ is a representable bifibration, as desired.
\end{proof}

\begin{notation}
Let $\kappa$ be an uncountable regular cardinal and let $\Mcal$ be a commutative algebra in $\Pr^L_\kappa$. We denote by 
\[
\theta'^\kappa_\Mcal: \Mcal\modd(\Pr^L_\kappa) \rightarrow \Algbrd(\Mcal).
\]
the functor classifying the projection $p^\kappa$ from proposition \ref{prop repr y corepr bifbr}.
\end{notation}

\begin{remark}\label{remark theta prime kappa}
Let $\kappa$ be an uncountable regular cardinal and let $\Mcal$ be a commutative algebra in $\Pr^L_\kappa$. The functor $\theta'^\kappa_\Mcal$ can be obtained as the composition
\[
\Mcal\modd(\Pr^L_\kappa) \xrightarrow{(p^\kappa_1)^R}  \LMod'^\kappa|_\Mcal \xrightarrow{p^\kappa_2} \Algbrd(\Mcal).
\]
Composing $(p^\kappa_1)^R$ with the inclusion of $\LMod'^\kappa|_\Mcal$ inside $\widehat{\LMod'}|_{\Mcal\modd(\Pr^L_\kappa)}$ yields a section of the cocartesian fibration
\[
p'^\kappa_1 : \widehat{\LMod'}|_{\Mcal\modd(\Pr^L_\kappa)} \rightarrow \Mcal \modd(\Pr^L_\kappa) .
\]
By corollary \ref{coro Lmodprime repr bifibr}  the projection $p'^\kappa_1$ admits a right adjoint. It follows that there is a lax commutative triangle
\[
\begin{tikzcd}[row sep = 3em]
 \Mcal \modd(\Pr^L_\kappa) \arrow{r}{(p^\kappa_1)^R} \arrow{dr}[swap, name = diag]{(p'^\kappa_1)^R} & \LMod'^\kappa|_\Mcal \arrow{d}{} \arrow[Rightarrow, to=diag, shorten <>=8pt] \\
& \widehat{\LMod'}|_{\Mcal\modd(\Pr^L_\kappa)}.
\end{tikzcd}
\]
Composing with the projection $\LMod'|_{\Mcal\modd(\Pr^L_\kappa)} \rightarrow \widehat{\Algbrd}(\Mcal)$ we obtain a natural transformation
\[
\theta'^\kappa_\Mcal \rightarrow \theta'_\Mcal|_{\Mcal\modd(\Pr^L_\kappa)}
\]
of functors $\Mcal \modd(\Pr^L_\kappa)  \rightarrow \widehat{\Algbrd}(\Mcal)$. For each object $\Ccal$ in $\Mcal\modd(\Pr^L_\kappa)$, the morphism of algebroids 
\[
\theta'^\kappa_\Mcal(\Ccal) \rightarrow \theta'_\Mcal|_{\Mcal\modd(\Pr^L_\kappa)}(\Ccal)
\]
is cartesian for the projection $\widehat{\Algbrd}(\Mcal) \rightarrow \widehat{\Cat}$, and lies above the inclusion of the full subcategory of $\kappa$-compact objects inside $\Ccal$.
\end{remark}

\begin{notation}\label{notation thetaM kappa}
Let $\kappa$ be an uncountable regular cardinal and let $\Mcal$ be a commutative algebra in $\Pr^L_\kappa$. Consider the composite functor
\[
\Mcal\modd(\Pr^L_\kappa) \xrightarrow{\theta'^\kappa_\Mcal} \Algbrd(\Mcal) \xrightarrow{} \Algbrd(\Mcal)_{\Spc}
\]
where the second map is the colocalization functor. It follows from proposition \ref{prop theta el caso de spc} together with the description of $\theta'^\kappa_\Mcal$ from remark \ref{remark theta prime kappa} that the above composite map factors through $\Cat^\Mcal$. We denote by 
\[
\theta^\kappa_\Mcal : \Mcal\modd(\Pr^L_\kappa) \rightarrow \Cat^\Mcal_{\Spc}
\]
the resulting functor.
\end{notation}

\begin{remark}
Let $\kappa$ be an uncountable regular cardinal and let $\Mcal$ be a commutative algebra in $\Pr^L_\kappa$. It follows from remark \ref{remark theta prime kappa} that there is a natural transformation
\[
\theta^\kappa_\Mcal \rightarrow \theta_\Mcal
\]
such that for every object $\Ccal$ in $\Mcal\modd(\Pr^L_\kappa)$, the functor of enriched categories $\theta^\kappa_\Mcal(\Ccal) \rightarrow \theta_\Mcal(\Ccal)$ exhibits $\theta^\kappa_\Mcal(\Ccal)$ as the full subcategory of $\theta_\Mcal(\Ccal)$ on those objects which correspond to $\kappa$-compact objects in $\Ccal$. 
\end{remark}

\begin{lemma}\label{lemma hay presentability}
Let $\kappa$ be an uncountable regular cardinal and let $\Mcal$ be a commutative algebra in $\Pr^L_\kappa$. Then  the projection 
\[
p_2^\kappa: \LMod'^\kappa|_{\Mcal} \rightarrow \Algbrd(\Mcal)
\]
 preserves colimits.
\end{lemma}
\begin{proof}
It follows from a combination of \cite{HTT} proposition 5.5.7.10 and \cite{HA} lemma 4.8.4.2 and that $\Pr^L_\kappa$ is a presentable symmetric monoidal category. Using \cite{HA} corollary 4.2.3.7 we see that $\Mcal\modd(\Pr^L_\kappa)$ is also presentable.

Using proposition \ref{prop repr y corepr bifbr} we see that the fibers of the cocartesian fibration $p_1^\kappa$ are presentable and for every map $F: \Ccal \rightarrow \Dcal$ in $\Mcal\modd(\Pr^L_\kappa)$ the induced functor
\[
F_! : \LMod'^\kappa_{\Mcal} (\Ccal) \rightarrow \LMod'^\kappa_{\Mcal} (\Dcal)
\]
preserves colimits. Hence $p_1^\kappa$ admits all relative colimits. Note that the projection map $\LMod'^\kappa|_\Mcal \rightarrow \Algbrd(\Mcal)$ sends $p_1^\kappa$-cocartesian arrows to invertible arrows, and for every $\Ccal$ in $\Mcal\modd(\Pr^L_\kappa)$ the functor $\LMod'^\kappa_{\Mcal} (\Ccal) \rightarrow \Algbrd(\Mcal)$ preserves colimits. Applying \cite{HTT} proposition 4.3.1.9 and 4.3.1.10 we have that $p_2^\kappa$ itself preserves colimits.
\end{proof}

\begin{proposition}\label{prop have left adjoint}
Let $\kappa$ be an uncountable regular cardinal and let $\Mcal$ be a commutative algebra in $\Pr^L_\kappa$. Then the functor $\theta^\kappa_\Mcal$ admits a left adjoint.
\end{proposition}
\begin{proof}
Thanks to the adjoint functor theorem, it suffices now to show that $\theta^\kappa_\Mcal$ is accessible and preserves limits. Since the inclusion $\Cat^\Mcal \rightarrow \Algbrd(\Mcal)$ creates limits and sufficiently filtered colimits, it suffices to show that the functor $\theta'^\kappa_\Mcal$ is accessible and limit preserving. The fact that $\theta'^\kappa_\Mcal$ is accessible follows from the description of $\theta'^\kappa_\Mcal$ from remark \ref{remark theta prime kappa}, together with lemma \ref{lemma hay presentability}.

It remains to prove that $\theta'^\kappa_\Mcal$ is limit preserving. Recall that the projection $\Algbrd(\Mcal) \rightarrow \Cat$ admits all relative limits. We claim that the composite map
\[
\Mcal\modd(\Pr^L_\kappa) \xrightarrow{\theta'^\kappa_\Mcal} \Algbrd(\Mcal) \rightarrow \Cat
\]
is limit preserving. Examining the commutative diagram from remark \ref{remark lmod fits into} shows that the composition of the above map with the inclusion of $\Cat$ into $\widehat{\Cat}$ admits a factorization as follows:
\[
\Mcal\modd(\Pr^L_\kappa) \xrightarrow{\mu} \Funct([1], \widehat{\Cat}) \xrightarrow{\ev_0} \widehat{\Cat}
\]
The map $\mu$ is such that the composition
\[
\Mcal\modd(\Pr^L_\kappa) \xrightarrow{\mu} \Funct([1], \widehat{\Cat}) \xrightarrow{\ev_1} \widehat{\Cat}
\]
recovers the canonical projection obtained by composing the following series of forgetful functors:
\[
\Mcal\modd(\Pr^L_\kappa) \rightarrow \Pr^L_\kappa \rightarrow \widehat{\Cat}
\]
The description of $\theta'^\kappa_\Mcal$ from remark \ref{remark theta prime kappa} shows that $\mu \ev_0$ is the composite functor
\[
\Mcal\modd(\Pr^L_\kappa) \rightarrow \Pr^L_\kappa = \Cat^{\rex(\kappa)} \rightarrow \widehat{\Cat}
\]
where the middle equivalence is given by passage to $\kappa$-compact objects (\cite{HTT} proposition 5.5.7.8), and the last arrow is the usual forgetful functor. This composition is indeed limit preserving, so it follows that the composition of $\theta'^\kappa_\Mcal$ with the forgetful functor to $\Cat$ is limit preserving, as we claimed.

Consider now a limit diagram $X^\lhd : \Ical^\lhd \rightarrow \Mcal\modd(\Pr^L_\kappa)$. Denote by $\ast$ the initial object of $\Ical$. Let $Y$ be the category of objects for the $\Mcal$-algebroid $\theta'^\kappa_\Mcal(X^\lhd(\ast))$. In other words, $Y$ is the full subcategory of the category underlying the $\Mcal$-module $X^\lhd(\ast)$ on the $\kappa$-compact objects. 

Note that the composite map
\[
\Ical^\lhd \xrightarrow{X^\lhd} \Mcal\modd(\Pr^L_\kappa) \xrightarrow {\theta'^\kappa_\Mcal} \Algbrd(\Mcal) \rightarrow \Cat
\]
factors through $\Cat_{Y/}$. Therefore we have that $\theta'^\kappa_\Mcal X^\lhd$ factors through  $\Algbrd(\Mcal) \times_{\Cat} \Cat_{Y/}$. Since the projection $\Algbrd(\Mcal) \times_{\Cat} \Cat_{Y/} \rightarrow \Cat_{Y/}$ is a cartesian fibration and $Y$ is initial there, we have that the fiber of $\Algbrd(\Mcal)$ over $Y$ is a colocalization of $\Algbrd(\Mcal) \times_{\Cat} \Cat_{Y/}$. From this we may construct a diagram
\[
Z^\lhd : \Ical^\lhd \rightarrow \Algbrd_Y(\Mcal)
\]
equipped with a natural transformation
\[
\begin{tikzcd}[row sep = 3em]
 \Ical^\lhd \arrow{r}{Z^\lhd} \arrow{dr}[swap, name = diag]{\theta'^\kappa_\Mcal X^\lhd} & \Algbrd_Y(\Mcal) \arrow{d}{} \arrow[Rightarrow, to=diag, shorten <>=8pt] \\
& \Algbrd(\Mcal)
\end{tikzcd}
\]
such that the induced map $Z^\lhd(\ast) \rightarrow X^\lhd(\ast)$ is an isomorphism, and for every object $i$ in $\Ical$ the morphism of algebroids $Z^\lhd(i) \rightarrow X^\lhd(i)$ is cartesian over $\Cat$. Using \cite{HTT} propositions 4.3.1.9 and 4.3.1.10 we see that in order to show that $\theta'^\kappa_\Mcal X^\lhd$ is a limit diagram it suffices to show that that $Z^\lhd$ is a limit diagram.

Observe that $\LMod'^\kappa|_\Mcal$ also has the structure of a cartesian fibration over $\Cat$, and the projection $\LMod'^\kappa|_\Mcal \rightarrow \Algbrd(\Mcal)$ is a morphism of cartesian fibrations over $\Cat$. Repeating the above procedure, we may write $Z^\lhd$ as $p_2^\kappa W^\lhd$ where 
\[
W^\lhd: \Ical^\lhd \rightarrow (\LMod'^\kappa|_\Mcal)_Y
\]
is a diagram which comes equipped with a natural transformation
\[
\begin{tikzcd}[row sep = 3em]
 \Ical^\lhd \arrow{r}{W^\lhd} \arrow{dr}[swap, name = diag]{(p_1^\kappa)^R X^\lhd} & (\LMod'^\kappa|_\Mcal)_Y \arrow{d}{} \arrow[Rightarrow, to=diag, shorten <>=8pt] \\
& \LMod'^\kappa|_\Mcal
\end{tikzcd}
\]
such that the induced map $W^\lhd(\ast) \rightarrow (p_1^\kappa)^R X^\lhd(\ast)$ is an isomorphism, and for every object $i$ in $\Ical$ the morphism $W^\lhd(i) \rightarrow (p_1^\kappa)^R X^\lhd(i)$ is cartesian over $\Cat$.

We now observe that a diagram in $\Algbrd_Y(\Mcal)$ is a limit diagram if and only if its images in $\Algbrd_{\lbrace a, b \rbrace}(\Mcal)$ are limit diagrams for every map $\lbrace a, b \rbrace \rightarrow Y$, where $\lbrace a, b \rbrace$ denotes a two-element set. Fix one such map and let $U^\lhd$ be the composition of $W^\lhd$ with the induced morphism
\[
(\LMod'^\kappa|_\Mcal)_Y \rightarrow (\LMod'^\kappa|_\Mcal)_{\lbrace a, b \rbrace}.
\]
Our task is to show that the composite map
\[
\Ical^\lhd \xrightarrow{U^\lhd} (\LMod'^\kappa|_\Mcal)_{\lbrace a, b \rbrace} \rightarrow \Algbrd_{\lbrace a, b \rbrace}(\Mcal)
\]
is a limit diagram. Recall now from \cite{Hinich} that $\Algbrd_{\lbrace a, b \rbrace}(\Mcal)$ is the category of algebras in a certain presentable monoidal category $\Mcal_{\lbrace a, b \rbrace}$. The category $(\LMod'^\kappa|_\Mcal)_{\lbrace a, b \rbrace}$ fits into a pullback square
\[
\begin{tikzcd}
(\LMod'^\kappa|_\Mcal)_{\lbrace a, b \rbrace} \arrow{d}{} \arrow{r}{} & (\LMod'^\kappa|_{\Mcal_{\lbrace a, b \rbrace}})_{[0]} \arrow{d}{} \\
\Mcal\modd(\Pr^L_\kappa) \arrow{r}{} & \Mcal_{\lbrace a, b \rbrace}\modd(\Pr^L_\kappa)
\end{tikzcd}
\]
where the bottom horizontal arrow maps a $\Mcal$-module $\Ccal$ to the $\Mcal_{\lbrace a, b \rbrace}$-module $\Funct(\lbrace a, b \rbrace, \Ccal)$. We note that the right vertical arrow admits a factorization through the category 
\[
 \Mcal_{\lbrace a, b \rbrace}\modd(\Pr^L_\kappa)_{\Mcal_{\lbrace a, b \rbrace} /}
 \]
  of pointed $\Mcal_{\lbrace a , b \rbrace}$-modules. The resulting projection has a fully faithful section
  \[
   S: \Mcal_{\lbrace a, b \rbrace}\modd(\Pr^L_\kappa)_{\Mcal_{\lbrace a, b \rbrace} /} \rightarrow (\LMod'^\kappa|_{\Mcal_{\lbrace a, b \rbrace}})_{[0]} 
  \]
  whose image consists of those triples $(A, \Dcal, M)$ of an $\Mcal_{\lbrace a , b \rbrace}$-module $\Dcal$, an algebra $A$ in $\Mcal_{\lbrace a , b \rbrace}$, and a $\kappa$-compact $A$-module $M$ in $\Dcal$ for which the structure map $A \otimes M \rightarrow M$ exhibits $A$ as the endomorphism object of $M$. The description of endomorphism objects from \cite{Hinich} proposition 6.3.1 shows that in fact the composite map 
\[
\Ical^\lhd \xrightarrow{U^\lhd}  (\LMod'^\kappa|_\Mcal)_{\lbrace a, b \rbrace} \rightarrow (\LMod'^\kappa|_{\Mcal_{\lbrace a, b \rbrace}})_{[0]}
\]
factors through the image of $S$. Furthermore, the resulting diagram $\Ical^\lhd \rightarrow   \Mcal_{\lbrace a, b \rbrace}\modd(\Pr^L_\kappa)$ is the image of $X^\lhd$ under the functor
\[
\Mcal\modd(\Pr^L_\kappa) \rightarrow \Mcal_{\lbrace a, b \rbrace}\modd(\Pr^L_\kappa)
\] 
and is therefore a limit diagram. The fact that $\theta'^\kappa_\Mcal X^\lhd$ is a limit diagram now follows from the fact that the composite map
\[
\Mcal_{\lbrace a, b \rbrace}\modd(\Pr^L_\kappa)_{\Mcal_{\lbrace a, b \rbrace} /} \xrightarrow{S} (\LMod'^\kappa|_{\Mcal_{\lbrace a, b \rbrace}})_{[0]} \rightarrow \Alg(\Mcal_{\lbrace a, b \rbrace})
\]
is the functor that sends a pointed $\Mcal$-module to the endomorphism object of the basepoint, which admits a left adjoint.
\end{proof}

\begin{remark}
Let $\kappa$ be an uncountable regular cardinal and let $\Mcal$ be a commutative algebra in $\Pr^L_\kappa$. Passing to $\kappa$-compact objects induces an equivalence between the symmetric monoidal category $\Pr^L_\kappa$ and the symmetric monoidal category $\Cat^{\rex(\kappa)}$ of small categories admitting $\kappa$-small colimits, and functors which preserve those colimits. In particular, we have an equivalence between $\Mcal\modd(\Pr^L_\kappa)$ and $\Mcal^{\kappa\text{-comp}}\modd(\Cat^{\rex(\kappa)})$, where   $\Mcal^{\kappa\text{-comp}}$ denotes the full subcategory of $\Mcal$ on the $\kappa$-compact objects. 

From this point of view, the left adjoint to $\theta^\kappa_\Mcal$ maps small $\Mcal$-enriched categories into $\kappa$-cocomplete categories tensored over $\Mcal^{\kappa\text{-comp}}$.  We think about this as a version of the functor of free $\kappa$-cocompletion in the context of enriched category theory.
\end{remark}

\subsection{Bimodules}\label{subsection bimod}
 We now discuss the notion of bimodule between algebroids.

\begin{notation}
For each $\BM$-operad $\Mcal$ we denote by $\Mcal_l$ and $\Mcal_r$ the $\Assos^-$ and $\Assos^+$ components of $\Mcal$. We denote by $\Mcal_m$ the fiber of $\Mcal$ over the bimodule object in $\BM$.
\end{notation}

\begin{definition}
Let $\Mcal$ be a $\BM$-operad.  Let $\Acal, \Bcal$ be algebroids in $\Mcal_l$ and $ \Mcal_r$ respectively, with categories of objects $X$ and $Y$ respectively. An $\Acal-\Bcal$-bimodule   in $\Mcal$ is a $\BM_{X, Y}$-algebra in $\Mcal$, whose underlying $\Assos_X$ and $\Assos_Y$ algebras are identified with $\Acal$ and $\Bcal$.
\end{definition}

\begin{remark}
Let $\Mcal$ be a $\BM$-operad. Let $\Acal, \Bcal$ be algebroids in $\Mcal_l$ and $ \Mcal_r$ respectively, with categories of objects $X$ and $Y$ respectively. A bimodule $\Pcal$ between them assigns to each pair of objects $(x, y)$ in $X\times Y^\op$ an object $\Pcal(x, y)$ in $\Mcal_m$. For every $n \geq 0, m \geq 0$, and every sequence of arrows 
\[
x'_0 = x_0 \leftarrow x'_1, x_1, \ldots, x_{n-1} \leftarrow x'_n , x_n \leftarrow x'_{n+1} = x_{n+1}
\]
in $X$ and 
\[
y'_0 = y_0 \leftarrow y'_1,  y_2, \ldots, y_{m-1} \leftarrow y'_m , y_m \leftarrow y'_{m+1} = y_{m+1}
\]
in $Y$, the bimodule $\Pcal$ induces a multimorphism
\[
\lbrace \Acal( x'_1, x_1 ),\ldots,  \Acal(x'_n, x_n), P(x'_{n+1}, y_0), \Bcal(y'_1, y_1) , \ldots , \Bcal(y'_m, y_m) \rbrace \rightarrow \Pcal(x'_0, y_{m+1}) 
\]
in $\Mcal$. In the case when $\Mcal$ is a $\BM$-monoidal category (in other words, $\Mcal_l$ and $\Mcal_r$ are monoidal categories and $\Mcal_m$ is a bimodule between them), this induces  a morphism
\[
 \Acal( x'_1, x_1 ) \otimes \ldots \otimes  \Acal(x'_n, x_n) \otimes  P(x'_{n+1}, y_0) \otimes \Bcal(y'_1, y_1) \otimes  \ldots \otimes \Bcal(y'_m, y_m) \rightarrow \Pcal(x'_0, y_{m+1}) .
\]
In particular, in the cases when $n = 1, m = 0$ or $n = 0, m = 1$ and all arrows are identities we obtain, for each pair of objects $x_1, x_2$ in $X$ and each pair of objects $y_1, y_2$ in $Y$, a map
\[
\Acal(x_1, x_2) \otimes \Pcal(x_2, y_1) \rightarrow \Pcal(x_1, y_1)
\]
(the left action of $\Acal$ on $\Pcal$) and 
\[
\Pcal(x_2, y_1) \otimes \Bcal(y_1, y_2) \rightarrow  \Pcal(x_2, y_2).
\]
(the right adjoint of $\Bcal$ on $\Pcal$). These actions commute with each other, and are compatible with the units and composition of $\Acal$ and $\Bcal$, up to homotopy.
\end{remark}

\begin{example}
Let $\Mcal$ be an associative operad. Then every algebroid $\Acal: \Assos_X \rightarrow \Mcal$ defines, by precomposition with the projection $\BM_{X, X} \rightarrow \Assos_X$ arising from the equivalence of example \ref{remark diagonal BMX to assos}, an $\Acal-\Acal$-bimodule $\Pcal$ in $\Mcal$. This has the following properties:
\begin{itemize}
\item For every pair of objects $x, x'$ in $X$ we have $\Pcal(x', x) = \Acal(x', x)$.
\item For every triple of objects $x, x', x''$ in $X$, the action map
\[
\lbrace \Acal(x'', x') , \Pcal(x', x) \rbrace  \rightarrow \Pcal(x'', x)
\]
is equivalent, under the identifications of the previous item, to the composition map
\[
\lbrace \Acal(x'', x') , \Acal(x', x) \rbrace \rightarrow \Acal(x'', x).
\]
\item For every triple of objects $x, x', x''$ in $X$, the action map
\[
\lbrace \Pcal(x', x) , \Acal(x, x'') \rbrace \rightarrow \Pcal(x', x'')
\]
is equivalent, under the identifications of the first item, to the composition map
\[
\lbrace \Acal(x', x) , \Acal(x, x'') \rbrace \rightarrow \Acal(x', x'').
\]
\end{itemize}
We call $\Pcal$ the diagonal bimodule of $\Acal$. We will usually use the notation $\Acal(-, -)$ for $\Pcal$, and in the case when $\Acal$ is an $\Mcal$-enriched category, we instead write $\Hom_{\Acal}(-, -)$.
\end{example}

\begin{remark}
Let $\Mcal$ be an associative operad. As discussed in \cite{Hinich}, an $\Mcal$-algebroid $\Acal$ with category of objects $X$ defines an associative algebra $\widetilde{\Acal}$ in a certain associative operad $\Mcal_X$, and an $\Acal-\Acal$-bimodule is the same data as an $\widetilde{\Acal}-\widetilde{\Acal}$-bimodule in $\Mcal_X$. Under this dictionary, the diagonal bimodule of $\Acal$ corresponds to the diagonal bimodule of $\widetilde{\Acal}$. Via the folding equivalence of \cite{Hinich} section 3.6, the diagonal bimodule of $\widetilde{\Acal}$ defines a left $ \widetilde{\Acal} \boxtimes \widetilde{\Acal}^\op$-module in $\Mcal_{X \times X^\op}$. If $\Mcal$ is a presentable symmetric monoidal category, this defines a morphism 
\[
\Acal \otimes \Acal^\op \rightarrow \overline{\Mcal}\]
 with the property that each pair of objects $(y, x)$ gets mapped to $\Acal(y, x)$. We call this the Hom functor of $\Acal$. This determines a morphism of algebroids
 \[
 \Acal \rightarrow \Funct(\Acal^\op, \overline{\Mcal})
 \]
 which is the Yoneda embedding. It was show in \cite{Hinich} corollary 6.2.7 that this map is fully faithful. Note that thanks to proposition \ref{prop funct cuando uno es especial} the algebroid $\Funct(\Acal^\op, \overline{\Mcal})$ is in fact an enriched category. We conclude that any $\Mcal$-algebroid admits a fully faithful embedding into an $\Mcal$-enriched category.
\end{remark}

\begin{construction}
Consider the functor $\Alg_{\BM_{-,-}}(-)$ defined by the composition
\[
\Cat^\op \times \Cat^\op \times \Op_{\BM} \xrightarrow{\BM_{-,-} \times \id_{\Op_{\BM}}} \Op_{\BM}^\op \times \Op_{\BM} \xrightarrow{\Alg_-(-)} \Cat.
\]
For each $\BM$-algebroid $\Mcal$ we denote by $\BMod(\Mcal)$ the total category of the cartesian fibration associated to the functor $\Alg_{\BM_{-,-}}(\Mcal): \Cat^\op \times \Cat^\op \rightarrow \Cat$. This comes equipped with a forgetful functor $\BMod(\Mcal) \rightarrow \Algbrd(\Mcal_l) \times \Algbrd(\Mcal_r)$. For each pair of algebroids $\Acal$ in $\Mcal_l$ and $\Bcal$ in $\Mcal_r$ we denote by $_{\Acal}\kr \BMod_{\Bcal}(\Mcal)$ the fiber  over $(\Acal, \Bcal)$, and call it the category of $\Acal-\Bcal$-bimodules.

The assignment $\Mcal \mapsto \BMod(\Mcal)$ defines a functor $\BMod(-): \Op_{\BM} \rightarrow \widehat{\Cat}$. Let $\BMod$ be the total category of the associated cocartesian fibration. Note that this fits into a commutative square
\[
\begin{tikzcd}
\BMod \arrow{r}{} \arrow{d}{} & \Algbrd \times \Algbrd \arrow{d}{} \\
\Cat \times \Cat \times \Op_{\BM} \arrow{r}{} & \Cat \times \Cat \times \Op_{\Assos} \times \Op_{\Assos}.
\end{tikzcd}
\]
Here the left vertical arrow and right vertical arrow are the two-sided fibrations classified by the functors $\Alg_{\BM_{-, -}}(-)$ and $\Algbrd_-(-)$, and the horizontal arrows are the functors of ``forgetting the bimodule". 
\end{construction}

\begin{proposition}
Let $\Mcal$ be a presentable $\BM$-monoidal category. Then the projection $\BMod(\Mcal) \rightarrow \Algbrd(\Mcal_l) \times \Algbrd(\Mcal_r)$ is a cartesian fibration. Moreover, a morphism $F: (\Acal, \Pcal, \Bcal) \rightarrow (\Acal', \Pcal', \Bcal')$ is cartesian if and only if for every pair of objects $x$ in $\Acal$ and $y$ in $\Bcal$, the induced map $\Pcal(x, y) \rightarrow \Pcal'(F(x), F(y))$ is an equivalence.
\end{proposition}
\begin{proof}
This follows from a variation of the argument of proposition \ref{prop restrict scalars lmod}.
\end{proof}



\tableofcontents

\section{Enriched adjunctions and conical limits}

Let $\Mcal$ be a presentable symmetric monoidal category. Let  $F: \Ccal \rightarrow \Dcal$ be a functor of $\Mcal$-enriched categories.  Given a functor $G: \Dcal \rightarrow \Ccal$, a natural transformation $\epsilon: FG \rightarrow \id_\Dcal$ is said to exhibit $G$ as right adjoint to $F$ if for every pair of objects $c$ in $\Ccal$ and $d$ in $\Dcal$ the induced map
\[
\Hom_\Ccal(c, G(d)) \xrightarrow{\epsilon(d) \circ F_*(-) } \Hom_\Dcal(F(c), d) 
\]
is an isomorphism.

In \ref{subsection local} we study a local version of the notion of right adjoint to $F$, where the object $G(d)$ may only be well defined for specific values of $d$. We prove a criterion guaranteeing the existence of local adjoints to $F$ in the case when $F$ is obtained as the limit of a family of functors $F_i : D(i) \rightarrow D'(i)$. Given another $\Mcal$-enriched category $\Jcal$ and an object $d$ in $\Funct(\Jcal, \Dcal)$, we show that the right adjoint to $F_* : \Funct(\Jcal, \Ccal) \rightarrow \Funct(\Jcal, \Dcal)$ at a functor $d$ exists provided that the right adjoint to $F$ at $\ev_j d$ exists for all $j$ in $\Jcal$.

In \ref{subsection global} we study the notion of adjoint functors between $\Mcal$-enriched categories, as a special case of the notion of adjunction in a $2$-category. We show that $F$ admits a right adjoint if and only if it admits local right adjoints at every object in $\Dcal$.

In \ref{subsection conical} we specialize to the case when $F$ is the diagonal map $\Dcal \rightarrow \Funct(\Ical_\Mcal, \Dcal)$, where $\Ical_\Mcal$ is the  $\Mcal$-enriched category obtained from a category $\Ical$ by pushforward along the unit map $\Spc \rightarrow \Mcal$. A local right adjoint at an object $X$ in $\Funct(\Ical_\Mcal, \Dcal)$ is called a conical limit of $X$. We show that the data of a conical limit defines in particular an extension of $X$ to a diagram $\Ical^\lhd_\Mcal \rightarrow \Dcal$. We study the interactions of the notion of conical limits with changes in the enriched category - in particular, we are able to conclude that a conical limit in $\Dcal$ defines in particular a  limit diagram in the category underlying $\Dcal$. Specializing our discussion of local adjoints we obtain basic results on the existence of conical limits on limits of enriched categories, and in enriched categories of functors.

In \ref{subsection modules} we study how the notions of enriched adjunctions and conical limits interact with the procedure of enrichment of presentable modules over $\Mcal$. We show that if $F: \Ccal \rightarrow \Dcal$ is a morphism of presentable modules then the induced functor of $\Mcal$-enriched categories admits a right adjoint, and it also admits a left adjoint provided that $F$ admits a left adjoint which strictly commutes with the action of $\Mcal$. As a consequence, we are able to conclude that if $\Dcal$ is a presentable module of $\Mcal$, then the induced $\Mcal$-enriched category admits all small conical limits and colimits. In particular, the canonical enrichment of $\Mcal$ over itself is conically complete and cocomplete.

\subsection{Local adjoints}\label{subsection local} We begin  with a general discussion of the notion of  locally defined adjoints for functors of enriched categories.

\begin{definition}
Let $\Mcal$ be a presentable monoidal category. Let $F: \Dcal \rightarrow \Dcal'$ be a functor of $\Mcal$-categories, and let $d'$ be an object of $\Dcal'$. Let $d$ be an object of $\Dcal$ equipped with a morphism $\epsilon: F(d) \rightarrow d'$ in the category underlying $\Dcal'$. We say that the pair $(d, \epsilon)$ is right adjoint to $F$ at $d'$ if for every object $e$ in $\Dcal$ the composite functor
\[
\Hom_\Dcal(e, d) \xrightarrow{F_*} \Hom_{\Dcal'}(F(e), F(d)) \xrightarrow{\epsilon} \Hom_{\Dcal'}(F(e), d')
\]
is an equivalence.
\end{definition}

\begin{remark}\label{local right adj como cart map}
Let $G: \Mcal' \rightarrow \Mcal$ be a colimit preserving monoidal functor between presentable monoidal categories. Let $F: \Dcal \rightarrow \Dcal'$ be a functor of $\Mcal$-categories, and let $(d, \epsilon)$ be right adjoint to $F$ at an object $d'$ in $\Dcal'$. Then $(d,\epsilon)$ is also right adjoint to $(G^R)_! F : (G^R)_! \Dcal \rightarrow (G^R)_!\Dcal'$ at $d'$.

In the particular case when $\Mcal' = \Spc$ and $G$ is the unit map for $\Mcal$ in $\Pr^L$, the functor $G^R$ is the functor $(\tau_\Mcal)_!$ that sends each $\Mcal$-enriched category to its underlying category. We thus see that $(d, \epsilon)$ is right adjoint to $(\tau_\Mcal)_! F$ at  $d'$.

Let $p: \Ecal \rightarrow [1]$ be the cocartesian fibration associated to the functor $(\tau_\Mcal)_! F$. The pair $(d, \epsilon)$ induces a morphism $\alpha$ between $(0, d)$ and $(1, d')$ in $\Ecal$. The condition that $(d,\epsilon)$ be right adjoint to $(\tau_\Mcal)_!F$ at  $d'$ is equivalent to the condition that $\alpha$ be a $p$-cartesian morphism. In particular, we see that the pair $(d, \epsilon)$ right adjoint to $F$ at $d'$ is unique if it exists.
\end{remark}

\begin{definition}\label{def local horiz right adj}
Let $\Mcal$ be a presentable monoidal category and consider
a commutative square of $\Mcal$-enriched categories
\[
\begin{tikzcd}
\Dcal \arrow{r}{F}  \arrow{d}{T} & \Dcal' \arrow{d}{T'} \\
\Ecal \arrow{r}{G} & \Ecal'
\end{tikzcd}.
\]
Let $d'$ be an object of $\Dcal'$. We say that the above square is horizontally right adjointable at $d'$ if there is a pair $(d, \epsilon)$ right adjoint to $F$ at $d'$, and moreover the induced map $T'\epsilon: GTd \rightarrow T'd'$ is right adjoint to $G$ at $T'd'$.
\end{definition}

\begin{remark}\label{remark send cartesian to cartesian}
Let $\Mcal = \Spc$ and consider a commutative diagram as in definition \ref{def local horiz right adj}. This induces a morphism $\mathcal{T}$ of cocartesian fibrations over $[1]$ between the fibrations $\Ecal_F$, $\Ecal_G$ associated to $F$ and $G$. The square is horizontally right adjointable at $d'$ if there is a nontrivial cartesian arrow in $\Ecal_F$ with target $d'$, whose image under $\mathcal{T}$ is cartesian.
\end{remark}

\begin{proposition}\label{prop passage to limit}
Let $\Mcal$ be a presentable  monoidal category and let $\Ical$ be a category. Let $D, D': \Ical \rightarrow \Cat^\Mcal$ be functors, and let $\eta: D \rightarrow D'$ be a natural transformation. Denote by $F: \Dcal \rightarrow \Dcal'$ the limit of $\eta$. For each $i$ in $\Ical$ let $p_i: \Dcal \rightarrow D(i)$ and $p_i': \Dcal' \rightarrow D'(i)$ be the projections. Assume that for each arrow $\alpha: i \rightarrow j$ in $\Ical$, the square
\[
\begin{tikzcd}
D(i) \arrow{r}{\eta_i}  \arrow{d}{D(\alpha)} & D'(i) \arrow{d}{D'(\alpha)} \\
D(j) \arrow{r}{\eta_j} & D'(j)
\end{tikzcd}
\]
is horizontally right adjointable at $p_i'(d')$.  Then 
\begin{enumerate}[\normalfont (i)]
\item There exists a right adjoint to $F$ at $d'$.
\item A pair $(d, \epsilon)$ is right adjoint to $F$ at $d'$ if and only if $(p_i d, p_i'(\epsilon))$ is right adjoint to $\eta_i$ at $p_i'(d')$ for  every $i$ in $\Ical$.
\end{enumerate}
\end{proposition}
\begin{proof}
Consider first the case $\Mcal = \Spc$. Passing to cocartesian fibrations of the functors $\eta_i$ and $F$, we obtain a diagram 
\[
\Ecal_\eta: \Ical \rightarrow (\Cat^{\cocart}_{/[1]})_{\lbrace 1 \rbrace/}
\] 
with limit $\Ecal_F$, where the right hand side denotes the undercategory of the category of cocartesian fibrations and morphisms of cocartesian fibrations over $[1]$, under the cocartesian fibration $\lbrace 1 \rbrace \rightarrow [1]$.

 The adjointability of the square in the statement implies that the composition of $\Ecal_\eta$ with the forgetful functor $(\Cat^{\cocart}_{/[1]})_{\lbrace 1 \rbrace/} \rightarrow (\Cat_{/[1]})_{\lbrace 1 \rbrace/}$ factors through the subcategory $(\Cat_{/[1]})_{[1]/^{\cart}}$  of categories over $[1]$ equipped with a cartesian section, and functors which preserve this section. The case $\Mcal = \Spc$ of the proposition now follows from the fact that the projections 
  \[
 (\Cat^\cocart_{/[1]})_{\lbrace 1 \rbrace/} \rightarrow (\Cat_{/[1]})_{\lbrace 1 \rbrace/} \leftarrow (\Cat_{/[1]})_{[1]/^{\cart}} \rightarrow (\Cat_{/[1]})_{[1]/}
 \]
 create limits.

 We now consider the general case. By virtue of the above and  the uniqueness claim from remark \ref{local right adj como cart map}, it suffices to show that if $(d, \epsilon)$ is such that $(p_j d, p'_j  \epsilon)$ is right adjoint to $\eta_i$ at $p'_i d'$ for all $i$ then it is right adjoint to $F$ at $d'$. Let $e$ be an object in $\Dcal$. Note that there is a functor $R: \Ical \rightarrow \Funct([2], \Mcal)$ whose value on each index $i$ is given by
\[
\Hom_{D(i)}(p_ie, p_id) \xrightarrow{(\eta_i)_*} \Hom_{D'(i)}(\eta_i(e), \eta_i(d) ) \xrightarrow{p'_i\epsilon} \Hom_{D'(i)}(\eta_i(e), p'_i d')
\]
and which has a limit given by
\[
\Hom_\Dcal(e, d) \xrightarrow{F_*} \Hom_{\Dcal'}(F(e), F(d)) \xrightarrow{\epsilon} \Hom_{\Dcal'}(F(e), d').
\]
For each $i$ in $\Ical$ the composition of the maps in $R(i)$ is an isomorphism, since $(p_id, p'_i\epsilon)$ is right adjoint to $\eta_i$ at $p'_id'$. We conclude that the composition of the maps in $\lim_\Ical R$ is an isomorphism, which means that $(d, \epsilon)$ is right adjoint to $F$ at $d'$, as desired.
\end{proof}

\begin{proposition}\label{prop right adj in functor cats}
Let $\Mcal$ be a presentable symmetric monoidal category. Let $\Jcal$ be an $\Mcal$-enriched category  $F: \Dcal \rightarrow \Dcal'$ be a functor of $\Mcal$-enriched categories. Let $d'$ be an object in $\Funct(\Jcal,\Dcal')$ and assume that for all objects $j$ in $\Jcal$ there exists a right adjoint to $F$ at $\ev_j d'$. Then
\begin{enumerate}[\normalfont (i)]
\item There exists a right adjoint to $F_*: \Funct(\Jcal, \Dcal) \rightarrow \Funct(\Jcal, \Dcal')$ at $d'$.
\item A pair $(d, \epsilon)$ is right adjoint to $F_*$ at $d'$ if and only if for every $j$ in $\Jcal$ the pair $(\ev_j d, \ev_j \epsilon)$ is right adjoint to $F$ at  $\ev_j d'$.
\end{enumerate}
\end{proposition}
\begin{proof}
Let $\Scal$ be the full subcategory of $\Cat^\Mcal$ on those $\Mcal$-enriched categories $\Jcal$ for which the proposition holds. We claim that $\Scal$ is closed under colimits in $\Cat^\Mcal$. Let $J: \Ical \rightarrow \Scal$ be a diagram, and let $\Jcal$ be its colimit in $\Cat^\Mcal$. We then have that the functor 
\[
F_*: \Funct(\Jcal, \Dcal) \rightarrow \Funct(\Jcal, \Dcal')
\]
 is obtained by passage to the limit of the functors 
 \[
 \eta_i: \Funct(J(i), \Dcal) \rightarrow \Funct(J(i), \Dcal')
 \] given by composition with $F$. Let $d'$ be an object in $\Funct(\Jcal, \Dcal')$ and assume that there exists a right adjoint to $F$ at $\ev_j d'$ for every $j$ in $\Jcal$.  The fact that $J(i)$ belongs to $\Scal$ for all $i$ implies that for every arrow $\alpha: i \rightarrow i'$ in $\Ical$ the square
\[
\begin{tikzcd}
\Funct(J(i'), \Dcal) \arrow{r}{\eta_{i'}} \arrow{d}{J(\alpha)^*} & \Funct(J(i'),\Dcal') \arrow{d}{J(\alpha)^*} \\ 
\Funct(J(i),\Dcal) \arrow{r}{\eta_{i}} & \Funct(J(i'),\Dcal')
\end{tikzcd}
\]
is horizontally right adjointable at $d'|_{J(i')}$. It then follows from proposition \ref{prop passage to limit} that there is indeed a right adjoint to $F_*$ at $d'$, and moreover a pair $(d, \epsilon)$ is right adjoint to $F_*$ at $d'$ if and only if the associated pair $(d|_{J(i)}, \epsilon|_{J(i)})$ is right adjoint to $\eta_i$ at $d'|_{J(i)}$ for all $i$ in $\Ical$. Using again the fact that $J(i)$ belongs to $\Scal$ for every $i$, we see that this happens if and only if $(\ev_j d, \ev_j \epsilon)$ is right adjoint to $F$ at $\ev_j d'$ for every $j$ in $\Jcal$ which is in the image of the map $J(i) \rightarrow \Jcal$ for some $i$ in $\Ical$. Since $\Jcal$ is the colimit of the objects $J(i)$, we have that the maps $J(i) \rightarrow \Jcal$ are jointly surjective, so we have that $\Jcal$ belongs to $\Scal$.

Since $\Cat^\Mcal$ is generated under colimits by cells, our result will follow if we show that for every $m$ in $\Mcal$ the enriched category underlying the cell $C_m$ belongs to $\Scal$. Let $d': C_m \rightarrow \Dcal'$ be an $m$-cell in $\Dcal'$, with source and target objects $d'_0$ and $d'_1$. Let $(d_0, \epsilon_0)$ and $(d_1, \epsilon_1)$ be right adjoint to $F$ at $d'_0$ and $d'_1$ respectively.

We claim that there is a unique enhancement of this data to a pair $(d, \epsilon)$ of an $m$-cell $d: C_m \rightarrow \Dcal$ and a morphism $\epsilon: F_*d \rightarrow d'$. Indeed, the data of a cell $d$ between $d_0$ and $d_1$ corresponds to a map $m \rightarrow \Hom_\Dcal(d_0, d_1)$. A map $\epsilon$ lifting $\epsilon_0$ and $\epsilon_1$ consists of an identification of $m$-cells $\epsilon_1 F_* d = d' \epsilon_0$. The fact that $(d_1, \epsilon_1)$ is right adjoint to $F$ at $d'_1$ implies that there is a unique such pair $(d, \epsilon)$, as claimed.

It remains to show that $(d, \epsilon)$ is right adjoint to $F_*$ at $d'$. Let $e: C_m \rightarrow \Dcal$ be another $m$-cell with source and target objects $e_0$ and $e_1$.   Recall from corollary \ref{coro homs en funct} that we have a cartesian square
\[
\begin{tikzcd}
\Hom_{\Funct(C_m, \Dcal)}( e, d) \arrow{r}{\ev_0} \arrow{d}{\ev_1} & \Hom_{\Dcal}(e_0, d_0) \arrow{d}{} \\
\Hom_\Dcal(e_1, d_1) \arrow{r}{} & \shom_\Mcal(m, \Hom_\Dcal(e_0, d_1))
\end{tikzcd}
\]
where the right and bottom arrows are given by composition with the cells $d$ and $e$, respectively. The functor $F$ induces a map from the above square to the cartesian square
\[
\begin{tikzcd}
\Hom_{\Funct(C_m, \Dcal')}(F_*e, F_*d) \arrow{r}{\ev_0} \arrow{d}{\ev_1} & \Hom_{\Dcal'}(Fe_0, Fd_0) \arrow{d}{} \\
\Hom_{\Dcal'}(Fe_1, Fd_1) \arrow{r}{} & \shom_\Mcal(m, \Hom_{\Dcal'}(Fe_0, Fd_1))
\end{tikzcd}
\]
where the right and bottom arrows are given by composition with $F_*d$ and $F_*e$. Finally, composition with $\epsilon$ yields a map from the above to the cartesian square
\[
\begin{tikzcd}
\Hom_{\Funct(C_m, \Dcal')}(F_*e, d') \arrow{r}{\ev_0} \arrow{d}{\ev_1} & \Hom_{\Dcal'}(Fe_0, d'_0) \arrow{d}{} \\
\Hom_{\Dcal'}(Fe_1, d'_1) \arrow{r}{} & \shom_\Mcal(m, \Hom_{\Dcal'}(Fe_0, d'_1))
\end{tikzcd}
\]
where the right and bottom arrows are given by composition with $d'$ and $F_*e$. We thus see that there is a commutative cube
\[
    \begin{tikzcd}[row sep=1.5em, column sep = 0.05em]
        & \Hom_{\Funct(C_m, \Dcal')}(F_*e, d') \arrow[dd, "\ev_1", pos=0.75] \arrow[rr, "\ev_0"] &&
    \Hom_{\Dcal'}(Fe_0, d'_0) \arrow[dd] \\
    \Hom_{\Funct(C_m, \Dcal)}( e, d) \arrow[rr,"\ev_0", pos=0.7] \arrow[ur, "\epsilon F_*"] \arrow[dd, "\ev_1"] &&
    \Hom_{\Dcal}(e_0, d_0) \arrow[dd] \arrow[ur,"\epsilon_0 F"] \\
 & \Hom_{\Dcal'}(Fe_1, d'_1) \arrow[rr] && \shom_\Mcal(m, \Hom_{\Dcal'}(Fe_0, d'_1)) \\
    \Hom_\Dcal(e_1, d_1) \arrow[rr] \arrow[ur, "\epsilon_1 F"] &&  \shom_\Mcal(m, \Hom_\Dcal(e_0, d_1)) \arrow[ur, "\epsilon_1 F"]
    \end{tikzcd}
\]
with cartesian front and back faces. Since $(d_0, \epsilon_0)$ and $(d_1, \epsilon_1)$ are right adjoint to $F$ at $d'_0$ and $d'_1$, the bottom left, bottom right, and top right diagonal arrows are isomorphisms. We conclude that the top left diagonal arrow is an isomorphism, which means that $(d, \epsilon)$ is right adjoint to $F_*$ at $d$, as desired.
\end{proof}

\begin{corollary}\label{coro check isos upon evaluation}
Let $\Mcal$ be a presentable symmetric  monoidal category. Let $\Jcal$ and $\Dcal$ be $\Mcal$-enriched categories. Then a morphism $\epsilon: d \rightarrow d'$ in $\Funct(\Jcal, \Dcal)$ is an isomorphism if and only if $\ev_j \epsilon$ is an isomorphism for every $j$ in $\Jcal$.
\end{corollary}
\begin{proof}
Specialize proposition \ref{prop right adj in functor cats} to the case when $F$ is the identity of $\Dcal$.
\end{proof}

\begin{corollary}\label{coro functorialidad right adjoints}
Let $\Mcal$ be a presentable symmetric monoidal category. Let $F: \Dcal \rightarrow \Dcal'$ be a functor of $\Mcal$-enriched categories, and assume that $F$ admits a right adjoint at $d'$ for every $d'$ in $\Dcal'$. Then there is a unique functor $F^R: \Dcal' \rightarrow \Dcal$ equipped with a morphism $\epsilon: F F^R \rightarrow \id_{\Dcal'}$ such that for every $d'$ in $\Dcal'$ the pair $(F^R(d'), \epsilon(d'))$ is right adjoint to $F$ at $\Dcal$.
\end{corollary}
\begin{proof}
Specialize proposition \ref{prop right adj in functor cats} to the case when $\Jcal = \Dcal'$ and take $(F^R,\epsilon)$ to be right adjoint to $F_*$ at $\id_{\Dcal'}$.
\end{proof}

\subsection{Global adjoints} \label{subsection global} We now discuss the notion of adjunction between functors of enriched categories. We will obtain this as a particular case of the general notion of adjunction in a $2$-category.

\begin{definition}\label{definition adjunction}
Let $\Dcal$ be a $2$-category. An arrow $\alpha: d \rightarrow e$ in $\Dcal$ is said to admit a right adjoint if there is an arrow $\alpha^R: e \rightarrow d$ and a pair of $2$-cells $\eta: \id_d \rightarrow \alpha^R \alpha$ and $\epsilon: \alpha \alpha^R \rightarrow \id_e$ satisfying the following two conditions:
\begin{itemize}
\item The composite $2$-cell
\[
\alpha = \alpha  \id_d \xrightarrow{\id_\alpha  \eta} \alpha \alpha^R \alpha \xrightarrow{\epsilon \id_\alpha} \id_e \alpha = \alpha
\]
is equivalent to the identity.
\item The composite $2$-cell
\[
\alpha^R = \id_d \alpha^R \xrightarrow{\eta \id_{\alpha^R}} \alpha^R \alpha \alpha^R \xrightarrow{\id_{\alpha^R} \epsilon} \alpha^R \id_e = \alpha^R
\]
is equivalent to the identity.
\end{itemize}
In this situation, we say that $\alpha^R$ is right adjoint to $\alpha$, and we call $\eta$ and $\epsilon$ the unit and counit of the adjunction, respectively. We say that $\alpha$ admits a left adjoint if it admits a right adjoint as a morphism in $\Dcal^{2\dsh\op}$.
\end{definition}

\begin{example}\label{example adjunctions cat}
In the case $\Dcal = \Catscr$, definition \ref{definition adjunction} recovers the usual notion of adjunction between functors of categories.
\end{example}

We refer to \cite{RV} for a proof of the following fundamental theorems.

\begin{theorem}
There exists a $2$-category $\Adj$ equipped with an epimorphism $L: [1] \rightarrow \Adj$ such that for every $2$-category $\Dcal$ composition with $L$ induces an equivalence between the space of functors $\Adj \rightarrow \Dcal$ and the space of functors $[1] \rightarrow \Dcal$ which pick out a right adjointable arrow in $\Dcal$.
\end{theorem}
\begin{theorem}
Let $\Dcal$ be a $2$-category. Then the following spaces are equivalent:
\begin{enumerate}[\normalfont (i)]
\item The space of arrows in $\Dcal$ which admit a right adjoint.
\item The space of triples $(\alpha, \alpha^R, \eta)$ of an arrow $\alpha: d \rightarrow e$ in $\Dcal$, an arrow $\alpha^R: e \rightarrow d$ in $\Dcal$, and a $2$-cell $\eta: \id_d \rightarrow \alpha^R \alpha$ which can be extended to an adjunction between $\alpha$ and $\alpha^R$.
\end{enumerate}
The equivalence is given by mapping a triple $(\alpha, \alpha^R, \eta)$ to $\alpha$.
\end{theorem}

We now specialize to the case when $\Dcal$ is the $2$-category of categories enriched in a presentable symmetric monoidal category

\begin{notation}
Let $\Mcal$ be a presentable symmetric monoidal category. We denote by $\Catscr^\Mcal$ the image of $\Cat^\Mcal$ under the composite functor
\[
\Cat^\Mcal\modd(\Pr^L) \xrightarrow{\theta_{\Cat^\Mcal}}\widehat{\Cat}^{\Cat^\Mcal} \xrightarrow{(\tau_\Mcal)_!} \widehat{\Cat}^{\Cat} \hookrightarrow \reallywidehat{2\kr\Cat}.
\]
We call $\Catscr^\Mcal$ the $2$-category of $\Mcal$-enriched categories. Note that the $1$-category underlying $\Catscr^\Mcal$ is equivalent to $\Cat^\Mcal$.
\end{notation}

\begin{example}
Let $n \geq 0$. Then $\Catscr^{n\kr\Cat}$ is the $2$-category underlying $(n+1)\Catscr$.
\end{example}

\begin{remark}\label{remark natural transf}
Let $\Mcal$ be a presentable symmetric monoidal category and let $\Ccal, \Dcal$ be $\Mcal$-enriched categories. Then the Hom category $\Hom_{\Catscr^\Mcal}(\Ccal, \Dcal)$ is the category underlying $\Funct(\Ccal, \Dcal)$.

Let $F, G: \Ccal \rightarrow \Dcal$ be functors, and $\eta: F \rightarrow G$ be a $2$-cell. Let $x, y$ be a pair of objects of $\Ccal$. Then there is an induced morphism
\[
C_{\Hom_\Ccal(x, y)} \otimes C_{1_\Mcal} \rightarrow \Dcal.
\]
Examining the description of the product of cells from remark \ref{remark product cells} we obtain a commutative diagram
\[
\begin{tikzcd}
\Hom_\Ccal(x, y)  \arrow{r}{F_*} \arrow{d}{G_*} & \Hom_{\Dcal}(Fx, Fy) \arrow{d}{\eta(y) \circ -} \\ \Hom_\Dcal(Gx, Gy) \arrow{r}{- \circ \eta(x)} & \Hom_\Dcal(Fx, Gy).
\end{tikzcd}
\]
\end{remark}

\begin{definition}
Let $\Mcal$ be a presentable symmetric monoidal category and let $F: \Ccal \rightarrow \Dcal$ be a functor of $\Mcal$-enriched categories. We say that $F$ admits a right adjoint if it admits a right adjoint when thought of as a morphism in the $2$-category $\Catscr^\Mcal$. 
\end{definition}

Our next goal is to show that a functor of enriched categories admits a right adjoint if and only if it admits local right adjoints at every point.

\begin{lemma}\label{lemma criterio adj}
Let $\Jcal$ and $\Dcal$ be $2$-categories. Let $\eta : F \rightarrow G$ be a morphism in $\Funct(\Jcal, \Dcal)$. Then $\eta$ admits a right adjoint if and only if for every morphism $\alpha: j \rightarrow j'$ in $\Jcal$, the commutative square
\[
\begin{tikzcd}
F(j) \arrow{d}{F \alpha} \arrow{r}{\eta_j} &  G(j) \arrow{d}{G\alpha} \\
F(j') \arrow{r}{\eta_{j'}} & G(j')
\end{tikzcd}
\]
is horizontally right adjointable.
\end{lemma}
\begin{proof}
Combine \cite{Haugsadj} theorem 4.6 and corollary 3.15.
\end{proof}

\begin{proposition}\label{prop global adj iff local adjs}
Let $\Mcal$ be a presentable symmetric monoidal category and let $F: \Ccal \rightarrow \Dcal$ be a functor of $\Mcal$-enriched categories. Then $F$ admits a right adjoint if and only if for every object $d$ in $\Dcal$ there exists a right adjoint for $F$ at $d$.
\end{proposition}
\begin{proof}
Assume first that $F$ admits a right adjoint $F^R$. Denote by $\eta: \id_{\Ccal} \rightarrow F^R F$ and $\epsilon: F F^R \rightarrow \id_{\Dcal}$ the unit and counit of the adjunction, respectively. Let $d$ be an object in $\Dcal$. We claim that the morphism
\[
\epsilon(d) : FF^R(d) \rightarrow \id_\Dcal(d) = d
\]
exhibits $F^R(d)$ as right adjoint to $F$ at $d$. Let $c$ be an object of $\Ccal$. We have to show that the map $V$ given by the composition
\[
\Hom_{\Ccal}(c, F^R(d)) \xrightarrow{F_*} \Hom_{\Dcal}(F(c), FF^R(d)) \xrightarrow{\epsilon(d)\circ -} \Hom_{\Dcal}(F(c), d)
\]
is an isomorphism. Observe that there is a map $W$ going in the opposite direction, given by the following composition:
\[
\Hom_{\Dcal}(F(c), d) \xrightarrow{(F^R)_*} \Hom_{\Ccal}(F^RF(c), F^R(d)) \xrightarrow{- \circ \eta(c)} \Hom_\Ccal(c, F^R(d))
\]
We claim that $V$ and $W$ are inverse equivalences. Observe that the map $WV$ is given by the composition
\[
\Hom_\Ccal(c, F^R(d)) \xrightarrow{(F^R F)_*} \Hom_\Ccal(F^RF(c), F^R F F^R(d) ) \xrightarrow{F^R\epsilon(d) \circ - \circ \eta(c) } \Hom_\Ccal(c, F^R(d)).
\]
Thanks to remark \ref{remark natural transf}, we can rewrite the above as the composition
\[
\Hom_\Ccal(c, F^R(d)) \xrightarrow{\eta(F^R(d)) \circ - } \Hom_\Ccal(c, F^RFF^R(d)) \xrightarrow{F^R \epsilon(d) \circ - } \Hom_\Ccal(c, F^R(d)).
\]
This is equivalent to the identity thanks to the second condition in definition \ref{definition adjunction}. The fact that $VW$ is equivalent to the identity follows from similar arguments. We conclude that $V$  is indeed an isomorphism, so $\epsilon(d)$ exhibits $F^R(d)$ as right adjoint to $F$ at $d$.

Assume now that for every object $d$ in $\Dcal$ there exists a right adjoint for $F$ at $d$.   To show that $F$ admits a right adjoint it suffices to show that its image under the Yoneda embedding 
\[
\Catscr^\Mcal \rightarrow \Funct((\Catscr^\Mcal)^{1\dsh\op}, \Catscr)
\]
admits a right adjoint.   Applying lemma \ref{lemma criterio adj}, we reduce to showing that for every functor $\alpha: \Jcal \rightarrow \Jcal'$ between $\Mcal$-enriched categories, the commutative square of categories
\[
\begin{tikzcd}
(\tau_\Mcal)_!\Funct(\Jcal', \Ccal) \arrow{d}{\alpha^*} \arrow{r}{F_*} &  (\tau_\Mcal)_!\Funct(\Jcal', \Dcal) \arrow{d}{\alpha^*} \\
(\tau_\Mcal)_!\Funct(\Jcal, \Ccal) \arrow{r}{F_*} & (\tau_\Mcal)_!\Funct(\Jcal, \Dcal)
\end{tikzcd}
\]
is horizontally right adjointable. Using proposition \ref{prop right adj in functor cats} together with remark \ref{local right adj como cart map} we see that the above commutative square of categories is horizontally right adjointable at every object of $ (\tau_\Mcal)_!\Funct(\Jcal', \Dcal)$. Our claim now follows from the fact that local and global adjointability agree for functors between categories.
\end{proof}

\subsection{Conical limits} \label{subsection conical}

We now specialize the notion of local adjoints to obtain a theory of conical limits and colimits.

\begin{notation}
Let $\Mcal$ be a presentable monoidal category, and $\Ical$ be a category. We denote by $\Ical_\Mcal$ the image of $\Ical$ under the functor $\Cat \rightarrow \Cat^\Mcal$ induced by pushforward along the unit map $\Spc \rightarrow \Mcal$. We note that for every $\Mcal$-enriched category $\Dcal$, there is a correspondence between functors $\Ical_\Mcal \rightarrow \Dcal$ and functors $\Ical \rightarrow (\tau_\Mcal)_!\Dcal$.
\end{notation}

\begin{definition}
Let $\Mcal$ be a presentable symmetric monoidal category. Let $\Ical$ be  category and let $\Dcal$ be an $\Mcal$-enriched category. Let $X: \Ical_\Mcal \rightarrow \Dcal$ be a functor. We say that $X$ admits a conical limit if the functor
\[
\Delta: \Dcal = \Funct(1_\Mcal, \Dcal) \rightarrow \Funct(\Ical_\Mcal, \Dcal)
\]
 of precomposition with the projection $\Ical \rightarrow [0]$ admits a right adjoint at $X$. In this case, we call its right adjoint at $X$ the conical limit of $X$. We say that $X$ admits a conical colimit if the induced diagram $X^\op: \Ical^\op \rightarrow \Dcal^\op$ admits a conical limit - in this case we define the conical colimit of $X$ to be the conical limit of $X^\op$. 
\end{definition}

\begin{remark}\label{cones vs arrows}
Let $\Ical$ be a category and denote by $\Ical^\lhd$ the category obtained from $\Ical$ by adjoining a final object. We have a pushout diagram in $\Cat$
\[
\begin{tikzcd}[column sep = large]
\Ical \arrow{r}{\id_{\Ical} \times \lbrace 0 \rbrace} \arrow{d}{} & \Ical \times [1] \arrow{d}{} \\ {[0]} \arrow{r}{} & \Ical^\lhd.
\end{tikzcd}
\]

Let $\Mcal$ be a presentable symmetric monoidal category. It follows from the above that for every $\Mcal$-enriched category $\Dcal$ there is a pullback diagram of spaces
\[
\begin{tikzcd}
\Hom_{\Cat^\Mcal}(\Ical_\Mcal^\lhd, \Dcal) \arrow{r}{} \arrow{d}{} & \Hom_{\Cat^\Mcal}(C_{1_\Mcal}, \Funct(\Ical_\Mcal, \Dcal)) \arrow{d}{} \\ \Hom_{\Cat^\Mcal}(1_\Mcal, \Dcal) \arrow{r}{} & \Hom_{\Cat^\Mcal}(\Ical_\Mcal, \Dcal)
\end{tikzcd}
\]
where the right vertical arrow is given by evaluation at the source, and to bottom horizontal arrow is the diagonal map. Hence we see that a pair $(d, \epsilon)$ of an object $d$ in $\Dcal$ and a morphism $\epsilon: \Delta d \rightarrow X$ in $\Funct(\Ical, \Dcal)$ is the same data as a diagram $X^\lhd : \Ical_\Mcal^\lhd \rightarrow \Dcal$. In particular, we have that a conical limit for a diagram $X: \Ical_\Mcal \rightarrow \Dcal$ can be identified with a particular kind of extension of $X$ to a diagram $X^\lhd: \Ical_\Mcal^\lhd \rightarrow \Dcal$.
\end{remark}

\begin{remark}\label{remark pushforward cones}
Let $\Mcal$ be a presentable symmetric monoidal category and let $F: \Dcal \rightarrow \Dcal'$ be a functor of $\Mcal$-enriched categories. Let $\Ical$ be a category. Then we have a commutative square
\[
\begin{tikzcd}
\Dcal \arrow{r}{\Delta} \arrow{d}{F} & \Funct(\Ical_\Mcal, \Dcal) \arrow{d}{F_*} \\
\Dcal' \arrow{r}{\Delta} & \Funct(\Ical_\Mcal, \Dcal').
\end{tikzcd}
\]
Let $d$ be an object in $\Dcal$ and $\epsilon: \Delta d \rightarrow X$ be a morphism in $\Funct(\Ical_\Mcal, \Dcal)$, associated to a diagram $X^\lhd: \Ical_\Mcal^\lhd \rightarrow \Dcal$ under the equivalence of remark \ref{cones vs arrows}. Then the induced pair $(Fd, F_* \epsilon)$ is associated to $F_* X^\lhd$.
\end{remark}

\begin{remark}\label{remark compatibility cones y adj}
Let $G: \Mcal \rightarrow \Mcal'$ be a colimit preserving monoidal functor between presentable symmetric monoidal categories and let $\Dcal$ be an $\Mcal'$ enriched category. It follows from proposition \ref{prop funct y changes} that we have an equivalence
\[
\Funct(-_{\Mcal}, (G^R)_! \Dcal) = (G^R)_!\Funct(-_{\Mcal'}, \Dcal)
\]
of functors  $\Cat \rightarrow \Cat^\Mcal$. Let $\Ical$ be a category. Evaluating the above equivalence at the projection
$\Ical \rightarrow [0]$ we obtain a commutative square
\[
\begin{tikzcd}[column sep = huge]
(G^R)_! \Dcal \arrow{r}{(G^R)_!\Delta} \arrow{d}{\id} & (G^R)_!\Funct(\Ical_{\Mcal'}, \Dcal) \arrow{d}{=} \\
(G^R)_! \Dcal \arrow{r}{\Delta} & \Funct(\Ical_{\Mcal}, (G^R)_!\Dcal).
\end{tikzcd}
\]
Let $d$ be an object in $\Dcal$, and $\epsilon: \Delta d \rightarrow X$ be a morphism in $\Funct(\Ical_\Mcal, \Dcal)$. Using the above square, we can also think about $\epsilon$ as a morphism in $\Funct(\Ical_{\Mcal}, (G^R)_!\Dcal)$. If the pair $(d, \epsilon)$ corresponds to a diagram $X^\lhd: \Ical^\lhd_{\Mcal'} \rightarrow \Dcal$ under the equivalence of remark \ref{cones vs arrows}, then the diagram $\Ical^\lhd_{\Mcal} \rightarrow (G^R)_!\Dcal$ obtained from $(d, \epsilon)$ via the above commutative square is equivalent to the image of $X^\lhd$ under the canonical equivalence
\[
\Hom_{\Cat^{\Mcal'}}(\Ical^\lhd_{\Mcal'}, \Dcal) =  \Hom_{\Cat^{\Mcal}}(\Ical_{\Mcal}^\lhd, (G^R)_! \Dcal).
\]
\end{remark}

We now explore the behavior of conical limits under changes in the  enriching category.
\begin{proposition}\label{remark colimits vs conical}
Let $G: \Mcal \rightarrow \Mcal'$ be a colimit preserving symmetric monoidal functor between presentable symmetric monoidal categories. Let $\Dcal$ be an $\Mcal$-enriched category and let $X'^\lhd: \Ical^\lhd_{\Mcal'} \rightarrow \Dcal$ be a conical limit diagram. Then the induced functor $X^\lhd: \Ical^\lhd_{\Mcal} \rightarrow (G^R)_! \Dcal$ is a conical limit diagram.
\end{proposition}
\begin{proof}
This is a direct consequence of the discussion in remarks \ref{local right adj como cart map} and \ref{remark compatibility cones y adj}.
\end{proof}
\begin{corollary}\label{coro conical is ordinary colim}
Let $\Mcal$ be a presentable symmetric monoidal category. Let $\Dcal$ be an $\Mcal$-enriched category and $X^\lhd: \Ical^\lhd_\Mcal \rightarrow \Dcal$ be a conical limit diagram in $\Dcal$. Then the induced diagram $\Ical^\lhd \rightarrow (\tau_\Mcal)_!\Dcal$ in the category underlying $\Dcal$, is a limit diagram. 
\end{corollary}
\begin{proof}
Apply proposition \ref{remark colimits vs conical} to the unit map $\Spc \rightarrow \Mcal$.
\end{proof}

\begin{proposition}
Let $i: \Mcal \rightarrow \Mcal'$ be a colimit preserving symmetric monoidal functor between presentable symmetric monoidal categories. Assume that $i$ is fully faithful and admits a strictly symmetric monoidal left adjoint. Let $\Dcal$ be an $\Mcal$-enriched category and let  $X^\lhd: \Ical_\Mcal^\lhd \rightarrow \Dcal$ be a conical limit diagram. Then $i_!X^\lhd : \Ical_{\Mcal'}^\lhd \rightarrow i_! \Dcal$ is a conical limit diagram
\end{proposition}
\begin{proof}
It follows from corollary \ref{coro compara internal functs} that we have an equivalence
\[
i_! \Funct(-_\Mcal, \Dcal) = \Funct(-_{\Mcal'} , i_!\Dcal)
\]
of functors $\Cat^\op \rightarrow \Cat^{\Mcal'}$.  Applying it to the projection $\Ical \rightarrow [0]$ we obtain a commutative square
\[
\begin{tikzcd}
i_! \Dcal \arrow{r}{i_! \Delta} \arrow{d}{\id} & i_! \Funct(\Ical_\Mcal, \Dcal) \arrow{d}{=} \\
i_! \Dcal \arrow{r}{\Delta} & \Funct(\Ical_{\Mcal'}, i_!\Dcal).
\end{tikzcd}
\]
Let $(d, \epsilon)$ be the right adjoint to $\Delta: \Dcal \rightarrow \Funct(\Ical_\Mcal, \Dcal)$   at $X = X^\lhd|_{\Ical_\Mcal}$. By virtue of remark \ref{local right adj como cart map}, we have that $(d, \epsilon)$ is also right adjoint to $i_! \Delta: i_! \Dcal \rightarrow i_! \Funct(\Ical_\Mcal, \Dcal)$.  Its image under the above equivalence is right adjoint to $\Delta: i_! \Dcal \rightarrow \Funct(\Ical_{\Mcal'}, \Dcal)$ at  $i_! X$ - in other words, it is a conical limit for $X$. Our result now follows from the fact that the associated limit diagram $\Ical^\lhd_{\Mcal'} \rightarrow i_!\Dcal$ is given by $i_! X^\lhd$.
\end{proof}

\begin{corollary}
Let $m \geq n \geq 1$ and let $\Dcal$ be an $n$-category, thought of as a category enriched in $(n-1)$-categories. Then a diagram $X^\lhd: \Ical^\lhd \rightarrow \Dcal$ is a conical limit diagram if and only if its image under the inclusion functor $i^{n, m}: \nCat \rightarrow m\kr\Cat$ is a conical limit diagram.
\end{corollary}

\begin{definition}
Let $\Mcal$ be a presentable symmetric monoidal category and let $F: \Dcal \rightarrow \Dcal'$ be a functor of $\Mcal$-enriched categories. We say that a conical limit diagram $X^\lhd: \Ical^\lhd \rightarrow \Dcal$ is preserved by $F$ if $FX^\lhd$ is a conical limit diagram in $\Dcal'$. Similarly, a conical colimit diagram $Y^\rhd: \Ical^\rhd \rightarrow \Dcal$ is said to be preserved by $F$ if $FY^\rhd$ is a conical colimit diagram in $\Dcal'$.
\end{definition}

\begin{remark}\label{remark pres colim}
Let $\Mcal$ be a presentable symmetric monoidal category and let $F: \Dcal \rightarrow \Dcal'$ be a functor of $\Mcal$-enriched categories. Let $X^\lhd: \Ical_\Mcal^\lhd \rightarrow \Dcal$ be a conical limit diagram. Then it follows from remark \ref{remark pushforward cones} that $X^\lhd$ is preserved by $F$ if and only if the commutative square
\[
\begin{tikzcd}
\Dcal \arrow{r}{\Delta} \arrow{d}{F} & \Funct(\Ical_\Mcal, \Dcal) \arrow{d}{F_*} \\
\Dcal' \arrow{r}{\Delta} & \Funct(\Ical_\Mcal, \Dcal')
\end{tikzcd}
\]
is horizontally right adjointable at $X$.
\end{remark}

We now specialize propositions \ref{prop passage to limit} and \ref{prop right adj in functor cats} to the case of conical limits.

\begin{proposition}
Let $\Mcal$ be a presentable symmetric monoidal category and let $\Kcal$ be a category. Let $D: \Kcal \rightarrow \Cat^\Mcal$ be a functor, and denote by $\Dcal$ its limit. For each $j$ in $\Kcal$ denote by $p_j: \Dcal \rightarrow D(j)$ the projection. Let $X : \Ical_\Mcal \rightarrow \Dcal$ be a diagram in $\Dcal$. Assume that for every $j$ in $\Kcal$ the diagram $p_j X : \Ical \rightarrow D(j)$ admits a conical limit, which is preserved by the functor $D(\alpha): D(j) \rightarrow D(j')$ for every arrow $\alpha: j \rightarrow j'$ in $\Kcal$. Then
\begin{enumerate}[\normalfont (i)]
\item The diagram $X$ admits a conical limit.
\item An extension $X^\lhd : \Ical^\lhd_\Mcal \rightarrow \Dcal$ is a conical limit diagram if and only if $p_j X^\lhd$ is a conical limit diagram in $D(j)$ for every $j$ in $\Kcal$.
\end{enumerate}
\end{proposition}
\begin{proof}
Combine proposition \ref{prop passage to limit} together with remarks  \ref{remark pushforward cones} and \ref{remark pres colim}.
\end{proof}

\begin{proposition}\label{prop limits in functor cat}
Let $\Mcal$ be a presentable symmetric monoidal category. Let $\Ical$ be a category, and let $\Jcal$ and $\Dcal$ be $\Mcal$-enriched categories. Let $X: \Ical_\Mcal \rightarrow \Funct(\Jcal, \Dcal)$ be a diagram, and assume that for every object $j$ in $\Jcal$, the diagram $\ev_j X : \Ical_\Mcal \rightarrow \Dcal$ admits a conical limit. Then 
\begin{enumerate}[\normalfont (i)]
\item There exists a conical limit for $X$.
\item An extension $X^\lhd: \Ical_\Mcal^\lhd \rightarrow \Funct(\Jcal, \Dcal)$ is a conical limit for $X$ if and only if for every object $j$ in $\Jcal$ the diagram $\ev_j X^{\lhd}$ is a conical limit.
\end{enumerate}
\end{proposition}
\begin{proof}
Apply proposition \ref{prop right adj in functor cats} to the diagonal map $\Delta: \Dcal \rightarrow \Funct(\Ical, \Dcal)$.
\end{proof}

\begin{definition}
Let $\Mcal$ be a presentable symmetric monoidal category and let $\Dcal$ be an $\Mcal$-enriched category.  Let $\Ical$ be a category. We say that $\Dcal$ admits all conical (co)limits of shape $\Ical$ if every diagram $X: \Ical \rightarrow \Dcal$ admits a conical (co)limit. We say that $\Dcal$ is conically (co)complete if it admits all conical (co)limits of shape $\Ical$ for every small category $\Ical$.
\end{definition}

\begin{corollary}
Let $\Mcal$ be a presentable symmetric monoidal category and let $\Dcal$ be an  $\Mcal$-enriched category which admits all conical limits of shape $\Ical$. Then for every $\Mcal$-enriched category $\Jcal$, the $\Mcal$-enriched category $\Funct(\Jcal, \Dcal)$ admits all conical limits of shape $\Ical$.
\end{corollary}
\begin{proof}
Follows directly from proposition \ref{prop limits in functor cat}.
\end{proof}

\begin{corollary}
Let $\Mcal$ be a presentable symmetric monoidal category and let $\Dcal$ be an $\Mcal$-enriched category which admits all conical limits of shape $\Ical$. Then there is a functor $\Funct(\Ical_\Mcal, \Dcal) \rightarrow \Dcal$ which maps each diagram $X: \Ical \rightarrow \Dcal$ to the value of its conical limit at the cone point of $\Ical^\lhd$.
\end{corollary}
\begin{proof}
Apply corollary \ref{coro functorialidad right adjoints} to the diagonal map $\Delta: \Dcal \rightarrow \Funct(\Ical, \Dcal)$.
\end{proof}

For later purposes we record the following basic consequence of proposition \ref{prop limits in functor cat}.

\begin{proposition}\label{prop colimits y epis}
Let $\Mcal$ be a presentable symmetric monoidal category. Let $\Dcal$ be an enriched category and let $f: \Jcal \rightarrow \Jcal'$ be an epimorphism of $\Mcal$-enriched categories which is surjective on objects. Let $X: \Ical_\Mcal \rightarrow \Funct(\Jcal, \Dcal)$ be a diagram admitting a conical colimit $\mathscr{X}$ which is preserved by the evaluation functors $\ev_j : \Funct(\Jcal, \Dcal) \rightarrow \Dcal$ for every object $j$ in $\Jcal$. Assume that for every arrow $\alpha: i \rightarrow i'$ in $\Ical$ the morphism $X(\alpha) : X(i) \rightarrow X(j)$ in  $\Funct( \Jcal, \Dcal)$ belongs to $\Funct(\Jcal', \Dcal)$. Then 
\begin{enumerate}[\normalfont (i)]
\item The functor $\mathscr{X}: \Jcal \rightarrow \Dcal$ factors through $\Jcal'$.
\item For every object $i$ in $\Ical$ the morphism $X(i) \rightarrow \mathscr{X}$ belongs to $\Funct(\Jcal', \Dcal)$.
\end{enumerate}
\end{proposition}
\begin{proof}
Note that since the map $f$ is an epimorphism, the map 
\[
f^*: \Funct(\Jcal', \Dcal) \rightarrow \Funct(\Jcal, \Dcal)
\] is indeed a monomorphism. Our assumptions imply that $X$ factors through the image of $f^*$. Let $X' : \Ical_\Mcal \rightarrow \Funct(\Jcal', \Dcal)$ be the induced diagram. Since $f$ is surjective on objects, we have that for every $j$ in $\Jcal'$ the diagram $\ev_j X'$ admits a conical colimit. It follows from (the dual version of) proposition \ref{prop limits in functor cat} that $X'$ can be extended to a conical colimit diagram $X'^\rhd : \Ical^\rhd_\Mcal \rightarrow \Funct(\Jcal', \Dcal)$. The diagram $f^*X'^\rhd$  is an extension of $X$ whose image under all evaluation functors is a conical colimit diagram. Applying proposition \ref{prop limits in functor cat} again we conclude that $f^*X'^\rhd$ is a conical limit diagram in $\Dcal$. Item (i) now follows from the fact that $\mathscr{X}$ is  equivalent to the value of  $f^*X'^\rhd$ at the cone point $\ast$ in $\Ical^\rhd$. Item (ii) is a consequence of the fact that the morphism $X(i) \rightarrow \mathscr{X}$ is equivalent to the image under $f^*$ of the morphism $X'^\rhd(i) \rightarrow X'^\rhd(\ast)$.
\end{proof}

\subsection{The case of presentable modules}\label{subsection modules}
We now study the interactions between adjunctions with the procedure of enrichment of modules over presentable symmetric monoidal categories.

\begin{proposition}\label{prop adjuntos enriched}
Let $\Mcal$ be a presentable symmetric monoidal category and let $F: \Ccal \rightarrow \Dcal$ be a morphism in $\Mcal\modd(\Pr^L)$. Then
\begin{enumerate}[\normalfont (i)]
\item The functor of enriched categories $\theta_\Mcal(F) : \theta_\Mcal(\Ccal) \rightarrow \theta_\Mcal(\Dcal)$ admits a right adjoint.
\item Assume that $F$ admits a left adjoint $F^L$, and that the canonical structure of oplax morphism of $\Mcal$-modules on $F^L$ is strict. Then $\theta_\Mcal(F)$ admits a left adjoint.
\end{enumerate}
\end{proposition}
\begin{proof}
We first prove item (i). Let $F^R: \Dcal \rightarrow \Ccal$ be the right adjoint to the functor underlying $F$. Let $d$ be an object in $\Dcal$ and let $\epsilon(d) : FF^R(d) \rightarrow d$ be the counit of the adjunction at $d$. We claim that $(F^R(d), \epsilon(d))$ is right adjoint to $\theta_\Mcal(F)$ at $d$. To see this, we have to show that for every $c$ in $\Ccal$ the composite map
\[
\Hom_{\theta_\Mcal(\Ccal)}(c, F^R(d) ) \xrightarrow{\theta_\Mcal(F)_*} \Hom_{\theta_\Mcal(\Dcal)}(F(c), FF^R(d) )  \xrightarrow{\epsilon(d)} \Hom_{\theta_\Mcal(\Dcal)}(F(c), d)
\]
is an isomorphism. It suffices to show that for every $m$ in $\Mcal$ the image of the above composition under the functor $\Hom_\Mcal(m, -)$ is an isomorphism. This is equivalent to
\[
\Hom_{\Ccal}(m\otimes c, F^R(d)) \xrightarrow{F_*} \Hom_\Dcal(m \otimes F(c), FF^R(d) ) \xrightarrow{\epsilon(d)} \Hom_\Dcal(m\otimes F(c), d)
\]
which is indeed an isomorphism, since $(F^R(d), \epsilon(d))$ is right adjoint to $F$ at $d$. Item (i) now follows from proposition \ref{prop global adj iff local adjs}.

We now prove item (ii). Let $d$ be an object in $\Dcal$ and let $\eta(d): d \rightarrow FF^L(d)$ be the counit of the adjunction at $d$. We claim that $(F^L(d), \eta(d))$ is left adjoint to $\theta_\Mcal(F)$ at $d$. To see this, we have to show that for every $c$ in $\Ccal$ the composite map
\[
\Hom_{\theta_\Mcal(\Ccal)} (F^L(d), c) \xrightarrow{\theta_\Mcal(F)_*}\Hom_{\theta_\Mcal(\Dcal)} ( F F^L(d), F(c) ) \xrightarrow{\eta(d)} \Hom_{\theta_\Mcal(\Dcal)} ( d, F(c) )
\]
is an isomorphism.  It suffices to show that for every $m$ in $\Mcal$ the image of the above composition under the functor $\Hom_\Mcal(m, -)$ is an isomorphism. This is equivalent to the composite map
\[
\Hom_{\Ccal} (m \otimes F^L(d), c) \xrightarrow{F_*}\Hom_{\Dcal} ( m \otimes F F^L(d), F(c) ) \xrightarrow{\id_m \otimes \eta(d)} \Hom_{\Dcal} ( m\otimes d, F(c) ).
\]
To show that the above is an equivalence, it suffices to show that $\id_m \otimes \eta(d)$ exhibits $m \otimes F^L(d)$ as left adjoint to $F$ at $m \otimes d$. This is implied by the fact that $F^L$ is a strict morphism of $\Mcal$-modules.
 \end{proof}
 
 \begin{remark}\label{remark taubar}
Let $\Mcal$ be a presentable symmetric monoidal category. Then the unit map $\Spc \rightarrow \Mcal$ induces a symmetric monoidal colimit preserving functor  of presentable symmetric monoidal categories $i: \Cat \rightarrow \Cat^\Mcal$. We can think about this as a morphism in $\Cat\modd(\Pr^L)$. Applying the functor $\theta_{\Cat}$ yields an enhancement  of $i$ to a symmetric monoidal functor of symmetric monoidal $2$-categories $\overline{i}: \Catscr \rightarrow \Catscr^\Mcal$. It follows from proposition \ref{prop adjuntos enriched} that $\overline{i}$ admits a right adjoint
\[
\overline{(\tau_\Mcal)_!} : \Catscr^\Mcal \rightarrow \Catscr.
\]
 \end{remark}
 
\begin{proposition}\label{prop overlinetau}
Let $\Mcal$ be a presentable symmetric monoidal category. Then the functor of categories $\Cat^\Mcal \rightarrow \Cat$ underlying  $\overline{(\tau_{\Mcal})_!}$ is equivalent to $\tau_{\Mcal}$.
\end{proposition} 
\begin{proof}
We continue with the notation from remark \ref{remark taubar}. Let 
\[
\epsilon: \overline{i} \overline{(\tau_{\Mcal})_!} \rightarrow \id_{\Catscr^\Mcal}
\] 
be the counit of the adjunction. We think about $\epsilon$ as a morphism in $\Funct(\Catscr^\Mcal, \Catscr^\Mcal)$. 

Consider the functor
\[
((-)^{\leq 1})_*: \Hom_{\twoCat}(\Catscr^\Mcal, \Catscr^\Mcal) \rightarrow \Hom_{\Cat}(\Cat^\Mcal, \Cat^\Mcal).
\]
This admits an enhancement to a functor
\[
\varphi: \Funct(\Catscr^\Mcal, \Catscr^\Mcal)^{\leq 1} \rightarrow \Funct(\Cat^\Mcal, \Cat^\Mcal)
\]
induced from the composite map
\[
\Cat^\Mcal \times \Funct(\Catscr^\Mcal, \Catscr^\Mcal)^{\leq 1} = (\Catscr^\Mcal \times \Funct(\Catscr^\Mcal, \Catscr^\Mcal))^{\leq 1}  \xrightarrow{\ev^{\leq 1}} (\Catscr^\Mcal)^{\leq 1} = \Cat^\Mcal.
\]

 The image of $\epsilon$ under $\varphi$ is a map
\[
\varphi(\epsilon): (\overline{i} \overline{(\tau_\Mcal)_!})^{\leq 1} = i(\overline{(\tau_\Mcal)_!})^{\leq 1}  \rightarrow \id_{\Cat^\Mcal}.
\]
For each object $x$ in $\Cat^{\Mcal}$, the morphism 
\[
\varphi(\epsilon)(x): i (\overline{(\tau_\Mcal)_!})^{\leq 1}(x) \rightarrow x
\]
can be identified with $\epsilon(x)$. Using corollary \ref{coro functorialidad right adjoints} we conclude that $\varphi(\epsilon)$ exhibits $(\overline{(\tau_\Mcal)_!})^{\leq 1}$ as right adjoint to $i$, as desired.
\end{proof} 
 
 \begin{corollary}\label{corollary adjoints of enriched vs ordinary}
 Let $\Mcal$ be a presentable symmetric monoidal category, and let $F: \Ccal \rightarrow \Dcal$ be a functor of $\Mcal$-enriched categories. Assume that $F$ admits a right adjoint $F^R$. Then $(\tau_\Mcal)_! (F^R)$ is right adjoint to $(\tau_\Mcal)_!(F)$. 
 \end{corollary}
 \begin{proof}
 This follows directly from proposition \ref{prop overlinetau}, since functors of $2$-categories preserve adjunctions.
 \end{proof}

We finish by studying the existence of conical limits and colimits in enriched categories arising from presentable modules.

\begin{notation}
Let $\Mcal$ be a presentable monoidal category and let $\Ical$ be a category. We denote by $\Ical_\Mcal'$ the image of $\Ical$ under the composite functor
\[
\Cat \xrightarrow{s} \Algbrd(\Spc) \xrightarrow{} \Algbrd(\Mcal)
\]
where the first map is the functor $s$ from construction \ref{construction inclusion cat}, and the second map is given by pushforward along the unit map $\Spc \rightarrow \Mcal$.
\end{notation}

\begin{lemma}\label{lemma modules over Icalprime}
Let $\Mcal$ be a presentable monoidal category. Let $\Dcal$ be a presentable $\Mcal$-module and let $\Ical$ be a category. Then the projection $\LMod \rightarrow \operatorname{Arr}_{\text{\normalfont oplax}}(\Cat)$ induces an equivalence
\[
\LMod_{\Ical'_\Mcal}(\Dcal) = \Funct(\Ical, \Dcal).
\] 
\end{lemma}
\begin{proof}
As in the proof of proposition \ref{prop restrict scalars lmod}, we let $\Mcal_\Ical$ be the $\Assos$-operad with the universal map $\Mcal_\Ical \times_{\Assos }\Assos_{\Ical} \rightarrow \Mcal$. Recall from \cite{Hinich} that this is a presentable monoidal category which acts on $\Funct(\Ical, \Dcal)$. The category $\LMod_{\Ical'_\Mcal}(\Dcal)$ is then the category of modules for an algebra in $\Mcal_\Ical$. As discussed in \cite{Hinich} 4.7, this algebra is in fact the unit in $\Mcal_X$. We conclude that its category of modules is equivalent to $\Funct(\Ical, \Dcal)$, as desired.
\end{proof}

\begin{lemma}\label{lemma restr of scalars no cambia}
Let $\Mcal$ be a presentable  monoidal category. Let $\Dcal$ be a presentable $\Mcal$-module and let $\Ical$ be a category. Then restriction of scalars along the canonical map $\Ical_\Mcal \rightarrow \Ical'_\Mcal$ induces an equivalence
\[
\LMod_{\Ical'_\Mcal}(\Dcal) = \LMod_{\Ical_\Mcal}(\Dcal).
\]
\end{lemma}
\begin{proof}
We continue with the notation from the proof of lemma \ref{lemma modules over Icalprime}. Recall from \cite{Hinich} 4.4.10 and 4.7.1 that the unit map $\Spc \rightarrow \Mcal$ induces a symmetric monoidal functor $\Spc_\Ical \rightarrow \Mcal_\Ical$, where $\Spc_\Ical$ is defined as $\Mcal_\Ical$. The monoidal category $\Spc_\Ical$ thus acts on $\Funct(\Ical, \Dcal)$ by restriction of scalars, and the category  $\LMod_{\Ical'_\Mcal}(\Dcal)$ is the category of modules over the algebra in $\Spc_\Ical$ associated to $\Ical'_{\Spc}$.  Similarly, we have that $\LMod_{\Ical_\Mcal}(\Dcal)$ is equivalent to the category of modules in $\Funct(\Ical^{\leq 0}, \Dcal)$ for the algebra in $\Spc_{\Ical^{\leq 0}}$ associated to $\Ical_{\Spc}$. 

It follows that to prove our lemma it suffices to assume that $\Mcal = \Spc$. We claim that for every category $\Ecal$ the induced functor
\[
\Hom_{\Cat}(\Ecal, \LMod_{\Ical'_\Mcal}(\Dcal)) \rightarrow \Hom_{\Cat}(\Ecal, \LMod_{\Ical_\Mcal}(\Dcal))
\]
is an equivalence. Using proposition \ref{prop theta el caso de spc} together with \cite{Hinich} proposition 6.3.7 we see that the above is equivalent to the canonical map
\[
\Hom_{\Algbrd(\Spc)}(\Ecal \times \Ical'_\Mcal, \theta'_{\Spc}(\Dcal)) \rightarrow \Hom_{\Algbrd(\Spc)}(\Ecal \times \Ical_\Mcal, \theta'_{\Spc}(\Dcal)).
\]
This is an isomorphism thanks to propositions \ref{prop el caso de spc} and \ref{prop compare defs}.
\end{proof}

\begin{proposition}\label{prop existen conicals}
Let $\Mcal$ be a presentable symmetric monoidal category and let $\Ccal$ be a presentable module over $\Mcal$. Then $\theta_\Mcal(\Ccal)$ is conically complete and cocomplete.
\end{proposition}
\begin{proof}
Let $\Ical$ be a small category. We have to show that the diagonal map
\[
\Delta: \theta_\Mcal(\Ccal ) \rightarrow \Funct(\Ical_\Mcal, \theta_\Mcal(\Ccal))
\]
admits both right and left adjoints.

 Let $h_\Mcal^R : \Algbrd(\Mcal ) \rightarrow \Algbrd(\Mcal)_{\Spc}$ be the colocalization functor. We note that $\Delta$ is the image under $h_\Mcal^R$ of the morphism of algebroids
 \[
 \Delta': \theta'_\Mcal(\Ccal) \rightarrow \Funct(\Ical_\Mcal, \theta'_\Mcal(\Ccal))
 \]
 of precomposition with the projection $\Ical_\Mcal \rightarrow 1_\Mcal$. Using the equivalence from \cite{Hinich} proposition 6.3.7, we see that the above map can be rewritten as
\[
\theta'_\Mcal(\Delta_{\text{mod}}) : \theta'_\Mcal(\Ccal) \rightarrow \theta'_\Mcal(\LMod_{\Ical_\Mcal}(\Ccal))
\]
 where: 
\begin{itemize}
\item The category $\LMod_{\Ical_\Mcal}(\Ccal)$ is equipped with the structure of presentable $\Mcal$-module by virtue of its realization as a category of left modules for an algebra in the $\Assos^-$-component of the $\BM$-monoidal category $\Funct_{\BM}(\BM_{\Ical^{\leq 0}, [0]}, \Ccal)$ (where we consider $\Ccal$ as a $\Mcal-\Mcal$-bimodule in the canonical way).
\item The functor $\Delta_{\text{mod}}$ denotes the functor of restriction of scalars
\[
\Ccal = \LMod_{1_\Mcal}(\Ccal) \rightarrow \LMod_{\Ical_\Mcal}(\Ccal)
\]
along the projection $\Ical_\Mcal \rightarrow 1_\Mcal$, equipped with its canonical structure of morphism of $\Mcal$-modules.
\end{itemize} 

Using lemma \ref{lemma restr of scalars no cambia} we may rewrite our map as
\[
\theta_\Mcal'(\Delta'_{\text{mod}}): \theta'_\Mcal(\Ccal) \rightarrow \theta_\Mcal'(\LMod_{\Ical_{\Mcal}'}(\Ccal))
\]
where $\Delta'_{\text{mod}}$ is defined as $\Delta_{\text{mod}}$, except that using the projection $\Ical'_\Mcal \rightarrow 1_\Mcal$ instead. 

Applying lemma \ref{lemma modules over Icalprime} we see that the above is equivalent to
\[
\theta'_\Mcal(\Delta_{\text{funct}}) : \theta'_\Mcal(\Ccal) \rightarrow \theta'_\Mcal (\Funct(\Ical, \Ccal))
\]
where $\Funct(\Ical, \Ccal)$ is equipped with its canonical structure of module over $\Mcal$, and $\Delta_{\text{funct}}$ denotes the diagonal functor  $\Ccal \rightarrow \Funct(\Ical, \Ccal)$.

We conclude that our original map $\Delta$ is equivalent to $\theta_\Mcal(\Delta_{\text{funct}})$. This admits both left and right adjoints thanks to proposition \ref{prop adjuntos enriched}.
\end{proof}



\tableofcontents

\section{Presentable $n$-categories}

Let $\Mcal$ be a presentable symmetric monoidal category. Then any presentable category $\Ccal$ with an action of $\Mcal$ compatible with colimits has a natural structure of $\Mcal$-enriched category. As we discussed in section 3, this procedure gives rise to a lax symmetric monoidal functor from the  symmetric monoidal category $\Mcal\modd_\pr$ of $\Mcal$-modules in $\Pr^L$ to the category $\Cat^{\Mcal}$ of $\Mcal$-enriched categories. Under certain conditions, a category equipped with an action of $\Mcal\modd_\pr$ inherits the structure of an $\Mcal\modd_\pr$-enriched category, and it therefore also has the structure of a category enriched in $\Cat^{\Mcal}$ - in other words, a $2$-category enriched in $\Mcal$. Our goal in this section is to  introduce higher versions of this procedure that allow one to produce $n$-categories enriched in $\Mcal$ for all values of $n$.

A naive iteration of the enrichment procedure described above runs into set theoretical difficulties: even though $\Mcal$ is presentable, the category $\Mcal\modd_\pr$ is no longer presentable in general as it is generally not locally small. It therefore does not make sense to consider its category of modules in $\Pr^L$. However, we still have that $\Mcal\modd_\pr$ admits all small colimits and its tensor product structure preserves all colimits, so one may consider the symmetric monoidal category $\Mcal\modd^2$ of $\Mcal\modd_\pr$-modules in the symmetric monoidal category $\widehat{\Cat}_\cocompl$ of large cocomplete categories. One can then define a functor from $\Mcal\modd^2$ to the category of $2$-categories enriched in a suitable completion of $\Mcal$.

If one attempts to repeat the above to obtain enriched $3$-categories one runs into a similar problem than in the $2$-categorical case. Again the category $\Mcal\modd^2$ is not presentable. However, $\Mcal\modd^2$ is now not large, but very large. It therefore does not make sense to consider its category of modules in $\widehat{\Cat}_\cocompl$. One way around this problem is to consider instead its category of modules in the category $\operatorname{CAT}_\cocompl$ of very large categories admitting all large colimits. This can be iterated to yield a theory that works for all values of $n$, however this requires one to work with an infinite sequence of nested universes and keep careful track of the relative sizes of various objects.

In this paper we pursue a different approach that is based on the observation that although $\Mcal\modd^2$ is very large, it is controlled by a large cocomplete category $\Mcal\modd^2_\pr$ of so called presentable $\Mcal\modd_\pr$-modules. It therefore makes sense to consider the symmetric monoidal category $\Mcal\modd^3$ of modules for $\Mcal\modd_\pr^2$ inside $\widehat{\Cat}_\cocompl$, which is again controlled by a subcategory $\Mcal\modd_\pr^3$ of presentable $\Mcal\modd_\pr^2$-modules. Iterating this reasoning yields symmetric monoidal categories $\Mcal\modd^n$ and $\Mcal\modd_\pr^n$ for all $n \geq 1$, which admit functors into the category of $n$-categories enriched in a suitable completion of $\Mcal$.

We begin in \ref{subsection presentable modules} by introducing the notion of presentable module over an arbitrary monoidal cocomplete category $\Ecal$. This is obtained as a special case of the notion of $\kappa$-compactness in presentable categories in the case when $\kappa = \kappa_0$ is the smallest large cardinal, applied to the very large presentable category $\Ecal\modd(\widehat{\Cat}_\cocompl)$. We show that if $A$ is an algebra object in $\Ecal$ then the category of right $A$-modules is a presentable left $\Ecal$-module. In the case when $\Ecal$ is a presentable monoidal category, we show that the category $\Ecal\modd_\pr$ of presentable $\Ecal$-modules agrees with the category $\Ecal\modd(\Pr^L)$ of $\Ecal$-modules in $\Pr^L$. When $\Ecal =\Spc$ this yields a characterization of $\Spc\modd_\pr = \Pr^L$ as the full subcategory of $\widehat{\Cat}_\cocompl$ on the $\kappa_0$-compact objects. 

In \ref{subsection higher mod} we iterate the above to obtain symmetric monoidal categories $\Ecal\modd^n$ and $\Ecal\modd^n_{\pr}$ attached to each symmetric monoidal cocomplete category $\Ecal$. Specializing to the case $\Ecal = \Spc$, we obtain a symmetric monoidal category $n\kr\Pr^L = \Spc\modd^n_\pr$ whose objects we call presentable $n$-categories. In the case when $\Ecal$ is the category of spectra, this yields a full subcategory $n\kr\Pr^L_{\St}$ of $n\kr\Pr^L$ whose objects we call presentable stable $n$-categories. Applying the results from \ref{subsection presentable modules} inductively, we are able to conclude that for every commutative ring spectrum $A$ there is an associated presentable stable $n$-category $A\modd^n$ of $A$-linear presentable stable $(n-1)$-categories.

In \ref{subsection ncat structure} we use the theory from section 3 to construct a lax symmetric monoidal functor $\psi_n: \Spc\modd^n \rightarrow \reallywidehat{\nCat}$. In other words, even though our approach to presentable $n$-categories is $1$-categorical in nature, we are able to upgrade presentable $n$-categories to honest $n$-categories.  We use the realization functor $\psi_n$ to obtain $(n+1)$-categorical enhancements $n\Prscr^L$ and $n\Catscr^L$ for $n\kr\Pr^L$ and $\Spc\modd^n$, respectively.

 In \ref{subsection conical en pres} we establish our first main result regarding the theory of presentable $n$-categories, theorem \ref{teo psi tiene conicals}: if $\Ccal$ belongs to $\Spc\modd^n$ then  $\psi_n(\Ccal)$ is a conically cocomplete $n$-category, and  moreover any limit that exists in $\Ccal$ is a conical limit in $\psi_n(\Ccal)$. In particular, we conclude that $n\Prscr^L$ is a conically cocomplete $n$-category.

A fundamental feature of the $2$-category $\Prscr^L$ is that colimits of right adjointable diagrams can be computed as limits after passage to right adjoints. In \ref{subsection passage} we establish a generalization of this statement (theorem \ref{teo passage}) : we show that for every object $\Ccal$ in $\Spc\modd^n$ the $n$-category $\psi_n(\Ccal)$ satisfies the so-called passage to adjoints property. In particular, we are able to conclude that $\psi_n(\Ccal)$ admits limits of left adjointable diagrams.

\subsection{Presentable modules over cocomplete monoidal categories}\label{subsection presentable modules}

We begin by discussing the notion of smallness for objects of a category with large colimits. This is a special case of the notion of $\kappa$-compactness from \cite{HTT} section 5.3.4, in the case when $\kappa$ is the smallest large cardinal.

\begin{notation}
Denote by $\kappa_0$ the smallest large cardinal. In other words, $\kappa_0$ is such that $\kappa_0$-small spaces are what we usually call small spaces. We denote by $\widehat{\Cat}_\cocompl$ the category of large categories admitting all small colimits, and functors which preserve small colimits. As usual, we denote by $\Pr^L$ the full subcategory of $\widehat{\Cat}_\cocompl$ on the presentable categories.
\end{notation}

\begin{definition}
Let $\Ccal$ be a very large, locally large category admitting all large colimits. An object $c$ in $\Ccal$ is said to be $\kappa_0$-compact if the functor $\Ccal \rightarrow \widehat{\Spc}$ corepresented by $c$ preserves large $\kappa_0$-filtered colimits. We say that $\Ccal$ is $\kappa_0$-compactly generated if it is generated under large colimits by its $\kappa_0$-compact objects, and the space of $\kappa_0$-compact objects in $\Ccal$ is large.
\end{definition}

\begin{remark}
Let $\Ccal$ be a very large, locally large category admitting all large colimits. Then the collection of $\kappa_0$-compact objects is closed under small colimits in $\Ccal$ (\cite{HTT} corollary 5.3.4.15). Moreover, if $\Ccal$ is $\kappa_0$-compactly generated then it can be recovered from its full subcategory of $\kappa_0$-compact objects by freely adjoining large $\kappa_0$-filtered colimits.
\end{remark}

\begin{proposition}\label{prop small objects in catcocompl}
The category $\widehat{\Cat}_{\normalfont \cocompl}$ is $\kappa_0$-compactly generated. Moreover, an object of $\widehat{\Cat}_{\normalfont \cocompl}$ is $\kappa_0$-compact if and only if it is a presentable category. 
\end{proposition}
\begin{proof}
The fact that $\widehat{\Cat}_\cocompl$ is presentable as a very large category follows from \cite{HA} lemma 4.8.4.2. Let $U: \widehat{\Cat}_\cocompl \rightarrow \widehat{\Cat}$ be the forgetful functor. It follows from \cite{HTT} proposition 5.5.7.11 that $U$ preserves $\kappa_0$-filtered colimits, and therefore its left adjoint $\Pcal: \widehat{\Cat} \rightarrow \widehat{\Cat}_\cocompl$ sends $\kappa_0$-compact objects to $\kappa_0$-compact objects. Note that $\widehat{\Cat}$ is $\kappa_0$-compactly generated by the small categories. Since $U$ is furthermore conservative we conclude that $\widehat{\Cat}_\cocompl$ is generated under large colimits by the objects of the form $\Pcal(\Ical)$ where $\Ical$ is a small category. In particular we conclude that $\widehat{\Cat}_\cocompl$ is $\kappa_0$-compactly generated.

We now show that $\Pr^L$ is closed under small colimits inside $\widehat{\Cat}_\cocompl$. Recall from \cite{HTT} section 5.5.7 that $\Pr^L$ is the union of its subcategories $\Pr^L_\kappa$ of $\kappa$-compactly generated categories and functors that preserve $\kappa$-compact objects, for $\kappa$ a small regular cardinal. The functor $\Pr^L_\kappa \rightarrow \Cat$ sending each object $\Ccal$ to its subcategory of $\kappa$-compact objects induces an equivalence between $\Pr^L_\kappa$ and the subcategory $\Cat^{\rex(\kappa), \id}$ of $\Cat$ consisting of idempotent complete categories with $\kappa$-small colimits and functors that preserve those colimits. It follows from \cite{HA} lemma 4.8.4.2 that $\Pr^L_\kappa$ is presentable and in particular it has all small colimits. 

For $\tau > \kappa$ a pair of regular cardinals, the inclusion $\Pr^L_\kappa \rightarrow \Pr^L_\tau$ is equivalent to the functor $\Cat^{\rex(\kappa),\id} \rightarrow \Cat^{\rex(\tau), \id}$ which freely adjoins $\tau$-small $\kappa$-filtered colimits. This is a left adjoint, and so it follows that the inclusion $\Pr^L_\kappa \rightarrow \Pr^L_\tau$ preserves small colimits. Similarly, for each $\kappa$ the inclusion $\Pr^L_\kappa \rightarrow \widehat{\Cat}_\cocompl$ is equivalent to the functor $\Cat^{\rex(\kappa),\id} \rightarrow \widehat{\Cat}_\cocompl$ that freely adjoins small $\kappa$-filtered colimits, and so it preserves small colimits.

Note that $\Pr^L$ is the colimit of the categories $\Pr^L_\kappa$ in $\widehat{\Cat}$, and in particular also in the category of very large categories.	 Moreover, the inclusion $\Pr^L \rightarrow \Catcocompl$ is induced by passing to the colimit the inclusions $\Pr^L_\kappa \rightarrow \Catcocompl$.  It follows from \cite{HTT} proposition 5.5.7.11  that  $\Pr^L$ admits small colimits and the inclusion $\Pr^L \rightarrow \Catcocompl$ preserves small colimits. In other words, $\Pr^L$ is closed under small colimits in $\widehat{\Cat}_{\cocompl}$, as we claimed.

It remains to show that $\Pr^L$ is generated under small colimits by the objects $\Pcal(\Ical)$ for $\Ical$ in $\Cat$. Indeed, for each small regular cardinal $\kappa$ the forgetful functor $\Pr^L_\kappa \rightarrow \Cat$ is conservative and preserves $\kappa$-filtered colimits, which implies that $\Pr^L_\kappa$ is generated under colimits by the image of $\Pcal|_{\Cat}$. Our claim now follows from the fact that $\Pr^L$ is the union of the subcategories $\Pr^L_\kappa$, and that the inclusions $\Pr^L_\kappa \rightarrow \Pr^L$ preserve small colimits.
\end{proof}

 Recall from \cite{HA} section 4.8 that $\widehat{\Cat}_{\cocompl}$ has a symmetric monoidal structure $\otimes$ such that for each triple of objects $\Ccal, \Dcal ,\Ecal$, the space of morphisms $\Ccal \otimes \Dcal \rightarrow \Ecal$ is equivalent to the space of functors $\Ccal \times \Dcal \rightarrow \Ecal$ which are colimit preserving in each coordinate. The unit of this symmetric monoidal structure is the category $\Spc$ of spaces. This symmetric monoidal structure is compatible with large colimits, in the sense that the functor 
 \[
 \otimes: \widehat{\Cat}_\cocompl \times \widehat{\Cat}_\cocompl \rightarrow \Catcocompl
 \]
  preserves large colimits in each variable. Moreover, the subcategory $\Pr^L$ is closed under tensor product and therefore $\Pr^L$ also inherits a symmetric monoidal structure which is compatible with small colimits.

\begin{notation}
Let $\Ecal$ be an algebra object in $\widehat{\Cat}_\cocompl$ (in other words, $\Ecal$ is a large cocomplete category equipped with a monoidal structure compatible with colimits). We denote by $\Ecal\modd^l$ (resp. $\Ecal\modd^r$) the category of left (resp. right) $\Ecal$-modules in $\widehat{\Cat}_\cocompl$.

 For each algebra $A$ in $\Ecal$ we denote by $A\modd^l$ (resp. $A\modd^r$) the category of left (resp. right) $A$-modules, thought of as an object of $\Ecal\modd^r$ (resp. $\Ecal\modd^l$). 

If $\Ecal$ is a commutative algebra, we will use the notation $\Ecal\modd$ for $\Ecal\modd^l = \Ecal\modd^r$. If $A$ is a commutative algebra in $\Ecal$, we will use the notation $A\modd$ instead of $A\modd^l = A\modd^r$.
\end{notation}

\begin{warning}
Let $\Ecal$ be an algebra object in $\widehat{\Cat}_\cocompl$, and assume that $\Ecal$ is furthermore presentable. Then $\Ecal\modd^l$ has two different meanings. If we think about $\Ecal$ as an  algebra in $\widehat{\Cat}_{\cocompl}$ then $\Ecal\modd^l$ denotes the category of cocomplete categories with an action of $\Ecal$. However if we think about $\Ecal$ as an algebra in $\Pr^L$ (which is itself a commutative algebra in $\widehat{\Cat}_{\cocompl}$) then $\Ecal\modd^l$ can be taken to mean the category of presentable categories with an action of $\Ecal$. We will consider modules in $\widehat{\Cat}_{\cocompl}$ by default, unless it is clear from context that $\Ecal$ is to be considered specifically as an algebra in $\Pr^L$. When in doubt, we will make the context clear in our notation by writing $\Ecal\modd^l(\Catcocompl)$ or $\Ecal\modd^l_{\pr} = \Ecal\modd^l(\Pr^L)$.
\end{warning}

\begin{proposition}\label{prop kappa0 compact Ecal modules}
Let $\Ecal$ be an algebra object in $\widehat{\Cat}_{\normalfont \cocompl}$. Then the category  $\Ecal\modd^l$ is $\kappa_0$-compactly generated. Moreover, the collection of free left $\Ecal$-modules $\Ecal \otimes \Ccal$ where $\Ccal$ is a presentable category is a collection of $\kappa_0$-compact generators for $\Ecal\modd^l$.
\end{proposition}
\begin{proof}
The forgetful functor $\Ecal\modd^l \rightarrow \widehat{\Cat}_\cocompl$ preserves large colimits by \cite{HA} corollary 4.2.3.5, and therefore its left adjoint preserves $\kappa_0$-compact objects. Combining this with proposition
\ref{prop small objects in catcocompl} we conclude that the free left $\Ecal$-module $\Ecal \otimes \Ccal$ is $\kappa_0$-compact for every presentable category $\Ccal$. The fact that these objects generate $\Ecal\modd^l$ under large colimits follows by observing that $\Ecal\modd^l$ is generated under colimits by free modules (\cite{HA} proposition 4.7.3.14), and that free modules live in the colimit-closure of the modules $\Ecal \otimes \Ccal$ with $\Ccal$ presentable. Using \cite{HA} 4.2.3.7 we conclude that $\Ecal\modd^l$ is $\kappa_0$-compactly generated. 
\end{proof}

\begin{definition}
Let $\Ecal$ be an algebra in $\widehat{\Cat}_{\normalfont \cocompl}$. We say that a left $\Ecal$-module is presentable if it is a $\kappa_0$-compact object in $\Ecal\modd^l$. We denote by $\Ecal\modd^l_\pr$ the full subcategory of $\Ecal\modd^l$ on the presentable $\Ecal$-modules. If $\Ecal$ comes equipped with the structure of a commutative algebra, we will use the notation $\Ecal\modd_\pr$ instead of $\Ecal\modd^l_{\normalfont \pr}$.
\end{definition}

\begin{remark}\label{remark generation}
Recall from the proof of proposition \ref{prop small objects in catcocompl} that $\Pr^L$ is generated under small colimits by presheaf categories. Since the category of small categories is generated under small colimits by $[1]$, we in fact have that $\Pr^L$ is generated under small colimits by $\Pcal([1])$. Using proposition \ref{prop kappa0 compact Ecal modules} we see that for any algebra $\Ecal$ in $\widehat{\Cat}_\cocompl$ the category $\Ecal\modd_\pr$ is the smallest category of left $\Ecal$-modules containing $\Ecal \otimes \Pcal([1]) = \Funct([1], \Ecal)$ and closed under small colimits.
\end{remark}

The following proposition provides an abundant source of presentable $\Ecal$-modules.

\begin{proposition}\label{prop modules pres}
Let $\Ecal$ be an algebra in $\widehat{\Cat}_\cocompl$, and let $A$ be an algebra in $\Ecal$. Then the category $A\modd^r$ of right $A$-modules equipped with its natural left $\Ecal$-module structure, is a presentable $\Ecal$-module.
\end{proposition}
\begin{proof}
Recall from \cite{HA} remark 4.8.4.8 that $A\modd^r$ is left dualizable as a left $\Ecal$-module - its left dual is the category $A\modd^l$ of left $A$-modules with its natural right $\Ecal$-module structure. We thus have an equivalence
\[
\Hom_{\Ecal\modd^l}(A\modd^r, -) = \Hom_{\widehat{\Cat}_\cocompl}(\Spc, A\modd^l \otimes_{\Ecal} -) : \Ecal\modd^l \rightarrow \widehat{\Spc}.
\]
It follows from a combination of \cite{HA} proposition 4.4.2.14 and \cite{HTT} proposition 5.5.7.11 that the right hand side preserves $\kappa_0$-filtered colimits, and thus $A\modd^r$ is $\kappa_0$-compact.
\end{proof}

\begin{remark}\label{remark functorialidad presentability} \label{remark endofunctor}
Let $\varphi: \Ecal \rightarrow \Ecal'$ be a morphism of algebras in $\widehat{\Cat}_\cocompl$. Then the extension of scalars functor
\[
\varphi_* = \Ecal' \otimes_\Ecal - : \Ecal\modd^l \rightarrow \Ecal'\modd^l
\]
preserves large colimits, and maps each free $\Ecal$-module $\Ecal \otimes \Ccal$ to the free $\Ecal'$-module $\Ecal' \otimes \Ccal$. It follows from proposition \ref{prop kappa0 compact Ecal modules} that $\varphi_*$ maps presentable $\Ecal$-modules to presentable $\Ecal'$-modules.

Similarly, given a pair of algebras $\Ecal, \Ecal'$ in $\widehat{\Cat}_\cocompl$, the exterior product functor
\[
\boxtimes : \Ecal\modd^l \times \Ecal'\modd^l \rightarrow (\Ecal \otimes \Ecal')\modd^l
\]
preserves colimits in each variable, and maps the pair $(\Ecal \otimes \Ccal, \Ecal' \otimes \Dcal)$ to $\Ecal \otimes \Ecal' \otimes \Ccal \otimes \Dcal$. Using proposition \ref{prop kappa0 compact Ecal modules} we conclude that $\boxtimes$ restricts to the subcategories of presentable objects.
 
 We therefore conclude that the assignment $\Ecal \mapsto \Ecal\modd^l_\pr$ has a natural enhancement to a lax symmetric monoidal functor
\[
\modd_\pr^l : \Alg(\widehat{\Cat}_\cocompl) \rightarrow \widehat{\Cat}_\cocompl.
\]

In particular, passing to commutative algebra objects yields an endofunctor
\[
\modd_\pr : \CAlg(\widehat{\Cat}_\cocompl) \rightarrow \CAlg(\widehat{\Cat}_\cocompl).
\]
In other words, if $\Ecal$ is a commutative algebra in $\widehat{\Cat}_\cocompl$, then the natural symmetric monoidal structure on $\Ecal\modd$ restricts to a symmetric monoidal structure on $\Ecal\modd_\pr$.
\end{remark}

\begin{example}
Let $\Ecal = \Spc$ with the cartesian symmetric monoidal structure. This is the trivial commutative algebra in $\Catcocompl$ and therefore $\Ecal\modd(\Catcocompl) =\Catcocompl$. Using proposition \ref{prop small objects in catcocompl} we see that $\Ecal\modd_\pr$ coincides with the category $\Pr^L$ of presentable categories. The symmetric monoidal structure on $\Ecal\modd_\pr$ agrees with the usual symmetric monoidal structure on $\Pr^L$.
\end{example}

\begin{remark}\label{remark Amod commutative}
Let $\Ecal$ be a commutative algebra in $\widehat{\Cat}_\cocompl$. Recall from \cite{HA} section 4.8.5 that the functor $\operatorname{-mod}^r:\Alg(\Ecal) \rightarrow \Ecal\modd$
has a natural symmetric monoidal structure. Combining propositions \ref{prop modules pres} and remark \ref{remark functorialidad presentability} we see that the above restricts to a symmetric monoidal functor $\Alg(\Ecal) \rightarrow \Ecal\modd_\pr$. In particular, we obtain a symmetric monoidal functor
\[
\modd: \CAlg(\Ecal) \rightarrow \CAlg(\Ecal\modd_\pr).
\]
\end{remark}

We finish by studying the interaction of the notion of presentability with extensions of commutative algebras.

\begin{proposition}\label{prop presentable modules over extension}
Let $\varphi: \Ecal \rightarrow \Ecal'$ be a morphism of commutative algebras in $\widehat{\Cat}_\cocompl$, so that $\Ecal'$ can alternatively be thought of as a commutative algebra in $\Ecal\modd$. Assume that $\Ecal'$ is presentable as an $\Ecal$-module.  Then the canonical equivalence of symmetric monoidal categories
\[
\Ecal'\modd(\widehat{\Cat}_\cocompl) = \Ecal'\modd(\Ecal\modd)
\]
restricts to an equivalence
\[
\Ecal'\modd_\pr = \Ecal'\modd(\Ecal\modd_\pr).
\]
\end{proposition}
\begin{proof}
By proposition \ref{prop kappa0 compact Ecal modules}, the subcategory $\Ecal'\modd_\pr$ of $\Ecal'\modd(\widehat{\Cat}_\cocompl)$ is generated under small colimits by the free modules $\Ecal' \otimes \Ccal$ with $\Ccal$ presentable. Similarly, $\Ecal\modd_\pr$ is generated under small colimits by the free modules $\Ecal \otimes \Ccal$ with $\Ccal$ presentable, and thus $\Ecal'\modd(\Ecal\modd_\pr)$ is generated under small colimits by the objects of the form 
\[
\Ecal' \otimes_\Ecal (\Ecal \otimes \Ccal) = \Ecal' \otimes \Ccal.
\] 
We thus see that $\Ecal'\modd_\pr$ and $\Ecal'\modd(\Ecal\modd_\pr)$ are generated under small colimits by the same collection of objects, and they therefore agree.
\end{proof}

\begin{corollary}\label{coro Emodd para presentable}
Let $\Ecal$ be a commutative algebra in $\Pr^L$. Then $\Ecal\modd_\pr$ is equivalent to the category $\Ecal\modd(\Pr^L)$ of presentable categories equipped with an $\Ecal$-module structure.
\end{corollary}
\begin{proof}
Apply proposition \ref{prop presentable modules over extension} in the case of the unit $\Spc \rightarrow \Ecal$.
\end{proof}

\subsection{Higher module categories}\label{subsection higher mod}
We now iterate the functor $\operatorname{-mod}_\pr$ to arrive at a notion of higher presentable modules over a cocomplete symmetric monoidal category.

\begin{notation}
For each $n \geq 0$ we let 
\[
\operatorname{-mod}_\pr^n : \CAlg(\widehat{\Cat}_\cocompl) \rightarrow \CAlg(\widehat{\Cat}_\cocompl)
\]
be the endofunctor obtained by composing the endofunctor $\modd_\pr$ of remark \ref{remark endofunctor} with itself $n$ times (in particular, if $n = 0$ we set $\operatorname{-mod}_\pr^0$ to be the identity endofunctor).

Denote by 
\[
\operatorname{-mod}: \CAlg(\widehat{\Cat}_\cocompl) \rightarrow \CAlg(\operatorname{CAT})
\]
 the functor that assigns to each commutative algebra $\Ecal$ in $\widehat{\Cat}_\cocompl$  the (very large) symmetric monoidal category $\Ecal\modd$. For each $n \geq 1$ we let $\operatorname{-mod}^n$ be the composite functor
\[
\CAlg(\widehat{\Cat}_\cocompl) \xrightarrow{\modd^{n-1}_\pr} \CAlg(\widehat{\Cat}_\cocompl) \xrightarrow{\modd} \CAlg(\operatorname{CAT}).
\]

Let $\Ecal$ be a commutative algebra in $\widehat{\Cat}_\cocompl$. We inductively define a symmetric monoidal functor $\operatorname{-mod}^n: \CAlg(\Ecal) \rightarrow \CAlg(\Ecal\modd_\pr^n)$ as follows:
\begin{itemize}
\item When $n = 0 $ we let $\operatorname{-mod}^0$ be the identity functor.
\item When $n > 1$ we let $\operatorname{-mod}^n$ be the composite functor
\[
\CAlg(\Ecal) \xrightarrow{\modd^{n-1}} \CAlg(\Ecal\modd_\pr^{n-1}) \xrightarrow{\modd} \CAlg(\Ecal\modd_\pr^n).
\]
where $\operatorname{-mod}$ denotes the functor of remark \ref{remark Amod commutative}.
\end{itemize}
\end{notation}

\begin{definition}
Let $\Ecal$ be a commutative algebra in $\Catcocompl$ and let $n \geq 0$. We call $\Ecal\modd^n_\pr$ the category of presentable $\Ecal$-linear $n$-categories. In the case $\Ecal = \Spc$ we use the notation $n\kr\Pr^L = \Spc\modd^n_\pr$, and call it the category of presentable $n$-categories.
\end{definition}

\begin{warning}
Let $\Ecal$ be a presentable symmetric monoidal category and let $n \geq 1$. By virtue of being a commutative algebra in $\Catcocompl$, we have associated categories $\Ecal\modd^n$ and $\Ecal\modd^n_\pr$ which are related to each other by passage to $\kappa_0$-compact objects, and Ind-$\kappa_0$-completion.  We can also think about $\Ecal$ as a commutative algebra in the cocomplete symmetric monoidal category $\Pr^L$, and attach to it the object $\Ecal\modd^n$ in $\Pr^L\modd^n_\pr$. These two notations are in conflict: as explained in remark \ref{remark linearity extension} below, the underlying symmetric monoidal category to the latter agrees in fact with $\Ecal\modd^n_\pr$. We hope it will be clear from context which version of the construction we are using at each time, and in case where it may be ambiguous we will use the notation $\Ecal\modd^n(\widehat{\Cat}_\cocompl)$ or $\Ecal\modd^n(\Pr^L)$ to specify which version of the two constructions we are using.  
\end{warning}

\begin{remark} \label{remark linearity extension}
Let $\varphi: \Ecal \rightarrow \Ecal'$ be a morphism of commutative algebras in $\widehat{\Cat}_\cocompl$, and assume that $\Ecal'$ is a presentable $\Ecal$-module. Denote by $\overline{\Ecal'}$ the associated commutative algebra in $\Ecal\modd_\pr$. An inductive application of proposition \ref{prop presentable modules over extension} shows that for every $n \geq 1$ the symmetric monoidal category $\Ecal'\modd^n_\pr$ is equivalent to the symmetric monoidal category underlying $\overline{\Ecal'}\modd^n$. In particular, $\Ecal'\modd^n_\pr$ can be enhanced to a commutative algebra in presentable $\Ecal$-linear $(n+1)$-categories.

Specializing the above to the case when $\Ecal = \Spc$ and $\Ecal'$ is a presentable symmetric monoidal category, we see that $\Ecal'\modd^n_\pr$ has a natural enhancement to a commutative algebra in presentable $(n+1)$-categories. Moreover, if $n > 1$ we have equivalences
\[
\Ecal'\modd^n_\pr = \Ecal'\modd^n(\Pr^L) =( \Ecal'\modd^{n-1}(\Pr^L))\modd(n\kr\Pr^L) = \Ecal'\modd^{n-1}_\pr\modd(n\kr\Pr^L).
\]
In other words, a presentable $\Ecal$-linear $n$-category is the same data as a presentable $n$-category equipped with the action of the presentable symmetric monoidal $n$-category of $\Ecal$-linear $(n-1)$-categories.
\end{remark}

Specializing to the case when $\Ecal$ is the category of spectra we obtain a notion of presentable stable $n$-category.

\begin{definition}
Denote by $\Sp$ the symmetric monoidal category of spectra. For each $n \geq 0$ we denote by $n\kr\Pr^L_{\St}$ the category $\Sp\modd^n_\pr = \Sp\modd^n(\Pr^L)$, and call it the category of presentable stable $n$-categories. For each commutative ring spectrum $A$  we call $A\modd^n$ the presentable $n$-category of presentable $A$-linear $(n-1)$-categories.
\end{definition}

\begin{remark}
Recall from \cite{HA} section 4.8.2 that $\Sp$ is an idempotent commutative algebra in $\Pr^L$. It follows by induction that   $n\kr\Pr^L_{\St}$ is an idempotent commutative algebra in $(n+1)\Pr^L$ for every $n \geq 0$. We conclude that for every $n \geq 1$ the functor
\[
- \otimes (n-1)\kr\Pr^L_\St  :  n\kr\Pr^L \rightarrow n\kr\Pr^L_\St
\]
is a localization functor. We think about the above as the stabilization functor for presentable $n$-categories. The fact that it is a localization implies that being stable is a property of a presentable $n$-category.
\end{remark}

\subsection{The $n$-categorical structure}\label{subsection ncat structure}
Our next goal is to enhance the category $\Spc\modd^n $ to an $(n+1)$-category. We will first need to construct a lax symmetric monoidal functor  
\[
\psi_n: \Spc\modd^n \rightarrow \widehat{\nCat}.
\]
 This is accomplished by applying in an inductive way the functor that turns modules over presentable monoidal categories into enriched categories.

\begin{notation}
For each $n \geq 1$ denote by 
\[
\upsilon_{\Spc\modd^{n}}: \Spc\modd^{n} \rightarrow \widehat{\Cat}
\]
the lax symmetric monoidal functor given by the composition
\[
\Spc\modd^n = (n-1)\kr\Pr^L\modd(\widehat{\Cat}_\cocompl) \rightarrow \widehat{\Cat}_{\cocompl}\rightarrow \widehat{\Cat}
\] 
where the second and third arrow are the forgetful functors. The composition of $\upsilon_{\Spc\modd^n}$ with the symmetric monoidal functor $(-)^{\leq 0}: \widehat{\Cat} \rightarrow \widehat{\Spc}$ is the lax symmetric monoidal functor $\tau_{\Spc\modd^n}$ that maps each object of $\Spc\modd^n$ to its underlying space.
\end{notation}

\begin{construction}\label{construction lax functor}
Let 
\[
\psi_1: \Spc\modd = \widehat{\Cat}_\cocompl \rightarrow \widehat{\Cat}
\]
 be the forgetful functor, equipped with its canonical lax symmetric monoidal structure. In other words, we have $\psi_1 = \upsilon_1$. Let $n > 1$ and assume given a lax symmetric monoidal functor 
\[
\psi_{n-1}: \Spc\modd^{n-1} \rightarrow \reallywidehat{(n-1)\kr\Cat}
\]
 such that the composite lax symmetric monoidal functor
\[
\Spc\modd^{n-1} \xrightarrow{\psi_{n-1}} \reallywidehat{(n-1)\kr\Cat} \xrightarrow{(-)^{\leq 1}} \widehat{\Cat}
\]
is equivalent to $\upsilon_{\Spc\modd^{n-1}}$.

 Denote by $\widehat{\Pr}{}^L$  the category of very large presentable categories and (large) colimit preserving functors, and by $\operatorname{CAT}$ the category of very large categories. We let 
 \[
 \varphi_n: \Spc\modd^{n} \rightarrow  \operatorname{CAT}^{\Spc\modd^{n-1}}
\] 
  be the lax symmetric monoidal functor given by the following composition:
\[
\Spc\modd^n = (n-1)\kr\Pr^L \modd(\widehat{\Cat}_{\text{cocompl}}) \xrightarrow{\Ind_{\kappa_0}}  (\Spc\modd^{n-1})\modd(\widehat{\Pr}{}^L) \rightarrow  \operatorname{CAT}^{\Spc\modd^{n-1}} 
\]
Here the second arrow is the lax symmetric monoidal functor $\theta_{\Spc\modd^{n-1}}$ which sends each very large presentable module over $\Spc\modd^{n-1}$ to its underlying $\Spc\modd^{n-1}$-enriched category.

Note that the lax symmetric monoidal functor $(\tau_{\Spc\modd^{n-1}})_! \varphi_n$ is equivalent to the composite
\[
\Spc\modd^n \xrightarrow{\upsilon_{\Spc\modd^n}} \widehat{\Cat} \xrightarrow{\Ind_{\kappa_0}} \operatorname{CAT}.
\]

In particular, we observe that for each morphism $f: \Ccal \rightarrow \Dcal$ in $\Spc\modd^n$, the induced morphism $\varphi_n(f): \varphi_n(\Ccal) \rightarrow \varphi_n(\Dcal)$ restricts to the full enriched-subcategories of $\varphi_n(\Ccal)$ and $\varphi_n(\Dcal)$ on those objects which are $\kappa_0$-compact, when thought of as an object in $\Ccal$ or $\Dcal$. Furthermore, for every pair of objects $\Ccal, \Dcal$ in $\Spc\modd^n$, the morphism
\[
\varphi_n(\Ccal) \otimes \varphi_n(\Dcal) \rightarrow \varphi_n(\Ccal \otimes \Dcal)
\]
arising from the lax symmetric monoidal structure on $\varphi_n$, restricts to the subcategories of $\kappa_0$-compact objects.  We thus have a well defined lax symmetric monoidal functor
\[
\varphi'_n: \Spc\modd^n \rightarrow \widehat{\Cat}^{\Spc\modd^{n-1}}
\]
equipped with a monomorphism into $\varphi_n$, and which maps every object $\Ccal$ in $\Spc\modd^n$ to the full enriched subcategory of $\varphi_n(\Ccal)$ on those objects which are $\kappa_0$-compact when thought of as objects of $\Ccal$.

Let $\psi'_n$ be the composite lax symmetric monoidal functor
\[
\Spc\modd^n \xrightarrow{\varphi_n'} \widehat{\Cat}^{\Spc\modd^{n-1}} \xrightarrow{(\psi_{n-1})_!} \widehat{\Algbrd}(\reallywidehat{{(n-1)\kr\Cat}}).
\]
Note that for every object $\Ccal$ in $\Spc\modd^n$, the Segal space underlying $\psi'_n(\Ccal)$ is given by $(\psi_{n-1}^{\leq 0})_! \varphi'_n(\Ccal) = (\tau_{\Spc\modd^{n-1}})_! \varphi'_n(\Ccal)$ which is a complete Segal space since $\varphi'_n(\Ccal)$ is an enriched category. It follows that  $\psi'_n(\Ccal)$ is in fact also an enriched category. Hence $\psi'_n$ defines, by corestriction, a lax symmetric monoidal functor
\[
\psi_n : \Spc\modd^n \rightarrow \widehat{\Cat}^{\widehat{(n-1)\kr\Cat}} =  \widehat{\nCat}.
\]

Observe that we have equivalences of lax symmetric monoidal functors
\[
(\psi_n)^{\leq 1} = (-)^{\leq 0}_! \psi_n = ((\psi_{n-1})^{\leq 0})_! \varphi_{n}'  = (\tau_{\Spc\modd^{n-1}})_! \varphi'_n
\]
and the latter is obtained from $(\tau_{\Spc\modd^{n-1}})_! \varphi_n = \Ind_{\kappa_0} \upsilon_{\Spc\modd^n}$ by restricting to $\kappa_0$-compact objects. It follows that $(\psi_n)^{\leq 1}$ is equivalent, as a lax symmetric monoidal functor, to $\upsilon_{\Spc\modd^n}$. We conclude that this construction may be iterated to yield lax symmetric monoidal functors $\psi_n$ for all $n \geq 1$.
\end{construction}

\begin{remark}
Let $n \geq 2$ and let $\Ccal$ be an object of $\Spc\modd^n$. Unwinding construction \ref{construction lax functor} reveals that $\psi_n(\Ccal)$ is an $n$-category whose underlying category is the category underlying $\Ccal$. For each pair of objects $x, y$ of $\Ccal$ we have an equivalence
\[
\Hom_{\psi_n(\Ccal)}(x, y) = \psi_{(n-1)}( \shom_{\Ind_{\kappa_0}(\Ccal)}(x, y) )
\]
where the right hand side denotes the Hom object between $x$ and $y$, where we consider $\Ind_{\kappa_0}(\Ccal)$ as a module over $\Spc\modd^{n-1}$. 
\end{remark}

We now use the functors $\psi_n$ to obtain the desired $(n+1)$-categorical enhancement of $\Spc\modd^n$.

\begin{notation}
Let $n \geq 1$. We denote by $n\Catscr^L$ the (very large) symmetric monoidal $(n+1)$-category obtained by applying the composite lax symmetric monoidal functor
\[
(\Spc\modd^n)\modd(\widehat{\Pr}{}^L) \xrightarrow{\theta_{\Spc\modd^n}} \operatorname{CAT}^{\Spc\modd^n} \xrightarrow{(\psi_n)_!} \operatorname{CAT}^{\widehat{\nCat}} \hookrightarrow (n+1)\kr\operatorname{CAT}
\]
to the unit object of $(\Spc\modd^n)\modd(\widehat{\Pr}{}^L)$, where the first map is the map that turns presentable modules into enriched categories. 
\end{notation}

\begin{remark}
The symmetric monoidal category underlying $n\Catscr^L$ is $\Spc\modd^n$. Given two objects $\Ccal, \Dcal$ in $\Spc\modd^n$, we have an equivalence
\[
\Hom_{n\Catscr^L}(\Ccal, \Dcal) = \psi_n \shom_{\Spc\modd^n}(\Ccal, \Dcal).
\]
Given a third object $\Ecal$ in $\Spc\modd^n$, the composition map
\[
\Hom_{n\Catscr^L}(\Ccal, \Dcal) \times \Hom_{n\Catscr^L}(\Dcal, \Ecal) \rightarrow \Hom_{n\Catscr^L}(\Ccal, \Ecal)
\]
is obtained by applying the lax symmetric monoidal functor $\psi_n$ to the morphism
\[
\shom_{\Spc\modd^n}(\Ccal, \Dcal) \otimes \shom_{\Spc\modd^n}(\Dcal, \Ecal) \rightarrow \shom_{\Spc\modd^n}(\Ccal, \Ecal)
\]
which is associated to the composite map
\[
\Ccal \otimes \shom_{\Spc\modd^n}(\Ccal, \Dcal) \otimes \shom_{\Spc\modd^n}(\Dcal, \Ecal)  \rightarrow \Dcal \otimes  \shom_{\Spc\modd^n}(\Dcal, \Ecal)  \rightarrow \Ecal.
\]
\end{remark}

\begin{definition}
We denote by $n\Prscr^L$ the full $(n+1)$-subcategory of $n \Catscr^L$ on the presentable $n$-categories. We call $n\Prscr^L$ the $(n+1)$-category of presentable $n$-categories. We let $n\Prscr_{\St}^L$ be the full $(n+1)$-subcategory of $n\Prscr^L$ on the presentable stable $n$-categories. We call $n\Prscr_{\St}^L$ the $(n+1)$-category of presentable stable $n$-categories
\end{definition}

\begin{remark}
The categories underlying $n\Prscr^L$ and $n\Prscr_{\St}^L$ are $n\kr\Pr^L$ and $n\kr\Pr^L_{\St}$, respectively. Since these are closed under tensor products inside $\Spc\modd^n$, we conclude that $n\Prscr^L$ and $n\Prscr_{\St}^L$ inherit symmetric monoidal structures from $n\Catscr^L$.
\end{remark}

\subsection{Conical colimits in presentable $n$-categories}\label{subsection conical en pres}

Our next goal is to show that the realization functor $\psi_n$ takes value in conically cocomplete $n$-categories. This will be accomplished by combining proposition \ref{prop existen conicals} with a generalization of proposition \ref{remark colimits vs conical}, which depends on a variant of proposition \ref{prop funct y changes}.

\begin{lemma}\label{lemma Fshriek preserva limits}
Let $F: \Mcal \rightarrow \Mcal'$ be a lax symmetric monoidal functor  between presentable symmetric monoidal categories. Assume that $F$ is accessible and preserves limits, and that the lax symmetric monoidal functor $\tau_{\Mcal'} F$ is equivalent to $\tau_{\Mcal}$. Then the induced functors
\[
F_!: \Algbrd(\Mcal) \rightarrow \Algbrd(\Mcal')
\]
and 
\[
F_!: \Cat^{\Mcal} \rightarrow \Cat^{\Mcal'}
\]
are accessible and preserve limits.
\end{lemma}
\begin{proof}
The fact that $\tau_{\Mcal'} F$ is equivalent to $\tau_{\Mcal}$ implies that the functor
\[
F_! : \Algbrd(\Mcal) \rightarrow \Algbrd(\Mcal')
\]
restricts to the full subcategories of enriched categories. Since the inclusion of enriched categories inside algebroids is creates limits and sufficiently filtered colimits, to prove the lemma it suffices to show that the above functor is accessible and limit preserving. 

Recall that $F_!$ is a morphism of cartesian fibrations over $\Cat$. For every category $X$, the functor
\[
(F_!)_X : \Algbrd_X(\Mcal) \rightarrow \Algbrd_X(\Mcal')
\]
is accessible and  limit preserving thanks to \cite{HA} corollaries  3.2.2.4 and 3.2.3.1. We know that $ \Algbrd_X(\Mcal)$ and $\Algbrd_X(\Mcal')$ are presentable by \cite{HA} corollary 3.2.3.5, and therefore by the adjoint functor theorem we see that $(F_!)_X$ admits a left adjoint. Using \cite{HA} proposition 7.3.2.6 we see that $F_!$ itself admits a left adjoint, and the lemma follows.
\end{proof}

\begin{lemma}\label{lemma fshriek funct}
Let $F: \Mcal \rightarrow \Mcal'$ be a lax symmetric monoidal functor  between presentable symmetric monoidal categories. Assume that $F$ is accessible and preserves limits, and that the lax symmetric monoidal functor $\tau_{\Mcal'} F$ is equivalent to $\tau_{\Mcal}$. Let $\Ical$ be a category and let $\Dcal$ be an $\Mcal$-enriched category. Then there is an equivalence of $\Mcal'$-enriched categories 
\[
 F_!\Funct(\Ical_\Mcal, \Dcal) = \Funct(\Ical_{\Mcal'}, F_!(\Dcal))
\]
which is natural in $\Ical$, and which at the level of objects enhances the canonical equivalence 
\[
\Hom_{\Cat^{\Mcal}}(\Ical_\Mcal, \Dcal) = \Hom_{\Cat}(\Ical, (\tau_{\Mcal})_! \Dcal) = \Hom_{\Cat^{\Mcal'}}(\Ical_{\Mcal'}, F_!\Dcal).
\]
\end{lemma}
\begin{proof}
Denote by $(F_!)^L$ the left adjoint to $F_!: \Cat^\Mcal \rightarrow \Cat^{\Mcal'}$, which is guaranteed to exist by lemma \ref{lemma Fshriek preserva limits}. Let $\Ecal$ be an $\Mcal'$-algebroid. We have
\begin{align*}
\Hom_{\Cat^{\Mcal'}}(\Ecal, \Funct(\Ical_{\Mcal'}, F_!(\Dcal))) &= \Hom_{\Cat^{\Mcal'}}(\Ecal \otimes \Ical_{\Mcal'}, F_!(\Dcal)) \\
&= \Hom_{\Cat^{\Mcal}}((F_!)^L (\Ecal \otimes \Ical_{\Mcal'}), \Dcal).
\end{align*}
Similarly, we have
\begin{align*}
\Hom_{\Cat^{\Mcal'}}(\Ecal, F_! \Funct(\Ical_{\Mcal}, \Dcal ) ) &= \Hom_{\Cat^{\Mcal}}((F_!)^L (\Ecal) , \Funct(\Ical_{\Mcal}, \Dcal )) \\
&= \Hom_{\Cat^{\Mcal}}((F_!)^L (\Ecal) \otimes \Ical_{\Mcal} , \Dcal).
\end{align*}
To construct the desired equivalence it suffices to construct a functorial equivalence between $(F_!)^L(\Ecal \otimes \Ical_{\Mcal'})$ and $(F_!)^L(\Ecal) \otimes \Ical_\Mcal$. Note that since $(\tau_{\Mcal'})_! F_!$ is equivalent to $(\tau_\Mcal)_!$, we can rewrite the latter as $(F_!)^L(\Ecal) \otimes (F_!)^L(\Ical_{\Mcal'})$. We now observe that there is a natural morphism
\[
\eta_{\Ecal, \Ical}: (F_!)^L(\Ecal \otimes \Ical_{\Mcal'}) \rightarrow (F_!)^L(\Ecal) \otimes (F_!)^L(\Ical_{\Mcal'})
\]
obtained by adjunction from the composite map
\[
\Ecal \otimes \Ical_{\Mcal'} \xrightarrow{} F_! (F_!)^L(\Ecal) \otimes F_!(F_!)^L(\Ical_{\Mcal'}) \rightarrow F_!((F_!)^L(\Ecal) \otimes (F_!)^L(\Ical_{\Mcal'}) )
\]
where the first arrow is induced from the unit of the adjunction $(F_!)^L \dashv F_!$, and the second map is induced from the lax symmetric monoidal structure on $F_!$. 

Since the assignment $(\Ecal, \Ical)\mapsto \eta_{\Ecal, \Ical}$ is colimit preserving, to show that $\eta_{\Ecal, \Ical}$ is an isomorphism it suffices to consider the case $\Ical = [1]$ and $\Ecal = \overline{C}_m$ is the enriched category induced from the cell $C_m$ for some $m$ in $\Mcal'$. Note that we have $(F_!)^L(\overline{C}_m) = \overline{C}_{F^L m}$ where $F^L$ denotes the left adjoint to $F$. Using proposition \ref{prop product cells} together with the fact that $(F_!)^L$ preserves colimits, we obtain a pushout square of $\Mcal$-algebroids
\[
\begin{tikzcd}
\overline{C}_{F^L m} \arrow{r}{} \arrow{d}{} & \overline{C}_{F^L m , 1_{\Mcal}} \arrow{d}{} \\
\overline{C}_{1_\Mcal, F^L m} \arrow{r} & (F_!)^L(\overline{C}_m \otimes \overline{C}_{1_{\Mcal'}}).
\end{tikzcd}
\]
The fact that $\eta_{\overline{C}_m, [1]}$ is an isomorphism follows now from another application of proposition \ref{prop product cells}.
\end{proof}

\begin{lemma}\label{lemma pree pushf of adjoints}
Let $F: \Mcal \rightarrow \Mcal'$ be a lax symmetric monoidal functor  between presentable symmetric monoidal categories. Assume that  the lax symmetric monoidal functor $\tau_{\Mcal'} F$ is equivalent to $\tau_{\Mcal}$. Let $\Ccal$ be an $\Mcal$-enriched category. Let $G: \Ccal \rightarrow \Dcal$ be a functor between $\Mcal$-enriched categories. Let $d$ be an object in $\Dcal$ and let $(c, \epsilon)$ be right adjoint to $G$ at $d$. Then $(c, \epsilon)$ is also right adjoint to $F_! G$ at $d$.
\end{lemma}
\begin{proof}
We have to show that for every object $e$ in $(\tau_{\Mcal'})_! F_! \Dcal =(\tau_\Mcal)_!\Dcal$ the composite map
\[
\Hom_{F_! \Ccal}(e, c) \xrightarrow{(F_! G)_*} \Hom_{F_!\Dcal}(G(e), G(c)) \xrightarrow{\epsilon} \Hom_{F_! \Dcal}(G(e), d)
\]
is an isomorphism. This is equivalent to the image under $F$ of the composite map
\[
\Hom_{\Ccal}(e, c) \xrightarrow{G_*} \Hom_{\Dcal}(G(e), G(c)) \xrightarrow{\epsilon} \Hom_{\Dcal}(G(e), d).
\]
Our claim now follows from the fact that the above composite map is an isomorphism since $(c, \epsilon)$ is right adjoint to $G$ at $d$.	
\end{proof}

\begin{lemma}\label{lemma pushf of adjoints}
Let $F: \Mcal \rightarrow \Mcal'$ be a lax symmetric monoidal functor  between presentable symmetric monoidal categories. Assume that $F$ is accessible and preserves limits, and that the lax symmetric monoidal functor $\tau_{\Mcal'} F$ is equivalent to $\tau_{\Mcal}$. Let $\Ccal$ be an $\Mcal$-enriched category. Let $\Ical$ be a category and $X^\lhd: \Ical^\lhd \rightarrow (\tau_\Mcal)_! \Ccal$ be a conical  limit diagram in $\Ccal$. Then $X^\lhd$ defines also a conical limit diagram in $F_! \Ccal$. 
\end{lemma}
\begin{proof}
Let $X = X^\lhd|_\Ical$. The  diagram $X^\lhd$ defines a pair $(c, \epsilon)$ right adjoint to $\Delta: \Ccal \rightarrow \Funct(\Ical_\Mcal, \Ccal)$ at $X$. By lemma \ref{lemma pree pushf of adjoints} we have that $(c, \epsilon)$ is also right adjoint to $F_! \Delta: F_! \Ccal \rightarrow F_! \Funct(\Ical_\Mcal, \Ccal)$ at $X$. Using lemma \ref{lemma fshriek funct} we may identify this with the diagonal functor $F_! \Ccal \rightarrow \Funct(\Ical_{\Mcal'}, \Ccal)$. This identification is the identity on objects, so we conclude that $X$ also has a conical limit in $F_! \Ccal$.  The lemma now follows from  corollary \ref{coro conical is ordinary colim}.
\end{proof}

\begin{lemma}\label{lemma psi preserves limits}
Let $n \geq 1$. Then the functor $\psi_n : \Spc\modd^n \rightarrow \reallywidehat{\nCat}$ is accessible and preserves large limits.
\end{lemma}
\begin{proof}
We argue by induction. The case $n = 1$ is a direct consequence of the fact that the forgetful functor $\widehat{\Cat}_{\cocompl} \rightarrow \widehat{\Cat}$ preserves limits and $\kappa_0$-filtered colimits. 

Assume now that $n > 1$. We continue with the notation from construction \ref{construction lax functor}. Note that we have an equivalence
\[
\Spc\modd^n = (n-1)\kr\Pr^L\modd = (\Spc\modd^{n-1})\modd(\widehat{\Pr}^L_{\kappa_0})
\]
given by ind-$\kappa_0$-completion. From this point of view, the functor $\varphi'_n$ agrees with the functor $\theta_{\Spc\modd^{n-1}}^{\kappa_0}$ from  notation \ref{notation thetaM kappa}. It follows from proposition \ref{prop have left adjoint} that the functor $\varphi'_n$ is accessible and preserves large limits. Using lemma \ref{lemma Fshriek preserva limits} together with our inductive hypothesis we conclude that $\psi_n'$ is accessible and preserves large limits. The lemma now follows from the fact that the inclusion of $\reallywidehat{\nCat}$ inside $\widehat{\Algbrd}(\reallywidehat{(n-1)\kr\Cat})$ creates large limits and filtered colimits.
\end{proof}

\begin{theorem}\label{teo psi tiene conicals}
Let $n \geq 1$ and let $\Ccal$ be an object in $\Spc\modd^n$. Then the $n$-category  $\psi_n(\Ccal)$ is admits all small conical colimits. Furthermore, the inclusion $\Ccal \rightarrow \psi_n(\Ccal)$ maps limits in $\Ccal$ to conical limits in $\psi_n(\Ccal)$.
\end{theorem}
\begin{proof}
The case $n = 1$ is clear, so assume $n > 1$. We continue with the notation from construction \ref{construction lax functor}. By proposition \ref{prop existen conicals}, the enriched category $\varphi_n(\Ccal)$ admits all  large conical limits and colimits. It follows that its full subcategory $\varphi'_n(\Ccal)$ admits all small conical colimits, and that any small limit in $\Ccal$ defines a conical limit in $\varphi'_n(\Ccal)$.  Our result now follows from a combination of lemmas \ref{lemma pushf of adjoints} and  \ref{lemma psi preserves limits}.
\end{proof}

\begin{corollary}
The $(n+1)$-category $n\Prscr^L$ is conically cocomplete for each $n \geq 1$.
\end{corollary}
\begin{proof}
Apply theorem \ref{teo psi tiene conicals} to $\Ccal = n\kr\Pr^L$ inside $\Spc\modd^{n+1}$.
\end{proof}

\subsection{The passage to adjoints property}\label{subsection passage}

A fundamental feature of presentable $1$-categories is that they are stable under many of the usual constructions of category theory. The following  would imply that the world of presentable $n$-categories enjoys similar closure properties:

\begin{conjecture}\label{conjecture limits}
The category $n\kr\Pr^L$ has all small limits for all $n \geq 0$.
\end{conjecture}

Note that for every $n \geq 1$ the inclusion 
\[
n\kr\Pr^L \rightarrow \Spc\modd^n(\widehat{\Cat}_\cocompl) = \Ind_{\kappa_0}(n\kr\Pr^L)
\]
preserves all small limits that exist in $n\kr\Pr^L$. Therefore conjecture \ref{conjecture limits} is equivalent to the claim that the $\kappa_0$-small objects in $\Spc\modd^n(\widehat{\Cat}_\cocompl)$ are closed under small limits.

 Our next goal is to prove a weak form of conjecture \ref{conjecture limits} (stated below as proposition \ref{coro limits prl}) which guarantees the existence of limits of diagrams of right adjoints.

\begin{definition}\label{definition adjointability}
Let $\Ecal$ be a commutative algebra in $\widehat{\Cat}_\cocompl$ and let $F: \Ccal \rightarrow \Dcal$ be a morphism in $\Ecal\modd$. We say that $F$ is left adjointable if its underlying functor admits a left adjoint $F^L: \Dcal \rightarrow \Ccal$, and the canonical structure of oplax morphism of $\Ecal$-modules on $F^L$ is strict. We say that $F$ is right adjointable if its underlying functor admits a colimit preserving right adjoint $F^R: \Ccal \rightarrow \Dcal$, and the canonical structure of lax morphism of $\Ecal$-modules on $F^R$ is strict.

Let $X : \Ical \rightarrow \Ecal\modd$ be a diagram. We say that $X$ is left (resp. right) adjointable if for every arrow $\alpha$ in $\Ical$ the morphism $X(\alpha)$ is left (resp. right) adjointable. In this case, the induced diagram $X^L: \Ical^\op \rightarrow \Ecal\modd$ (resp. $X^R: \Ical^\op \rightarrow \Ecal\modd$) is said to arise from $X$ by passage to left (resp. right) adjoints.
\end{definition}

\begin{proposition}\label{lemma adjoint}
Let $\Ecal$ be a commutative algebra in $\widehat{\Cat}_\cocompl$ which is generated under small colimits by its dualizable objects. Let $F: \Ccal \rightarrow \Dcal$ be a morphism in $\Ecal\modd$. Then
\begin{enumerate}[\normalfont (i)]
\item $F$ is left adjointable if and only if its underlying functor is left adjointable.
\item $F$ is right adjointable if and only if its underlying functor has a colimit preserving right adjoint.
\end{enumerate}
\end{proposition}
\begin{proof}
We give a proof of item (i) - the proof of (ii) is completely analogous. The structure of morphism of $\Ecal$-modules on $F$ induces a commutative square 
\[
\begin{tikzcd}
\Ecal \times \Ccal \arrow{r}{ \otimes } \arrow{d}{ \id_\Ecal \times F} & \Ccal \arrow{d}{F} \\
\Ecal \times \Dcal \arrow{r}{\otimes } & \Dcal
\end{tikzcd}
\]
which we have to show is vertically left adjointable. Since the horizontal arrows preserve colimits in the $\Ecal$-variable and $F^L$ preserves colimits, we may restrict to showing that for every dualizable object $e$ in $\Ecal$ the induced commutative diagram
\[
\begin{tikzcd}
\Ccal \arrow{r}{e \otimes - } \arrow{d}{F} & \Ccal \arrow{d}{F} \\
\Dcal \arrow{r}{e \otimes -} & \Dcal
\end{tikzcd}
\]
is vertically left adjointable. Note that the horizontal arrows have right adjoints given by $e^\vee \otimes -$, and the fact that $F$ is a morphism of $\Ecal$-modules implies that the square is in fact horizontally right adjointable. Since the vertical arrows admit left adjoints we conclude that the above square is also vertically left adjointable, as desired.
\end{proof}

\begin{notation}
Let $\kappa$ be an uncountable regular cardinal. Denote by $\Pr^L_\kappa$ the full subcategory of $\Pr^L$ on the $\kappa$-compactly generated categories and functors which preserve $\kappa$-compact objects. Denote by $\Cat^{\rex(\kappa)}$ the full subcategory of $\Cat$ on those categories admitting $\kappa$-small colimits, and functors which preserve $\kappa$-small colimits.  For each presentable category $\Ccal$ we denote by $\Ccal^\kappa$ the full subcategory of $\Ccal$ on the $\kappa$-compact objects. Recall from \cite{HTT} proposition 5.5.7.10 that passage to $\kappa$-compact objects and ind-$\kappa$-completion are inverse equivalences between $\Pr^L_\kappa$ and $\Cat^{\rex(\kappa)}$.
\end{notation}
\begin{lemma}\label{lemma PRLkappa}
Let $\kappa$ be an uncountable regular cardinal. Then the inclusion $\Pr^L_\kappa \rightarrow \Pr^L$ creates $\kappa$-small limits.
\end{lemma}
\begin{proof}
Let $X^\lhd: \Ical^\lhd \rightarrow \Pr^L$ be a $\kappa$-small limit diagram, and assume that $X = X^\lhd|_\Ical$ factors through $\Pr^L_\kappa$. Denote by $\ast$ the initial object of $\Ical^\lhd$. We have to show that  $X^\lhd(\ast)$ is $\kappa$-compactly generated, and that an object in $X^\lhd(\ast)$ is $\kappa$-compact if and only if its projection to $X(i)$ is $\kappa$-compact for every $i$ in $\Ical$. It suffices moreover to consider the case of $\kappa$-small products, and pullbacks.

We begin with the case of $\kappa$-small products, so that $\Ical$ is a $\kappa$-small set, and we have $X^\lhd(\ast) = \prod_{i \in \Ical} X(i)$. For each $i$ in $\Ical$ the projection $X^\lhd(\ast) \rightarrow X(i)$ has a right adjoint, which is induced by the identity map $X(i) \rightarrow X(i)$ together with the map $X(i) \rightarrow X(j)$ that picks out the final object of $X(j)$ for all $j \neq i$. This preserves colimits indexed by contractible categories, and in particular $\kappa$-filtered colimits, so we conclude that the projection $X^\lhd(\ast) \rightarrow X(i)$ preserves $\kappa$-compact objects. Combining this with \cite{HTT} lemma 5.3.4.10  we see that an object in $X^\lhd(\ast)$ is $\kappa$-compact if and only if its projection to $X(i)$ is $\kappa$-compact for all $i$ in $\Ical$.

 It remains to show that $X^\lhd(\ast)$ is $\kappa$-compactly generated. Observe that the projections $X^\lhd(\ast) \rightarrow X(i)$ are jointly conservative, admit left adjoints, and preserve $\kappa$-filtered colimits. It follows that $X^\lhd(\ast)$ is generated under colimits by the sequences $(c_i)_{i \in \Ical}$ where $c_i$ is a $\kappa$-compact object of $X(i)$ for all $i$ and $c_i$ is initial for all but one index $i$ - and moreover these objects are $\kappa$-compact in $X^\lhd(\ast)$. This shows that $X^\lhd(\ast)$ is $\kappa$-compactly generated.

We now deal with the case of pullbacks, so that $\Ical^\lhd = [1] \times [1]$, and $X$ corresponds to a pullback diagram which we depict as follows:
\[
\begin{tikzcd}
\Ccal' \arrow{d}{p'} \arrow{r}{q'} & \Ccal \arrow{d}{p} \\
\Dcal' \arrow{r}{q} & \Dcal
\end{tikzcd}
\]
 It follows from \cite{HTT} lemma 5.4.5.7 that if $c'$ is an object of $\Ccal'$ such that $q'c'$ and $p'c'$ are $\kappa$-compact then $c'$ is $\kappa$-compact. It therefore remains to show that $\Ccal'$ is generated under colimits by those $\kappa$-compact objects $c'$ such that $q'c'$ and $p'c'$ are $\kappa$-compact.  This is essentially a consequence of the proof of \cite{HTT} proposition 5.4.6.6. We repeat the relevant part of the argument below, suitably adapted to our context.
 
Denote by $\Ccal''$ the full subcategory of $\Ccal'$ on those objects whose projections to $\Ccal$ and $\Dcal'$ are $\kappa$-compact. Let $c'$ be an arbitrary object of $\Ccal'$, and let $c = q' c'$, $d' = p' c'$ and $d = pq'c'$ . We have a pullback diagram
 \[
 \begin{tikzcd}
 \Ccal''_{/c'} \arrow{r}{f'} \arrow{d}{g'} & \Ccal^\kappa_{/c} \arrow{d}{g} \\ (\Dcal')_{/d'}^\kappa \arrow{r}{f} & \Dcal^\kappa_{/d}.
 \end{tikzcd}
 \]
 It follows from \cite{HTT} lemma 5.4.6.1 (applied to the cardinals $\omega << \kappa$) that $f$ and $g$ are $\omega$-cofinal (see \cite{HTT} definition 5.4.5.8). Applying \cite{HTT} lemma 5.4.6.5 we see that $\Ccal''_{/c'}$ is $\kappa$-filtered that $f'$ and $g'$ are $\omega$-cofinal, and therefore from \cite{HTT} lemma 5.4.5.12 we have that $f'$ and $g'$ are cofinal. Denote by $c''$ the colimit of the natural map $\Ccal''_{/c'} \rightarrow \Ccal'$. Since $f'$ and $g'$ are cofinal and $\Ccal, \Dcal'$ are $\kappa$-compactly generated we see that the canonical map $c'' \rightarrow c'$ becomes an isomorphism upon composition with $q'$ and $p'$. Therefore $c'$ is a colimit of objects of $\Ccal''$, as desired.
 \end{proof}

\begin{lemma}\label{lemma limits kcompl}
Let $\kappa$ be an uncountable regular cardinal and let $X: \Ical \rightarrow \Cat^{\rex(\kappa)}$ be a $\kappa$-small diagram. Let $X^\rhd: \Ical^\rhd \rightarrow \Cat^{\rex(\kappa)}$ be a colimit diagram for $X$, and denote by $\ast$ the final object of $\Ical^\rhd$. Assume that for every arrow $\alpha$ in $\Ical$ the induced functor $X(\alpha)$ admits a right adjoint which preserves  $\kappa$-small colimits. Then 
\begin{enumerate}[\normalfont (i)]
\item For every $i$ in $\Ical$ the induced functor $X(i) \rightarrow X(\ast)$ admits a right adjoint which preserves $\kappa$-small colimits.
\item The induced diagram $(X^\rhd)^R : (\Ical^\op)^\lhd \rightarrow \Cat^{\rex(\kappa)}$ is a limit diagram.
\end{enumerate}
\end{lemma}
\begin{proof}
Consider the functor $
\Ind_\kappa X^\rhd : \Ical^\rhd \rightarrow \Pr^L$. This is a colimit diagram since $\Ind_\kappa$ preserves colimits. Note that we have a natural monomorphism $X \rightarrow X^\rhd$. The adjoint functor theorem guarantees that for every arrow $\alpha$ in $\Ical^\rhd$ the functor $\Ind_\kappa X^\rhd(\alpha)$ has a right adjoint. Moreover, results from \cite{HTT} section 5.5.3 guarantee that the induced diagram
\[
(\Ind_\kappa X^\rhd)^R: (\Ical^\op)^\lhd \rightarrow \widehat{\Cat}
\]
is a limit diagram. 

Let $\alpha$ be an arrow in $\Ical$. By virtue of being right adjoint to a functor of $\kappa$-compactly generated categories which preserves $\kappa$-compact objects we have that $(\Ind_\kappa X^\rhd)^R(\alpha)$ preserves $\kappa$-filtered colimits. Since $X(\alpha)$ is itself  right adjointable we see that $(\Ind_\kappa X^\rhd)^R(\alpha)$ preserves $\kappa$-compact objects. Moreover, the fact that the right adjoint to $X(\alpha)$ preserves $\kappa$-small colimits implies that $(\Ind_\kappa X^\rhd)^R(\alpha)$ is in fact colimit preserving.

 Since the forgetful functor $\Pr^L \rightarrow \widehat{\Cat}$ creates small limits, we conclude that $(\Ind_\kappa X^\rhd)^R$ factors through $\Pr^L$. Since $(\Ind_\kappa X^\rhd)^R|_{\Ical^\op}$ factors through $\Pr^L_\kappa$ and $\Ical$ is $\kappa$-small we conclude using lemma \ref{lemma PRLkappa} that $(\Ind_\kappa X^\rhd)^R$ factors through $\Pr^L_\kappa$, and in fact defines a limit diagram in $\Pr^L_\kappa$. In particular, we have that for every $i$ in $\Ical$ the right adjoint to the functor $\Ind_\kappa X(i) \rightarrow \Ind_\kappa X(\ast)$ preserves $\kappa$-compact objects and is colimit preserving, which establishes item (i). Item (ii) follows from the fact that passage to $\kappa$-compact objects provides an equivalence $\Pr^L_\kappa = \Cat^{\rex(\kappa)}$.
\end{proof}

\begin{corollary} \label{proposition passage to right}
Let $\Ecal$ be a commutative algebra in $\widehat{\Cat}_\cocompl$, generated under colimits by its dualizable objects. Let $X: \Ical \rightarrow  \Ecal\modd$ be a right adjointable diagram, with $\Ical$ small. Let $X^\rhd : \Ical^\rhd \rightarrow \Ecal\modd$ be a colimit diagram for $X$. Then $X^{\rhd}$ is right adjointable, and the induced diagram $(X^\rhd)^R : (\Ical^\op)^\lhd \rightarrow \Ecal\modd$ is a limit diagram.
\end{corollary}
\begin{proof}
This is a direct consequence of proposition \ref{lemma adjoint} and lemma \ref{lemma limits kcompl} applied to $\kappa =\kappa_0$ after enlarging the universe.
\end{proof}

We now specialize the above discussion to the case when $\Ecal$ is $(n-1)\kr\Pr^L$.

\begin{remark}
Let $n > 1$. Recall from remark \ref{remark generation} that $(n-1)\kr\Pr^L$ is generated under small colimits by the free modules of the form $(n-2)\kr\Pr^L \otimes \Pcal(\Ccal)$, where $\Ccal$ is a small category. Since $\Pcal(\Ccal)$ is dualizable in $\widehat{\Cat}_\cocompl$ we conclude that $(n-2)\kr\Pr^L \otimes \Pcal(\Ccal)$ is dualizable in $(n-1)\kr\Pr^L$. Hence $(n-1)\kr\Pr^L$ is generated under small colimits by its dualizable objects.
\end{remark}

\begin{proposition}\label{lemma adj iff psin adj}
Let $n \geq 1$ and let $F: \Ccal \rightarrow \Dcal$ be a morphism in $\Spc\modd^n$. Then the functor of $n$-categories $\psi_n(F)$ admits a left adjoint if and only if $F$ is left adjointable.
\end{proposition}
\begin{proof}
The case $n = 1$ is clear, so assume that $n > 1$. Recall from our construction that the functor of $1$-categories underlying  $\psi_n(F)$  is equivalent to $F$. It follows that if $\psi_n(F)$ is left adjointable then $F$ is left adjointable. 

Assume now that $F$ is left adjointable. We continue with the notation from construction \ref{construction lax functor}. Note that the functor of very large presentable categories $\Ind_{\kappa_0}F : \Ind_{\kappa_0}\Ccal \rightarrow \Ind_{\kappa_0}\Dcal$ is left adjointable. It follows from an application of proposition \ref{lemma adjoint} that $\Ind_{\kappa_0} F$ is a left adjointable morphism of $\Spc\modd^{n-1}$-modules.  By proposition \ref{prop adjuntos enriched}, we have that $\varphi_n(F)$ is a left adjointable functor of $\Spc\modd^{n-1}$-enriched categories. It follows from corollary \ref{corollary adjoints of enriched vs ordinary}, together with the fact that $\Ind_{\kappa_0}F$ and its left adjoint restrict to an adjunction on the subcategories of $\kappa_0$-compact objects,  that $\varphi_n'(F)$ is also left adjointable. The fact that $\psi_n(F)$ is left adjointable follows now from an application of lemma \ref{lemma pree pushf of adjoints}.
\end{proof}

\begin{proposition}\label{coro limits prl}
Let $n \geq 1$ and let $X : \Ical \rightarrow n\kr\Pr^L$ be a diagram with $\Ical$ small, such that for every arrow $\alpha$ in $\Ical$ the functor underlying $X(\alpha)$ admits a left adjoint. Then
\begin{enumerate}[\normalfont (i)]
\item The diagram $X$ is left adjointable, when thought of as a functor into $(n-1)\kr\Pr^L \modd$.
\item  Let $(X^L)^\rhd: (\Ical^\op)^\rhd \rightarrow n\kr\Pr^L$ be  a colimit diagram for $X^L$. Then $(X^L)^\rhd$ is right adjointable, and the induced diagram $((X^L)^\rhd)^R: \Ical^\lhd \rightarrow n\kr\Pr^L$ is a limit diagram for $X$.
\end{enumerate}
\end{proposition}
\begin{proof}
The case $n = 1$ follows from the results of \cite{HTT} section 5.5, so we assume that $n > 1$.  Item (i) then follows from an application of proposition \ref{lemma adjoint} and (ii) is a direct consequence of corollary \ref{proposition passage to right}.
\end{proof}

Our next goal is to recast proposition \ref{coro limits prl} in a language intrinsic to the $(n+1)$-categorical structure on $n\kr\Pr^L$. We will in fact be able to obtain a variant of it which works for any object in $\Spc\modd^n$.

\begin{notation}
Let $\Dcal$ be an $n$-category. We denote by $\Dcal^{\leq 1}$ the $1$-category underlying $\Dcal$, and by $(\Dcal^{\leq 1})^{\text{radj}}$ (resp. $\Dcal^{\leq 1})^{\text{ladj}}$ the subcategory of $\Dcal^{\leq 1}$ containing all objects, and only those morphisms which are right (resp. left) adjointable in $\Dcal$.
\end{notation}

\begin{definition}\label{def pass}
Let $\Dcal$ be an $n$-category. We say that $\Dcal$ satisfies the passage to adjoints property if the following conditions are satisfied:
\begin{itemize}
\item The category $(\Dcal^{\leq 1})^{\text{\normalfont radj}}$ has all colimits, and the inclusion $(\Dcal^{\leq 1})^{\text{\normalfont radj}} \rightarrow \Dcal$ preserves conical colimits.
\item The category $(\Dcal^{\leq 1})^{\text{\normalfont ladj}}$ has all limits, and the inclusion $(\Dcal^{\leq 1})^{\text{\normalfont ladj}} \rightarrow \Dcal$ preserves  conical limits.
\end{itemize}
\end{definition}

\begin{remark}
Let $\Dcal$ be an $n$-category. Then the passage to adjoints equivalence\footnote{The proof of our result does not rely on the existence of such an equivalence in the general case - however our choice of name for the passage to adjoints property is motivated by it.} from \cite{GR} induces an equivalence between the categories $(\Dcal^{\leq 1})^\text{radj}$ and $(\Dcal^{\leq 1})^{\text{ladj}}$. It follows that if $\Dcal$ satisfies the passage to adjoints property, then a right adjointable diagram $F: \Ical^\rhd \rightarrow \Dcal$ in $\Dcal$ is a conical colimit if and only if the diagram $F^R: (\Ical^\op)^\lhd \rightarrow \Dcal$ is a conical limit. 
\end{remark}

\begin{theorem}\label{teo passage}
Let $n \geq 2$ and let $\Ccal$ be an object in $\Spc\modd^n$. Then the $n$-category $\psi_n(\Ccal)$ satisfies the passage to adjoints property.
\end{theorem}

Our proof of theorem \ref{teo passage} requires a few lemmas.

\begin{lemma}\label{lemma criterio adjointability in psin}
Let $n \geq 2$ and let $\Ccal$ be an object in $\Spc\modd^n$. Let $\alpha: c \rightarrow d$ be a morphism in $\Ccal$. Denote by  $\shom_\Ccal$ the functor of Hom objects for $\Ccal$ as a module over $(n-1)\kr\Pr^L$. Then
\begin{enumerate}[\normalfont (i)]
\item The morphism $\alpha$ is left adjointable in $\psi_n(\Ccal)$ if and only if  for every morphism $\gamma: e \rightarrow e'$ in $\Ccal$, the commutative square of categories underlying the commutative square of $(n-2)\kr\Pr^L$-modules
\[
\begin{tikzcd}
\shom_\Ccal(e, c) \arrow{r}{\alpha_*} & \shom_\Ccal(e, d) \\ \shom_\Ccal(e', c) \arrow{u}{\gamma^*} \arrow{r}{\alpha_*} & \shom_\Ccal(e', d) \arrow{u}{\gamma^*}
\end{tikzcd}
\]
is horizontally left adjointable.
\item The morphism $\alpha$ is right adjointable in $\psi_n(\Ccal)$ if and only if for every morphism $\gamma: e \rightarrow e'$ in $\Ccal$, the commutative square of categories underlying the commutative square of $(n-2)\kr\Pr^L$-modules
\[
\begin{tikzcd}
\shom_\Ccal(d, e') \arrow{r}{\alpha^*} & \shom_\Ccal(c,e') \\ \shom_\Ccal(d, e) \arrow{u}{\gamma_*} \arrow{r}{\alpha^*} & \shom_\Ccal(c, e) \arrow{u}{\gamma_*}
\end{tikzcd}
\]
is horizontally left adjointable.
\end{enumerate}
\end{lemma}
\begin{proof}
This is a direct consequence of lemma \ref{lemma criterio adj} applied to the Yoneda embedding, as in the proof of proposition \ref{prop global adj iff local adjs}.
\end{proof}

\begin{lemma}\label{lemma limit of nat transf}
Let $n \geq 2$. Let $\Ical$ be a small category and let $X^\lhd, X'^\lhd: \Ical^\lhd \rightarrow \Spc\modd^{n-1}$ be two small limit diagrams. Denote by $\ast$ the initial object of $\Ical^\lhd$, and let $X = X^\lhd|_\Ical$ and $X' = X'^\lhd|_\Ical$. Let $\eta: X^\lhd \rightarrow X'^\lhd$ be a natural transformation, and assume that for every arrow $\alpha: i \rightarrow j$ in $\Ical$ the commutative square of categories
\[
\begin{tikzcd}
X(i) \arrow{r}{X(\alpha)} \arrow{d}{\eta_i} & X(j) \arrow{d}{\eta_j} \\
X'(i) \arrow{r}{X'(\alpha)} & X'(j)
\end{tikzcd}
\]
is horizontally left adjointable. Then
\begin{enumerate}[\normalfont (i)]
\item For every index $j$ in $\Ical$ the commutative square of categories
\[
\begin{tikzcd}
X^\lhd(\ast) \arrow{r}{} \arrow{d}{\eta_\ast} & X(j) \arrow{d}{\eta_j} \\
X'^\lhd(\ast) \arrow{r}{} & X'(j)
\end{tikzcd}
\]
is horizontally left adjointable.
\item Assume given an extension of $\eta|_{\Ical}$ to a natural transformation $\mu$ between functors $Y^\lhd, Y'^\lhd: \Ical^\lhd \rightarrow \Spc\modd^{n-1}$, and that for every index $j$ in $\Ical$ the commutative square of categories
\[
\begin{tikzcd}
Y^\lhd(\ast) \arrow{r}{} \arrow{d}{\mu_\ast} & X(j) \arrow{d}{\eta_j} \\
Y'^\lhd(\ast) \arrow{r}{} & X'(j)
\end{tikzcd}
\]
is horizontally left adjointable. Then the induced commutative square of categories
\[
\begin{tikzcd}
Y^\lhd(\ast) \arrow{r}{} \arrow{d}{\mu_\ast} & X^\lhd(\ast) \arrow{d}{\eta_\ast} \\
Y'^\lhd(\ast) \arrow{r}{} & X'^\lhd(\ast)
\end{tikzcd}
\]
\end{enumerate}
is horizontally left adjointable.
\end{lemma}
\begin{proof}
We first establish item (i). It follows from corollary \ref{proposition passage to right} that for every $j$ in $\Ical$ the functors $X^{\lhd}(\ast) \rightarrow X(j)$ and $X'^\lhd(\ast) \rightarrow X'(j)$ are left adjointable. Consider now the commutative square of categories
\[
\begin{tikzcd}
\Ind_{\kappa_0}X^\lhd(\ast) \arrow{r}{} \arrow{d}{\Ind_{\kappa_0}\eta_\ast} & \Ind_{\kappa_0}X(j) \arrow{d}{\Ind_{\kappa_0} \eta_j} \\
\Ind_{\kappa_0}X'^\lhd(\ast) \arrow{r}{} & \Ind_{\kappa_0}X'(j).
\end{tikzcd}
\]
Our claim will follow if we are able to show that this is horizontally left adjointable. Since the horizontal arrows admit left adjoints, it suffices in fact to show that this is vertically right adjointable. By virtue of lemma \ref{lemma criterio adj}, it suffices to show that the diagram
\[
G^\lhd: \Ical^\lhd \rightarrow \Funct([1], \operatorname{CAT}^{2\text{-cat}})
\]
induced from the natural transformation 
\[
\Ind_{\kappa_0} \eta : \Ind_{\kappa_0}X^\lhd \rightarrow \Ind_{\kappa_0} X'^\lhd
\]
factors through $\Funct(\Adj, \operatorname{CAT}^{2\text{-cat}})$, where $\operatorname{CAT}^{2\text{-cat}}$ denotes the $2$-category of very large categories. Combining propositions \ref{prop existen conicals} and \ref{prop limits in functor cat} with lemma \ref{lemma PRLkappa} we see that $G^\lhd$ is in fact a conical limit diagram.

For each arrow $\alpha: i \rightarrow j$ in $\Ical$, the commutative square  of categories
\[
\begin{tikzcd}[column sep = huge, row sep = large]
\Ind_{\kappa_0} X(i) \arrow{r}{\Ind_{\kappa_0} X(\alpha)} \arrow{d}{\Ind_{\kappa_0} \eta_i} & \Ind_{\kappa_0} X(j) \arrow{d}{\Ind_{\kappa_0} \eta_j} \\
\Ind_{\kappa_0} X'(i) \arrow{r}{\Ind_{\kappa_0} X'(\alpha)} & \Ind_{\kappa_0} X'(j)
\end{tikzcd}
\]
is horizontally left adjointable thanks to our hypothesis on $\eta$. By the adjoint functor theorem, we also know that the vertical arrows are right adjointable. It follows that the above square is also vertically right adjointable. Another application of lemma \ref{lemma criterio adj} shows that $G^\lhd|_{\Ical}$ factors through $\Funct(\Adj, \operatorname{CAT}^{2\text{-cat}})$. Applying proposition \ref{prop limits in functor cat} we conclude that $G^\lhd$ factors through 
$\Funct(\Adj, \operatorname{CAT}^{2\text{-cat}})$, as we claimed.

We now prove item (ii). Applying corollary \ref{proposition passage to right} we see that $(X^\lhd)^L$ is a colimit diagram in $\Spc\modd^{n-1}$ and hence we have a morphism $(X^\lhd)^L \rightarrow (Y^\lhd)^L$ which restricts to the identity on $\Ical^\op$. Consider now the induced morphism 
\[
 \Ind_{\kappa_0}(X^\lhd)^L \rightarrow  \Ind_{\kappa_0}(Y^\lhd)^L.
\] 
Passing to right adjoints, we obtain a natural transformation $ \Ind_{\kappa_0}Y^\lhd \rightarrow  \Ind_{\kappa_0}X^\lhd$  which is the identity on $\Ical$. Since $ \Ind_{\kappa_0}X^\lhd$ is a limit diagram by lemma \ref{lemma PRLkappa}, we conclude that the induced functor 
\[
\Ind_{\kappa_0}Y^\lhd(\ast) \rightarrow  \Ind_{\kappa_0} X^\lhd(\ast)
\]
is obtained by ind-$\kappa_0$-completion of the map $Y^\lhd(\ast) \rightarrow X^\lhd(\ast)$ induced from $\mu$. It follows that the functor $Y^\lhd(\ast) \rightarrow X^\lhd(\ast)$ is left adjointable. Similarly, we have that the functor $Y'^\lhd(\ast) \rightarrow X'^\lhd(\ast)$ is horizontally left adjointable. 

As before, our claim would follow if we are able to show that the commutative square of categories
\[
\begin{tikzcd}
\Ind_{\kappa_0} Y^\lhd(\ast) \arrow{r}{} \arrow{d}{\Ind_{\kappa_0} \mu_\ast} &  \Ind_{\kappa_0} X^\lhd(\ast) \arrow{d}{\Ind_{\kappa_0} \eta_\ast} \\
\Ind_{\kappa_0} Y'^\lhd(\ast) \arrow{r}{} &\Ind_{\kappa_0}  X'^\lhd(\ast)
\end{tikzcd}
\]
is horizontally left adjointable. Since the horizontal arrows admit left adjoints, it suffices to show that the square is in fact vertically right adjointable. Consider the diagram
\[
H^\lhd: \Ical^\lhd \rightarrow \Funct([1], \operatorname{CAT}^{2\text{-cat}})
\]
induced from the natural transformation 
\[
\Ind_{\kappa_0} \mu : \Ind_{\kappa_0}Y^\lhd \rightarrow \Ind_{\kappa_0} Y'^\lhd.
\]
This extends the functor $G^\lhd|_\Ical$. By lemma  \ref{lemma criterio adj}, we may reduce to showing that the induced natural transformation $H^\lhd \rightarrow G^\lhd$ factors through $\Funct(\Adj, \operatorname{CAT}^{2\text{-cat}})$. This is a consequence of the fact that $G^\lhd$ is a conical limit diagram in $\Funct(\Adj, \operatorname{CAT}^{2\text{-cat}})$.
\end{proof}

\begin{proof}[Proof of theorem \ref{teo passage}]
To simplify notation, we set $\Ccal^{\text{radj}} = (\psi_n(\Ccal)^{\leq 1})^{\text{radj}}$ and $\Ccal^{\text{ladj}} = (\psi_n(\Ccal)^{\leq 1})^{\text{ladj}}$.

Let $X: \Ical \rightarrow \Ccal^{\text{radj}}$ be a diagram, and let $X^\rhd: \Ical^\rhd \rightarrow \Ccal$ be an extension of $X$ to a colimit diagram in $\Ccal$. We claim that $X^\rhd$ factors through $\Ccal^{\text{radj}}$. Let $\gamma: e \rightarrow e'$ be a morphism in $\Ccal$. We have limit diagrams
\[
\shom_\Ccal(X^\rhd, e) : (\Ical^\rhd)^\op \rightarrow \Spc\modd^{n-1}
\]
and 
\[
\shom_\Ccal(X^\rhd, e') : (\Ical^\rhd)^\op \rightarrow \Spc\modd^{n-1}.
\]
The morphism $\gamma$ induces a natural transformation
\[
\shom_\Ccal(X^\rhd, \gamma): \shom_\Ccal(X^\rhd, e) \rightarrow \shom_\Ccal(X^\rhd, e').
\]

Combining lemma \ref{lemma criterio adjointability in psin} and part (i) of lemma \ref{lemma limit of nat transf} we see that for every $i$ in $\Ical$ the commutative square of categories
\[
\begin{tikzcd}
\shom_\Ccal(X^{\rhd}(\ast), e) \arrow{d}{\gamma_*} \arrow{r}{} & \shom_\Ccal(X(i), e) \arrow{d}{\gamma_*}\\
\shom_\Ccal(X^{\rhd}(\ast), e') \arrow{r}{} & \shom_\Ccal(X(i), e')
\end{tikzcd}
\]
is horizontally left adjointable. Using lemma \ref{lemma criterio adjointability in psin} once more  we conclude that $X^\rhd$ factors through $\Ccal^\radj$. 

Similarly, combining lemma \ref{lemma criterio adjointability in psin} with part (ii) of lemma \ref{lemma limit of nat transf} we see that if $Y^\rhd$ is another extension of $X$ which factors through $\Ccal^\radj$, then the induced morphism $X^\rhd(\ast) \rightarrow Y^\rhd(\ast)$ belongs to $\Ccal^\radj$. This means that $X^\rhd$ is in fact also a colimit diagram in $\Ccal^\radj$, and therefore we see that the inclusion $\Ccal^\radj \rightarrow \Ccal$ creates colimits. It now follows from theorem \ref{teo psi tiene conicals} that $\psi_n(\Ccal)$ satisfies the first condition in definition \ref{def pass}.

We now show that $\psi_n(\Ccal)$ satisfies the second condition in definition \ref{def pass}. Thanks to theorem \ref{teo psi tiene conicals}, it suffices to show that $\Ccal^\ladj$ has arbitrary products and fiber products, and that these are preserved by its inclusion into $\Ccal$. An argument analogous to the case of colimits reduces one to showing that $\Ccal$ has arbitrary products and fiber products of diagrams in $\Ccal^\ladj$.

We consider the case of pullbacks - the proof for arbitrary products is analogous. Let $\Ical$ be the category with objects $0, 1, 2$ and nontrivial arrows $0 \rightarrow 2 \leftarrow 1$. Let $X : \Ical \rightarrow \Ccal^\ladj$ be a diagram. Extend $X$ to a limit diagram $X^\lhd : \Ical^\lhd \rightarrow \Ind_{\kappa_0} \Ccal$. Repeating our argument for the case of colimits shows that $X^\lhd$ defines in fact a limit diagram in the wide subcategory $(\Ind_{\kappa_0} \Ccal)^\ladj$ of $\Ind_{\kappa_0} \Ccal$  containing those arrows  which are left adjointable when considered inside the $n$-category $(\psi_{n-1})_! \varphi_n(\Ccal)$. 

Consider now the diagram $X^L : \Ical^\op \rightarrow \Ccal^\radj$ obtained from $X$ by passage to left adjoints. As in the case of $X$, we may extend $X^L$ to a diagram $(X^L)^\rhd: (\Ical^\op)^\rhd \rightarrow  (\Ind_{\kappa_0} \Ccal)^\radj$ which is both a colimit diagram  in $(\Ind_{\kappa_0} \Ccal)^\radj$ and $\Ind_{\kappa_0} \Ccal$. Since pushout squares in $(\Ind_{\kappa_0} \Ccal)^\radj$ and pullback squares in $(\Ind_{\kappa_0} \Ccal)^\ladj$ are in one to one correspondence by passing to adjoints of all arrows involved, we obtain an equivalence between $(X^L)^\rhd(\ast)$ and $X^\lhd(\ast)$. Since $\Ccal$ is closed under small colimits inside $\Ind_{\kappa_0} \Ccal$ we conclude that $X^\lhd(\ast)$ in fact belongs to $\Ccal$. Hence $X^\lhd$ gives the desired extension of $X$ to a limit diagram in $\Ccal$.
\end{proof}


\bibliographystyle{amsalpha}
\bibliography{References}
 
\end{document}